\documentclass[a4paper]{report}
\usepackage{amssymb,amsmath,amsthm,amsfonts,mathrsfs,amscd,verbatim}
\usepackage{setspace}
\usepackage[dvips]{graphicx,color}
\usepackage{epsfig}

\newtheorem{thm}{Theorem}[section]
\newtheorem{prop}[thm]{Proposition}

\newtheorem{lem}[thm]{Lemma}

\theoremstyle{definition}
\newtheorem{defn}[thm]{Definition}

\theoremstyle{remark}
\newtheorem{rem}[thm]{Remark}
\newtheorem{ex}[thm]{Example}

\title{Extremal metrics and K-stability}
\date{}

\begin{document}

\begin{center}
  \vspace*{1in}
  {\LARGE\bf Extremal metrics and K-stability}
  \par
  \vspace{1.5in}
  {\large G\'abor Sz\'ekelyhidi}
  \par
  \vspace{0.1in}
  {\large Imperial College, University of London}
  \par
  \vfill
  A thesis presented for the degree of Doctor of Philosophy of the
  University of London and the Diploma of Membership of Imperial College
\end{center}
\pagebreak

\begin{center}
  \vspace*{1in}
  {\LARGE\bf Declaration}
  \vspace*{0.5in}
\end{center}

{The material presented in this thesis is the author's own,
except where it appears with attribution to others. }

\pagebreak
\vspace*{1.6in}
\begin{center}{\large\bf Abstract}
\end{center}

In this thesis we study the relationship between the existence of
canonical metrics on a complex manifold and stability in the sense of
geometric invariant theory. We introduce a modification of K-stability
of a polarised variety
which we conjecture to be equivalent to the existence of an extremal
metric in the polarisation class. 
A variant for a complete extremal metric on the complement of a
smooth divisor is also given. On toric surfaces we prove a 
Jordan-H\"older type theorem for decomposing
semistable surfaces into stable pieces. 
On a ruled surface we compute the infimum of the
Calabi functional for the unstable polarisations, exhibiting a
decomposition analogous to the Harder-Narasimhan filtration of an
unstable vector bundle. 

\vfill
\newpage

\begin{spacing}{1.1}
\tableofcontents
\end{spacing}
\vfill\newpage 

\chapter*{Introduction}
\addcontentsline{toc}{chapter}{Introduction}

The subject of this thesis is finding canonical metrics on K\"ahler
manifolds. The first result of this form is the classical uniformisation
theorem, which states that every compact Riemann surface admits a metric
of constant curvature, unique if we prescribe the total area. In higher
dimensions a condition analogous to prescribing the total area is fixing
the K\"ahler class of a metric. 
In~\cite{Cal82} Calabi introduced the functional 
\[ \int_M S(\omega)^2\frac{\omega^n}{n!}, \]
for K\"ahler metrics $\omega$ on $M$ in a fixed cohomology class, where
$S(\omega)$ is the scalar curvature of $\omega$. He proposed finding
critical points of this functional (the Calabi functional) as candidates
for a canonical metric in the K\"ahler class. 
Such a metric is called an
\emph{extremal metric} and the main problem is their 
uniqueness and existence.
It has been shown in~\cite{CT05} that any two extremal metrics in a
K\"ahler class are related by a holomorphic automorphism. The question
of existence is still open. 

Calabi showed that the Euler-Lagrange
equation of the variational problem 
is that the gradient of the scalar curvature is a holomorphic
vector field. A special case is that of K\"ahler-Einstein metrics since
these have constant scalar curvature. In this case the first Chern class
of the manifold is proportional to the K\"ahler class of the metric, and so it
has to be negative definite, zero, or positive definite. 
In the case when the first Chern class is
negative or zero, Yau~\cite{Yau78} (also Aubin~\cite{Aub78} in the
negative case) showed
that the variety
admits a unique K\"ahler-Einstein metric, solving a conjecture of Calabi. 
The case of positive first Chern class proved to be more difficult and is
still not completely resolved. Yau conjectured that in this case the
existence of a K\"ahler-Einstein metric is related to the stability of
the underlying variety in the sense of Mumford's geometric invariant
theory~\cite{MFK94}. 
Tian made great progress towards understanding this (see~\cite{Tian97})
giving an analytic ``stability'' condition which is equivalent to the 
existence of a
K\"ahler-Einstein metric. This condition is the properness of the
Mabuchi functional, which is an energy functional on the K\"ahler class
whose critical points are K\"ahler-Einstein metrics. 
In~\cite{Tian97} 
Tian also defined the algebro-geometric notion of K-stability (not
exactly the same
as what we call K-stability), which is
satisfied when the Mabuchi functional is proper. 

In~\cite{Don97}, Donaldson showed that the scalar curvature arises as a
moment map for a suitable infinite dimensional symplectic action (see
also Fujiki~\cite{Fuj92}). This put earlier results into a new context,
and explained on a formal level why the existence of a
K\"ahler-Einstein metric, or more generally a metric of constant scalar
curvature (cscK), is related to the stability of the variety. 
Moreover it made
it possible to formulate precise conjectures. In particular
in~\cite{Don02} Donaldson generalised Tian's definition of K-stability
by giving an algebro-geometric definition of the Futaki invariant, and
conjectured that it is equivalent to the existence of a cscK metric. 

The definition of K-stability is roughly the following (see
Section~\ref{sec:Kstab} for details). Given a polarised variety (this
means we have chosen a K\"ahler class which is the first Chern class of
an ample line bundle), we
consider degenerations of it into possibly singular schemes. For each
such test-configuration we define a number called
the generalised Futaki invariant, and the variety is K-stable if this
number is positive for all non-trivial test-configurations. The idea is
that the Futaki invariant controls the asymptotic behaviour of the
Mabuchi functional as we tend to a degenerate metric, so that properness
of the Mabuchi functional corresponds to the Futaki invariant being
positive for all nontrivial degenerations. This means that an important
problem is to study the metric behaviour of such an algebro-geometric
test-configuration. In~\cite{Tian97} Tian studied the case where the
central fibre is normal. In Section~\ref{sec:metricdegen} we will study the
case of deformation to the normal cone of the zero section of a ruled
manifold, so the central fibre has a normal crossing singularity. 

An interesting testing ground for these ideas is the case of toric
varieties, which was developed by Donaldson. In~\cite{Don02} he showed
that in the case of toric surfaces K-stability
implies that the Mabuchi functional is bounded from below (this was later
extended by Zhou-Zhu~\cite{ZZ06} to show properness of the Mabuchi
functional) and a minimising sequence has a subsequence that converges
in a weak sense. In~\cite{Don05_1} Donaldson proved interior estimates for
the cscK equation for toric surfaces. 

Unfortunately recent examples in~\cite{ACGT3} show that positivity of
the Futaki invariant for algebraic test-configurations may not be
enough to ensure the existence of a cscK metric. One approach suggested
in~\cite{RT06} is to allow more general test-configurations with
polarisations which are real linear 
combinations of line bundles, or with
non-algebraic central fibres. 
In Section~\ref{sec:uniformKstab} we suggest an alternative 
way of strengthening the
definition of K-stability to what we call uniform K-stability, and
then in
Section~\ref{sec:toricsurface} we show that a K-polystable toric surface is
uniformly K-polystable.

So far we have only considered cscK metrics, and 
it is natural to ask whether one can give a stability 
criterion for the existence
of general extremal metrics. Such a criterion was proposed by the author
in~\cite{GSz04} (see also Mabuchi~\cite{Mab04_1} for a different
definition). Given the interpretation of the scalar
curvature as a moment map, 
what one needs to do is to find a stability
criterion satisfied by the orbit of a critical point of the norm squared
of the moment map in general. The norm squared of the moment map was
studied by Kirwan~\cite{Kir84}, but not exactly from this point of view,
so we develop the finite dimensional theory in
Chapter~\ref{chap:finiteGIT}.

One advantage of extending the search for canonical metrics from cscK to general
extremal metrics is that we have interesting explicit examples such as the ruled
surfaces constructed in~\cite{TF97} and the more general constructions
in~\cite{ACGT3}. In Section~\ref{sec:extremalruled} we will use such
explicit constructions to give complete extremal metrics on a ruled
surface. These will then be used in Section~\ref{sec:extremallwr} to
determine the infimum of the Calabi functional for the unstable
polarisations. The infimum is achieved
by a degenerate metric, where the variety splits up into pieces which
either admit a complete extremal metric, or collapse. 
In general if a variety is unstable one
expects that there is such a decomposition into stable pieces in analogy
with the Harder-Narasimhan filtration of an unstable vector bundle. 

\subsection*{Chapter summary}
In Chapter~\ref{chap:finiteGIT} we develop the finite
dimensional theory of stability and the moment map. Apart from
Section~\ref{sec:modulus}, the material in this chapter is contained in
slightly different form in the work of Kirwan~\cite{Kir84}. The main
result is an extension of the Kempf-Ness theorem. 
\theoremstyle{plain}
\newtheorem*{thmstab}{Theorem~\ref{thm:stab}}
\begin{thmstab} 
  A point $x$ in $X$ is in the $G$-orbit of a critical point of
  $\Vert\mu\Vert^2$, if and only if it is polystable relative to a
  maximal torus 
  which fixes it. 
\end{thmstab}
\noindent Here $G$ is a reductive group acting on a complex variety $X$, and $\mu$
is the moment map for the action of a maximal
compact subgroup of $G$. 
In Section~\ref{sec:modulus} we 
introduce the notion of the modulus of stability of a stable point and
show that one can give a lower bound for the first eigenvalue of the
derivative of the moment map in terms of this modulus
(see Theorem~\ref{thm:eigenval}). In the final section we show how the
theory works out in the case of a torus action. 

Chapter~\ref{chap:extremal} is a review of some well-known results 
about extremal
metrics. In Section~\ref{sec:scalmoment} we recall that the scalar
curvature arises as a moment map for an infinite dimensional group
action. Together with the results in the first chapter, this gives the
motivation for the subsequent results, in particular the definition of
K-stability in the next chapter.

In Chapter~\ref{chap:stabvar} we recall the definition of K-polystability
from~\cite{Don02} and we introduce the notion of relative
K-polystability which is the suitable generalisation to the case of 
extremal metrics
with non-constant scalar curvature. In Section~\ref{sec:uniformKstab} we
introduce the notion of uniform K-polystability. This addresses the
problem mentioned above that K-polystability may not be enough to ensure
the existence of a cscK metric, but it is still to be seen whether
uniform K-polystability is the correct notion. We then consider the case
of a pair $(X,D)$ where $D\subset X$ is a divisor and define a variant of
K-polystability for this situation. The aim is to find a condition for 
$X\setminus D$ to admit a complete extremal metric, ie. a complete
metric such that the gradient of its scalar curvature is a holomorphic
vector field, which is asymptotically hyperbolic near $D$ 
(see Section~\ref{sec:kstabpair} for the definition
of the class of metrics we consider). 
In Section~\ref{sec:ruledexample} we illustrate the
definition of K-polystability on a ruled surface by finding
destabilising test-configurations for certain polarisations. This will be
complemented in Section~\ref{sec:extremalruled} where we construct
extremal metrics on this ruled surface (both compact metrics and
complete metrics on the complement of a section) for the other
polarisations. 

In Section~\ref{sec:lwrcalabi} we recall Donaldson's
theorem in~\cite{Don05} which gives a lower bound for the Calabi functional
in terms of a destabilising
test-configuration. This gives a fairly simple proof of the fact that a
variety that admits a cscK metric must be K-semistable. We give the
following refinement of the theorem. 
\newtheorem*{thmextremallwr}{Theorem~\ref{thm:extremallwr}}
\begin{thmextremallwr}
Let $T$ be a maximal torus of automorphisms of a polarised variety
$(X,L)$ with
corresponding extremal vector field $\chi$. Suppose there is a
test-configuration for $(X,L)$ 
compatible with $T$ such that the modified Futaki
invariant $F_\chi(\alpha)<0$ for
the $\mathbf{C}^*$-action $\alpha$ induced on the central fibre.
Then for any metric $\omega\in 2\pi c_1(L)$,
\[ \Vert S(\omega)-\hat{S}\Vert^2_{L^2}\geq 2\cdot (2\pi)^n
\frac{F_\chi(\alpha)^2}{
\Vert\alpha\Vert^2} + \Vert\chi\Vert_{L^2}^2. \]
\end{thmextremallwr}
\noindent Here $\hat{S}$ is the average scalar curvature. 
This theorem shows that a polarised variety that admits an extremal
metric is relatively K-semistable since if $\omega$ is an extremal
metric then $\Vert S(\omega)-\hat{S}\Vert^2_{L^2} = \Vert\chi
\Vert_{L^2}^2$ (see Section~\ref{sec:FMprelim}).

In Chapter~\ref{chap:toric} we study toric varieties. First we
generalise the toric test-configurations defined in~\cite{Don02} to
bundles of toric varieties and compute their Futaki invariants
(Theorem~\ref{thm:toricbtc}). We will use this in the next chapter to
define test-configurations for a ruled manifold. In
Section~\ref{sec:toricsurface} we concentrate on toric surfaces and
prove two results. The first is that a K-polystable toric surface is
uniformly K-polystable, which relies on
\newtheorem*{thmineq1}{Proposition~\ref{prop:ineq1}}
\begin{thmineq1}
  Given a convex polygon $P$ there exists a constant $C$ such 
  that for all non-negative
  continuous convex functions $f$ on $P$, 
  \[ \Vert f\Vert_{L^2(P)}\leq C\int_{\partial P} f\, d\sigma.\]
\end{thmineq1}
We then use the notion of measure majorisation from convex geometry to
study semistable surfaces, and prove
\newtheorem*{thmsemidecomp}{Theorem~\ref{thm:semidecomp}}
\begin{thmsemidecomp}
  A K-semistable polygon $P$ has a canonical decomposition into
  subpolygons $Q_i$ each of which is either K-polystable, or a
  parallelogram with two opposite edges lying on edges of $P$. 
\end{thmsemidecomp}
\noindent A subpolygon $Q_i$ defines a pair $(X_i,D_i)$, with the divisor
corresponding to the edges of $Q_i$ lying in the interior of $P$.
K-polystability of $Q_i$ is interpreted as K-polystability of the pair
$(X_i,D_i)$. 

In the final chapter we study ruled manifolds using the explicit
construction of metrics from momentum profiles due to
Hwang-Singer~\cite{HS02}. In Section~\ref{sec:metricdegen} we construct
a sequence of metrics which model the deformation to the normal cone of
a section, and show that the derivative of the Mabuchi functional along
this degeneration tends to the Futaki invariant of the corresponding
test-configuration. We then restrict attention to a ruled surface. First
we use momentum profiles to construct extremal metrics on it, as well as
complete extremal metrics on the complement of the zero or infinity
section. We see that we obtain the same restrictions on the polarisation
as in the stability computation in Section~\ref{sec:ruledexample}. In
the last section we show that the infimum of the Calabi functional for
the unstable polarisations is achieved by degenerate metrics assembled
from the complete extremal metrics we have constructed. 

\subsection*{Acknowledgements}
I would like to thank my supervisor Simon Donaldson for his generosity
in sharing his ideas and his patience when explaining them. I would also
like to thank the members of the geometry group at Imperial for creating
a great atmosphere for doing mathematics. I am grateful to my PhD
examiners Michael Singer and Richard Thomas for their 
many helpful comments. 

I acknowledge funding received
from the Overseas Research Council, the Department of Mathematics at
Imperial and EPSRC during the past three years. 

\chapter{Finite dimensional GIT}\label{chap:finiteGIT}

This chapter contains some background on the finite dimensional
theory of geometric invariant theory and symplectic quotients. 
The basic references are
Mumford-Fogarty-Kirwan~\cite{MFK94} and Kirwan~\cite{Kir84} (see also
Thomas~\cite{Thomas06}). Essentially the only novelty is in
Section~\ref{sec:modulus}, the results in the rest of the chapter can
be obtained from the theory of Kirwan~\cite{Kir84}.

The aim of geometric invariant theory (GIT) is to define a quotient
variety $X/G$ when an algebraic group $G$ acts on an algebraic variety
$X$. It is
natural to require functions over $X/G$ to be given by $G$ invariant
functions over $X$, and this requirement gives a simple definition for
the quotient. The difficulty is to understand what the quotient variety
parametrises. In other words, we would like to understand the projection
map from $X$ to $X/G$. There will be certain bad (unstable) orbits where
this map is not defined, and also some semistable orbits which become
identified with each other. 
This will be discussed in Section~\ref{sec:stability}.

In symplectic geometry there is also a way of constructing quotients.
Here we start with a symplectic manifold $M$ with symplectic form
$\omega$, and a compact group $K$ acting on $M$, preserving $\omega$. In the
case where $M$ is also an algebraic variety, then the Kempf-Ness theorem,
discussed in Section~\ref{sec:kempfness}
relates the symplectic quotient by $K$ to the GIT quotient by the
complexification of $K$. 

A central role is played by the norm squared of the moment map, which in
the infinite dimensional setting is the Calabi functional. In
Section~\ref{sec:normsquared} we show how
one can characterise orbits of critical points of this functional using
stability,
generalising the Kempf-Ness theorem.  

In Section~\ref{sec:modulus} we introduce a notion we call the modulus
of stability which measures how far a point is from being unstable, and
we prove some simple results about it. This is used in
Chapter~\ref{chap:stabvar} to motivate the definition of uniform
K-stability. 

In Section~\ref{sec:torus} we illustrate the above theory in the case of
a torus action, where everything can be seen quite explicitly. While a
torus action may seem very special, the Cartan decomposition implies
that many questions about general actions can be reduced to a torus
action.

\section{Stability}\label{sec:stability}

To give precise definitions let $(X,L)$ be a smooth
complex projective variety with an
ample line bundle, in other words a \emph{polarised variety}. The graded
ring of functions over $X$ is defined to be 
\[ R(X) = \bigoplus_{k=0}^\infty H^0(X,L^k). \]
Suppose a complex reductive group $G$ acts on $X$ by holomorphic
automorphisms. Suppose we can lift this action of $G$ to a holomorphic
action on $L$. A choice of such a 
lifting is called a \emph{linearisation} of the action. This
induces an action on $R(X)$, and we write $R(X)^G$ for the algebra of
invariant functions. One can show that $G$ being reductive
implies that this is a finitely generated algebra, so we can form the
variety
\[ X/G = \text{Proj}\, R(X)^G.\]

While this definition of the quotient space is very simple, what we need
to understand is what its points represent. The
inclusion $R(X)^G\to R(X)$ induces a rational map $X\dashrightarrow
X/G$. The map is not defined at points $x\in X$ where every invariant
section in $R(X)^G$ vanishes. 

\begin{defn} A point $x\in X$ is called \emph{unstable} if every
  non-constant
  element of $R(X)^G$ vanishes at $x$. It is called \emph{semistable} if
  it is not unstable.
\end{defn}

If we denote the set of semistable points by $X^{ss}$, we now have a map
$X^{ss}\to X/G$. 

\begin{defn} A point 
  $x\in X$ is called \emph{polystable} if there exists an
  element $f$ of $R(X)^G$ which does not vanish at $x$, the set $X_f$
  where $f$ does not vanish is
  affine, and the action of $G$ on $X_f$ is closed (the orbit of each
  point is closed). 
  If in addition $x$ has discrete isotropy group then it is called
  \emph{stable}.
\end{defn}

We call a $G$-orbit (poly/semi)-stable if a point in the orbit is. This
does not depend on which point we choose. 
One can show that the closure of each semistable orbit contains a unique
polystable orbit and the quotient $X/G$ parametrises the polystable orbits. 
The following alternative characterisation of polystable and semistable 
points is often
useful.

\begin{prop} A point $x\in X$ is polystable if and only if for a choice of
  non-zero lift
  $\hat{x}\in L$, the orbit $G\hat{x}$ is closed in $L$. It is
  semistable if and only if the closure of the orbit $G\hat{x}$ does not
  intersect the zero section of $L$. 
\end{prop}

A central result in geometric invariant theory is the Hilbert-Mumford
numerical criterion for stability. It says that the stability of a point
can be determined by studying its orbits under one-parameter subgroups.
Let $\lambda:\mathbf{C}^*\to G$ be a nontrivial 
one-parameter subgroup and $x\in
X$. Since $X$ is projective, we can define
\[ x_0 = \lim_{t\to 0} \lambda(t)x .\]
We obtain an induced $\mathbf{C}^*$ action on the fibre $L_{x_0}$, which
has a weight $-w(x,\lambda)$.  
\begin{thm}[Hilbert-Mumford criterion] The point $x$ is 
  \begin{enumerate}
    \item stable if and only if $w(x,\lambda) > 0$ for all $\lambda$,
    \item semistable if and only if $w(x,\lambda) \geq 0$ for all
      $\lambda$,
    \item polystable if and only if $w(x,\lambda) \geq 0$ for all
      $\lambda$ with equality only if $\lambda$ fixes $x$.
  \end{enumerate}
\end{thm}

\noindent We will prove this theorem in the case of a torus action in
Section~\ref{sec:torus}.

\begin{ex} \label{ex:snp1}
  Let $X=S^n\mathbf{P}^1$, the space of unordered
  $n$-tuples of points on $\mathbf{P}^1$. We can identify such an
  $n$-tuple of points with a homogeneous polynomial of degree $n$, ie.
  with a section of $\mathcal{O}(n)$, unique up to scaling. Thus
  $X=\mathbf{P}H^0(\mathcal{O}(n))$. Let $SL(2,\mathbf{C})$ act on $X$
  via the natural action induced by the isomorphism
  $H^0(\mathcal{O}(n))\cong S^n(\mathbf{C}^2)$. 
  
  Let us test whether a given $f\in H^0(\mathcal{O}(n))$ 
  is stable for this action. 
  Choose a $\mathbf{C}^*$ subgroup of $SL(2,\mathbf{C})$ and diagonalise
  it:
  \[ \lambda\mapsto \left( \begin{array}{cc}
              \lambda^k & 0 \\
	      0 & \lambda^{-k} 
	  \end{array} \right),\]
  in $[x:y]$ coordinates on $\mathbf{P}^1$, for some $k\geq0$. In these
  coordinates we can write $f=\sum_{i=0}^n a_ix^iy^{n-i}$. As
  $\lambda\to0$, the monomials $x^iy^{n-i}$ with $2i-n\leq0$ do not tend to
  zero. Thus the closure of the orbit as $\lambda\to0$ does not contain
  the origin, unless $a_i=0$ for all $i\leq n/2$; that is as long as $f$
  does not vanish to order greater than $n/2$ at the point $[0:1]$. Changing the
  one-parameter subgroup corresponds to changing coordinates, so we can
  conclude that $f$ is semi-stable as long as it has no roots of
  multiplicity greater than $n/2$. Similarly, $f$ is polystable if it has no
  roots of multiplicity at least $n/2$ or if it has two roots of
  multiplicity $n/2$. Finally, $f$ is stable if it is
  polystable and has at least 3 distinct roots, since in this case the
  stabiliser is trivial. 
\end{ex}

\section{Kempf-Ness theorem}\label{sec:kempfness}
Let us now consider taking quotients in the symplectic category. Let
$(X,\omega)$ be a symplectic manifold with symplectic form $\omega$, and
suppose that a compact group $K$ acts on $X$, preserving $\omega$. Write
$\mathfrak{k}$ for the Lie algebra of $K$. 
To define the symplectic quotient we need a \emph{moment map} for the
action of $K$. This is a $K$-equivariant map $\mu: X\to\mathfrak{k}^*$,
such that for each $\xi\in\mathfrak{k}$ the function
$\langle\mu,\xi\rangle$ is a Hamiltonian for the vector field on $X$
induced by $\xi$. In other words,
\[ d\langle\mu,\xi\rangle = \omega(\sigma(\xi),\cdot),\]
where $\sigma:\mathfrak{k}\to \mbox{Vect}(X)$ is the infinitesimal action.
We will see shortly that a choice of moment map for the action is
equivalent to a choice of linearisation of the action in GIT. 
Given a moment map $\mu$, the symplectic quotient is defined to be
$\mu^{-1}(0)/K$. This is a symplectic manifold if $0$ is a regular value
of $\mu$ and $K$ acts properly on $\mu^{-1}(0)$. 

In order to relate this to the GIT quotient, we need some compatibility
between the two setups. Suppose that $X$ is a K\"ahler variety,
and let $L$ be an ample line bundle over $X$ endowed with a Hermitian
metric with curvature form $-i\omega$. Suppose that the action of
$K$ on $X$ preserves both the symplectic and holomorphic structures.
Given an element $\xi\in\mathfrak{k}$ which induces a holomorphic vector
field $v_\xi$ on $X$, we define a holomorphic vector field $\hat{v}_\xi$
on $L$ by
\[ \hat{v}_\xi = \tilde{v}_\xi + i\langle\mu,\xi\rangle t,\]
where $\tilde{v}_\xi$ is the horizontal lift of $v_\xi$ and $t$ is the
canonical 
vertical vector field on $L$. This gives an
infinitesimal action of $\mathfrak{k}$ on $L$, which we can extend to
the complexification $\mathfrak{g}$. Let us suppose that this
infinitesimal action can be integrated to an action of $G$. We are now
in the setup of GIT, with a complex reductive group acting on a pair
$(X,L)$. 

\begin{ex}\label{ex:moment}
  Suppose that the line bundle $L$ induces an embedding
  $X\hookrightarrow \mathbf{P}^n$, and the K\"ahler metric on $X$ is
  the pullback of the Fubini-Study metric. Suppose the group $K$ acts
  on $X$ via a representation
  \[ \rho : K\to U(n+1). \]
  In this case we can write down a moment map for the action on
  $\mathbf{P}^n$:
  \[ \mu(x).a = -i\frac{\overline{\hat{x}^t}\rho_*(a)\hat{x}}
  {\Vert\hat{x}\Vert^2},\]
  for all $a\in\mathfrak{k}$, where $\hat{x}\in\mathbf{C}^{n+1}\setminus
  \{0\}$ is a lifting of $x\in\mathbf{P}^n$. The moment map for the
  action on $X$ is just the restriction of this map to $X$. 

  The compatible linearisation is obtained by looking at the
  complexified representation $G\to GL(n+1,\mathbf{C})$. The total space
  of the line bundle $\mathcal{O}_{\mathbf{P}^n}(-1)$ is just the blowup
  of $\mathbf{C}^{n+1}$ in the origin, so we obtain an action on this
  line bundle. This induces an action on its dual, which when restricted
  to $X$ gives $L$.
\end{ex}

\begin{thm}[Kempf-Ness]\label{prop:kempf}
  A $G$-orbit contains a zero of the moment map if
  and only if it is polystable. A $G$-orbit is semistable if and only if
  its closure contains a zero of the moment map. 
\end{thm}

The key idea in the proof of this theorem is to consider the following
\emph{norm functional} on the $G$-orbit of a point $x\in X$. Choose a
non-zero lift $\hat{x}\in L_x$, and define
\[\begin{aligned}
  \phi : G/K&\to\mathbf{R} \\
  [g] &\mapsto -\log \Vert g\cdot\hat{x}\Vert.
\end{aligned}\]
Let $\xi\in\mathfrak{k}$, and consider the restriction of $\phi$ to the
geodesic $\exp(it\xi)$,
\[ f(t) = -\log\Vert\exp(it\xi)\cdot\hat{x}\Vert.\]
Computing the derivative of $\phi$ in the
direction $i\xi$, we find
\begin{eqnarray*} 
  f^\prime(0) &=& \langle\mu(g\cdot x),\xi\rangle,\\
  f^{\prime\prime}(0) &=& \Vert\sigma_x(\xi)\Vert^2,
\end{eqnarray*}
where $\sigma_x:\mathfrak{g}\to T_xX$ is the infinitesimal action.
This means that $\phi$ is convex along geodesics, and $g$ is a critical
point of $\phi$ if and only if $\mu(g\cdot x)=0$. Thinking of $\phi$ as
a function on the $G$-orbit $G\cdot x$ now, we see that a critical point
exists if and only if the $G$-orbit $G\cdot\hat{x}$ in $L_x$ is closed,
ie. $x$ is polystable. 

\begin{ex} Let us consider Example~\ref{ex:snp1} again, this time from
  the symplectic point of view. The symplectic form on
  $S^n\mathbf{P}^1=\mathbf{P}^n$ induced by the standard symplectic form on
  $\mathbf{P}^1$ is just the standard symplectic form on $\mathbf{P}^n$.
  If we denote the moment map for the action of $SU(2)$ on
  $\mathbf{P}^1$ by $\mu$, then the
  moment map for the action of $SU(2)$ on $S^n\mathbf{P}^1$ is given by
  \[ \begin{aligned} 
    \mu_n:S^n\mathbf{P}^n &\to\mathfrak{su}(2)^*,\\
    (x_1,\ldots,x_n)&\mapsto \mu(x_1)+\ldots+\mu(x_n).
  \end{aligned}
  \]
  
  \noindent We can embed $\mathbf{P}^1$ as a coadjoint orbit in
  $\mathfrak{su}(2)^*$, and the moment map for the action of $SU(2)$ is
  just this embedding. Given an invariant inner product on
  $\mathfrak{su}(2)^*$, this orbit is a sphere, and we can see that the
  moment map $\mu_n$ simply gives the centre of mass of the $n$-tuple of
  points. Zeros of the moment map correspond to balanced
  configurations, which have centre of mass zero. 
  The Kempf-Ness theorem in this case says that an $n$-tuple is
  polystable if
  and only if we can move the points to a balanced configuration by
  applying a transformation in $SL(2,\mathbf{C})$. 
\end{ex}
  
\section{Norm squared of the moment map} \label{sec:normsquared}
We use the notation from the previous section. 
Let us now choose a rational invariant inner product on $\mathfrak{k}$.
By rational we mean that for a maximal torus $T\subset K$ with Lie algebra
$\mathfrak{t}\subset\mathfrak{k}$, the inner product takes integral
values on the kernel of the exponential map $\mathfrak{t}\to T$. 
Let us define the function 
\[ \begin{split} 
   f:X\to\mathbf{R}\\ 
   f(x)=\Vert\mu(x)\Vert^2.
 \end{split}\]
The aim of this section is to study critical
points of this function and generalise the Kempf-Ness theorem to
characterise $G$-orbits of critical points of $f$ using a stability
condition. In the following
proposition we identify $\mathfrak{k}$ with its dual using the inner
product. 

\begin{prop} \label{lem:rat}
  A point $x\in X$ is a critical point of $f$ if and only if the
  vector field on $X$ induced by $\mu(x)$ vanishes at $x$. Moreover when
  $\mu(x)$ is non-zero, it generates a circle subgroup of $K$. 
\end{prop}
\begin{proof}
  To prove the first statement we differentiate $f$. Write $v$ for the
  vector field on $X$ induced by $\mu(x)$. For a tangent vector $w$ at
  $x$, 
  \[ df_x(w) = 2\langle d\mu(w),\mu(x)\rangle.\]
  Since $\langle\mu, \mu(x)\rangle$ is a Hamiltonian for $v$, we have 
  \[ df_x(w) = 2\omega(v,w).\]
  Therefore $x$ is a critical point if and only if $\omega(v,w)$
  evaluated at $x$ is zero for all $w$, ie. if $v$ vanishes at $x$. 

  To prove the second statement let $\beta=\mu(x)$ and denote
  by $T$ the closure of the subgroup of $K$ generated by $\beta$. This
  is a compact connected Abelian Lie group, hence it is a torus. 
  Letting $\mathfrak{t}$ be the Lie algebra of $T$, the moment map
  $\mu_T$ for the action of $T$ on $X$ is given by the composition of
  $\mu$ with the orthogonal projection from $\mathfrak{k}$ to
  $\mathfrak{t}$. Since by definition, $\beta\in\mathfrak{t}$, we have
  that $\mu(x)=\mu_T(x)$. Let $v_1,\ldots,v_k$ be an integral basis for
  the kernel of the exponential map from $\mathfrak{t}$ to $T$. Because
  of the rationality assumption on the inner product, what we need
  to show is that $\langle\mu_T(x),v_i\rangle$ is rational for all
  $i$, since then the orbit of $\mu(x)$ closes up to an $S^1$ orbit.
  Since $f_i=\langle\mu_T,v_i\rangle$ is the Hamiltonian function
  for the vector field induced by $v_i$, we know that $v_i$ acts on the
  fibre $L_x$ via $2\pi f_i(x)\underline{t}$. Since $\exp(v_i)=1$,
  we find that $f_i(x)$ must be an integer. 
\end{proof}

Now we define the subgroups of $G$ which will feature in the
stability condition. For a torus $T$ in $G$ with Lie algebra
$\mathfrak{t}$,  define two subalgebras of
$\mathfrak{g}$:
\begin{equation}\label{eq:subalg} 
\begin{aligned}
  \mathfrak{g}_T &:= \{\alpha\in\mathfrak{g}\,\vert\,
  [\alpha,\beta]=0\quad\text{for all }\beta\in\mathfrak{t}\} \\
  \mathfrak{g}_{T^\perp} &:= \{\alpha\in\mathfrak{g}_T\,\vert\,
  \langle\alpha,\beta\rangle = 0\quad\text{for all }
  \beta\in\mathfrak{t}\} \subset\mathfrak{g}_T.
\end{aligned}
\end{equation}

Denote the corresponding connected subgroups by $G_T$ and $G_{T^\perp}$.
Then $G_T$ is the identity component of the centraliser of $T$ and $G_{T^\perp}$ is a subgroup isomorphic to the quotient
of $G_T$ by $T$. It is a closed subgroup of $G_T$ by the following
Lemma and induction on the dimension of $T$. 

\begin{lem} Let $H$ be a compact Lie group with Lie algebra
  $\mathfrak{h}$ endowed with a rational invariant inner product. Let
  $\beta\in\mathfrak{h}$ be in the centre of $\mathfrak{h}$, and suppose
  $\beta$ generates a circle subgroup of $H$. Write $H_{\beta^\perp}$
  for the connected subgroup of $H$ generated by the Lie algebra
  \[ \mathfrak{h}_{\beta^\perp} :=\{ \alpha\in\mathfrak{h}\,\vert\, \langle
  \alpha,\beta\rangle=0\}. \]
  Then $H_{\beta^\perp}$ is closed. 
\end{lem}
\begin{proof} 
  We will use the result of Malcev~\cite{Mal45} stating that a subgroup
  of a Lie group corresponding to a Lie subalgebra is closed if and only
  if it contains the closure of all of its one-parameter subgroups.

  Let $\alpha\in\mathfrak{h}_{\beta^\perp}$. We need to show that if
  $h\in H$ is in the closure of the one-parameter subgroup
  generated by $\alpha$, then $h=\exp(\gamma)$ for some
  $\gamma\in\mathfrak{h}_{\beta^\perp}$. Let us denote the closure of
  the subgroup generated by $\alpha$ and $\beta$ by $T$. This is a
  compact connected Abelian group, so it is a torus and it contains $h$.
  Let $\mathfrak{t}$ be the Lie algebra of $T$, and let
  $\mathfrak{t}_{\beta^\perp}$ be the subalgebra of elements orthogonal
  to $\beta$. Since $\beta$ is a rational element (it generates a
  circle) and the inner product is rational, we can choose a rational
  basis for $\mathfrak{t}_{\beta^\perp}$, so the subgroup of $T$
  generated by $\mathfrak{t}_{\beta^\perp}$ is closed. In particular $h$
  is in this subgroup, so $h=\exp(\gamma)$ for some $\gamma\in
  \mathfrak{t}_{\beta^\perp}\subset\mathfrak{h}_{\beta^\perp}$.
\end{proof}

Working on the level of the compact subgroup $K$, if
$\mathfrak{t}\subset\mathfrak{k}$, then the same formulae as in
Equation~\ref{eq:subalg}
define Lie algebras $\mathfrak{k}_T, \mathfrak{k}_{T^\perp}$ and
subgroups $K_T, K_{T^\perp}$ of $K$, such that
\begin{eqnarray*}
\mathfrak{k}_T = \mathfrak{k}\cap\mathfrak{g}_T,&\quad
\mathfrak{k}_{T^\perp} =
\mathfrak{k}\cap\mathfrak{g}_{T^\perp}\\
K_T = K\cap G_T,&\quad K_{T^\perp} = K\cap G_{T^\perp}.
\end{eqnarray*}

We can now write down the stability condition that we need.

\begin{defn} 
  Let $T$ be a torus in $G$ fixing $x$. We say that $x$ is \emph{polystable
  relative to $T$} if it is polystable for the action of $G_{T^\perp}$ on
  $(X,L)$. 
\end{defn}

The main result of this section is the following. 

\begin{thm}\label{thm:stab}
  A point $x$ in $X$ is in the $G$-orbit of a critical point of
  $f$ if and only if it is polystable relative to a maximal torus
  in $G_x$, where $G_x$ is the stabiliser of $x$. 
\end{thm}

Before giving the proof, consider the effect of varying the maximal compact
subgroup of $G$. If we replace $K$ by a conjugate
$gKg^{-1}$ for some $g\in G$ and we replace $\omega$ by $(g^{-1})^*\omega$,
then we obtain a new compact group acting by symplectomorphisms. The associated
moment map $\mu_g$ is related to $\mu$ by
\begin{equation} \label{eq:moment}
\mu_g(gx) = \mathrm{ad}_g\mu(x) \in \mathrm{ad}_g\mathfrak{k},
\end{equation}

\noindent where we identify the Lie algebra of $gKg^{-1}$ with
$\mathrm{ad}_g\mathfrak{k}\subset\mathfrak{g}$. Using the inner
product on $\mathrm{ad}_g\mathfrak{k}$ induced by the bilinear form on
$\mathfrak{g}$, define the function $f_g(x) = \Vert\mu_g\Vert^2$. This
satisfies $f_g(gx)=f(x)$ by (\ref{eq:moment}) and the
$\mathrm{ad}$-invariance of the bilinear form, so in particular the
critical points of $f_g$
are obtained by applying $g$ to the critical points of $f$.

\begin{proof}[Proof of Theorem~\ref{thm:stab}]
  Suppose first that $x$ is in the $G$-orbit of a critical point of $f$.
  By replacing $K$ with a conjugate if necessary, we can assume that $x$
  itself is a critical point, so $\mu(x)$ fixes $x$. If $\mu(x)=0$ then
  Proposition~\ref{prop:kempf} implies that $x$ is polystable. If
  $\mu(x)\not=0$ then by
  Lemma~\ref{lem:rat} we obtain a circle action fixing $x$, generated by
  $\beta=\mu(x)$. Choose a maximal torus $T$ fixing $x$, containing this
  circle. Since the moment map $\mu_{T^\perp}$ for the action of
  $K_{T^\perp}$ on $X$ is the composition of $\mu$ with the
  orthogonal projection from $\mathfrak{k}$ to
  $\mathfrak{k}_{T^\perp}$, we have that $\mu_{T^\perp}(x)=0$.
  Using Proposition~\ref{prop:kempf} this implies that $x$ is polystable for
  the action of $G_{T^\perp}$. 

  Conversely, suppose $x$ is polystable for the action of $G_{T^\perp}$ for
  a maximal torus $T$ which fixes $x$.  Choose a maximal compact
  subgroup $K$ of $G$ containing $T$. Then $K_{T^\perp}$ is a maximal
  compact subgroup of $G_{T^\perp}$ and using the assumption on $x$,
  Proposition~\ref{prop:kempf} implies that $y=gx$ is in the kernel of
  the corresponding moment map $\mu_{T^\perp}$ for some $g\in G_{T^\perp}$.
  Then, for the moment map corresponding to $K$, $\mu(y)$ is contained
  in $\mathfrak{t}$ (since $T$ fixes $y$ and $\mu$ is equivariant, we
  have $\mu(y)\in\mathfrak{k}_T$), and therefore fixes $y$. 
  This means that $y$ is a
  critical point of $f$ by Proposition~\ref{lem:rat}. 
\end{proof}

We will now reformulate this stability
condition using the Hilbert-Mumford numerical criterion. Write $G_x$ for
the stabiliser of $x$. Since
$G_x$ fixes $x$, the action on the fibre $L_x$ defines a map
$G_x\to\mathbf{C}^*$. The derivative at the identity gives a linear map
$\mathfrak{g}_x\to \mathbf{C}$ which we denote by $-F_x$ in order to
match with the sign of the Futaki invariant defined later.  We 
say that $-F_x(\alpha)$ is the weight of the action of $\alpha$ on
$L_x$. According to the numerical criterion
we have the following necessary and
sufficient condition for a point $x$ to be polystable: for all one-parameter
subgroups $t\mapsto\exp(t\alpha)$ in $G_{T^\perp}$, the weight on
the central fibre $L_{x_0}$ is negative, or equal to zero if
$\exp(t\alpha)$ fixes $x$. Here $x_0$ is defined to be $\lim_{t\to
0}\exp(t\alpha)x$. In other words, the condition is that 
\[F_{x_0}(\alpha) \geqslant 0 ,\]

\noindent with equality if and only if $x$ is fixed by the one-parameter
subgroup. 

It is in\-con\-venient to restrict atten\-tion to one-parameter subgroups in
$G_{T^\perp}$ because the orthogonality condition is not a natural
one for test-configurations which we will introduce later. 
We would therefore like to be able to
consider one-parameter subgroups in $G_T$ and adapt the numerical
criterion. For a one-parameter subgroup in $G_T$ generated by
$\alpha\in\mathfrak{k}_T$ we consider the one-parameter
subgroup in $G_{T^\perp}$ generated by the orthogonal
projection\footnote{If this does not generate a
one-parameter subgroup then we can approximate it with elements of
$\mathfrak{k}_T$ that do.}
 of $\alpha$
onto $\mathfrak{k}_{T^\perp}$, which we denote by
$\overline{\alpha}$. We have 
\[ \overline{\alpha} = \alpha - \sum_{i=1}^k\langle\alpha,\beta_i\rangle\beta_i,\]
where $\beta_1,\ldots,\beta_k$ is an orthonormal basis for
$\mathfrak{t}$. 
Since $[\alpha,\mathfrak{t}]=0$ and $x$ is fixed by $T$, the central fibre
for the two one-parameter groups generated by $\alpha$ and
$\overline{\alpha}$ is the same, the only difference is the weight of
the action on this fibre. Since $F_{x_0}$ is linear, we obtain
\[ F_{x_0}(\overline{\alpha}) = F_{x_0}(\alpha) -
\sum_{i=1}^k\langle\alpha, 
\beta_i\rangle F_{x_0}(\beta_i). \]

The \emph{extremal vector field} $\chi$ is defined to be the
element in $\mathfrak{t}$ dual to the functional $F_x$ restricted to
$\mathfrak{t}$ under the inner product. In other words,
$F_x(\alpha)=\langle\alpha,\chi\rangle$ for all
$\alpha\in\mathfrak{t}$. This generates a one-parameter subgroup by the
same argument that was used in Proposition~\ref{lem:rat}. 
If we now choose the orthonormal basis
$\beta_i$ such that $\beta_1=\chi/\Vert\chi\Vert$, then the previous
formula reduces 
to 
\[ F_{x_0}(\overline{\alpha}) =
F_{x_0}(\alpha)-\langle\alpha,\chi\rangle. \]

If we define this expression to be $F_{x_0,\chi}(\alpha)$, then the
stability condition is 
equivalent to $F_{x_0,\chi}(\alpha)\geqslant 0$ 
for all one-parameter subgroups generated by
$\alpha\in\mathfrak{k}_{T}$ with equality only if the one-parameter
subgroup fixes $x$. We therefore obtain the following

\begin{thm}
  A point $x\in X$ is in the $G$-orbit of a critical point of $f$, if
  and only if for each one-parameter subgroup of $G$ generated by an
  element $\alpha\in\mathfrak{k}_T$ we have
  \[ F_{x_0,\chi}(\alpha)\geqslant 0, \]
  with equality only if $\alpha$ fixes $x$. Here $T$ is a maximal torus
  fixing $x$ and $\chi$ is the corresponding extremal vector field. 
\end{thm}

We now ask what the infimum of the function $f=\Vert\mu\Vert^2$ is
on a $G$-orbit. 

\begin{thm}\label{thm:destablwrbound}
  Let $x\in X$, and let $\alpha\in\mathfrak{k}_T$ generate
  a one-parameter subgroup such that
  the weight $F_{x_0,\chi}(\alpha)<0$. Then 
  \[ \inf_{g\in G} \Vert \mu(g\cdot x)\Vert^2\geq \Vert\chi\Vert^2 +
  \frac{F_{x_0,\chi}(\alpha)^2}{\Vert\alpha\Vert^2}.\]
\end{thm}
\begin{proof} Suppose $F_{x_0}(\alpha)<0$. We can arrange this by adding
  a multiple of $\chi$ to $\alpha$ if necessary. Consider the function
  \[ f(t) = \langle \mu(\exp(it\alpha)\cdot x),\alpha\rangle.\]
  Computing the derivative of $f$, we find
  \[ f^\prime(t) = \Vert\sigma_{\exp(it\alpha)\cdot x}(\alpha)\Vert^2
  \geq 0,\]
  so that $f$ is non-decreasing. Letting $t\to -\infty$ we get
  \[ f(t)\to\langle\mu(x_0),\alpha\rangle = -F_{x_0}(\alpha),\]
  but $f(0)=\langle\mu(x),\alpha\rangle$, so that we must have
  $\langle\mu(x),\alpha\rangle \geq -F_{x_0}(\alpha)$. This implies
  \begin{equation}\label{eq:ineq1}
    \Vert\mu(x)\Vert^2\geq\frac{F_{x_0}(\alpha)^2}{\Vert\alpha
  \Vert^2}.
\end{equation}

  We now need to modify $\alpha$ carefully to get the result we want.
  Let $\overline{\alpha}$ be the component of $\alpha$ orthogonal to
  $\chi$ (the same remark as above applies if this does not generate a
  one-parameter subgroup). Since by our assumption
  $F(\overline{\alpha})<0$, we can choose  
  a scalar $\lambda>0$ such that
  $F_{x_0}(\lambda\overline{\alpha})=-\Vert\lambda\overline{\alpha}
  \Vert^2$. Now define $\gamma=\lambda\overline{\alpha}-\chi$. As
  before, the central fibre for the one-parameter subgroup generated by
  $\gamma$ is $x_0$, just the weight is changed to 
  \[
  F_{x_0}(\lambda\overline{\alpha}-\chi)=-\Vert\lambda\overline{\alpha}
  \Vert^2-\Vert\chi\Vert^2,\]
  which is negative. Since
  \[
  \frac{F_{x_0}(\gamma)^2}{\Vert\gamma\Vert^2}=\Vert\lambda
  \overline{\alpha}\Vert^2+\Vert\chi\Vert^2=\frac{F_{x_0}(
  \overline{\alpha})^2}{\Vert\overline{\alpha}\Vert^2}
  +\Vert\chi\Vert^2,
  \]
  using Inequality~\ref{eq:ineq1} we obtain
  \[ \Vert\mu(x)\Vert^2\geq\frac{F_{x_0}(\overline{\alpha})^2}{
  \Vert\overline{\alpha}\Vert^2}+\Vert\chi\Vert^2.\]
  Finally, since $F_{x_0}(\overline{\alpha})=F_{x_0,\chi}(\alpha)$, and
  $\Vert\overline{\alpha}\Vert\leq\Vert\alpha\Vert$, we get the required
  inequality for $\Vert\mu(x)\Vert^2$. By replacing $x$ by $g\cdot x$ and
  $\alpha$ by $\mathrm{ad}_g(\alpha)$ we obtain the same inequality for
  $\Vert\mu(g\cdot x)\Vert^2$. 
\end{proof}

Note that in this theorem we get the strongest inequality if we choose
$\alpha$ orthogonal to the chosen torus of automorphisms since that
minimises $\Vert\alpha\Vert$. Note that if $x$ is a critical point, then
$\Vert\mu(x)\Vert=\Vert\chi\Vert$, so this result implies a weak version
of Theorem~\ref{thm:stab}. This is the form in which it will be used in
Chapter~\ref{chap:stabvar} to prove a necessary condition for a variety 
to admit an extremal metric.

\section{Modulus of stability}\label{sec:modulus}

Choose an invariant inner product on $\mathfrak{k}$. Let $x\in X$ be a
polystable point, and
write $\pi_x : \mathfrak{k}\to \mathfrak{k}_x$ for the orthogonal
projection onto the stabiliser of $x$. 
Define the \emph{modulus of stability} $\lambda$ of $x$ by
\[\lambda = \inf_\alpha \frac{w(x,\alpha)}{\Vert\alpha-\pi_x(\alpha)
\Vert}, \]
where the infimum is over all $\alpha\in\mathfrak{k}\setminus
\mathfrak{k}_x$ generating
one-parameter subgroups. This is an invariant of the orbit of $x$, and
measures how far this orbit is from being unstable. Note that $\lambda$
is strictly positive. To see this, note that we can restrict to $\alpha$ in the
orthogonal complement $\mathfrak{k}_x^\perp$ and by continuity we can
extend the function
\[ \psi(\alpha) = \frac{w(x,\alpha)}{\Vert\alpha\Vert} \]
to all non-zero $\alpha\in\mathfrak{k}_x^\perp$. The unit ball of
$\mathfrak{k}_x^\perp$ is
compact so $\psi$ achieves its infimum at some $\beta$ which may or may not
generate a one-parameter subgroup. To see that it does generate a
one-parameter subgroup, restrict
attention to the complex torus generated by $\beta$ and use
the arguments in Section~\ref{sec:torus}. Since $x$ is polystable,
$\psi(\beta)>0$ and so $\lambda>0$. 

\begin{prop} Let $x$ be polystable with modulus of stability $\lambda$. Let 
  $x_0$ be the limit of $x$ under a one-parameter subgroup which does
  not fix $x$. 
  Then $\Vert\mu(x_0)\Vert\geq\lambda$.
\end{prop}
\begin{proof}
  Let $\alpha\in\mathfrak{k}$ generate the one-parameter subgroup. We
  can assume that
  $\alpha$ is orthogonal to the stabiliser of $x$. By the definition of
  the weight, we have 
  \[\langle\mu(x_0),\alpha\rangle =
  -w(x,\alpha).\]
  From this we obtain 
  \[ \Vert\mu(x_0)\Vert\geq\frac{w(x,\alpha)}{\Vert\alpha\Vert}\geq
  \lambda. \]
\end{proof}

For each point $x\in X$, the infinitesimal action of $K$ induces a
linear map $\sigma_x:\mathfrak{k}\to T_x X$. Using the metrics on $X$ and
$\mathfrak{k}$ we form its adjoint $\sigma_x^*$. Suppose that the line bundle
$L$ over $X$ is very ample and induces an embedding
$X\subset\mathbf{P}^{n-1}$. 

\begin{thm}\label{thm:eigenval}
  Let $\mu(x)=0$ and let the modulus of stability of $x$ be
  $\lambda$. Assume for simplicity that $x$ has trivial stabiliser. 
  Then the smallest 
  eigenvalue of $\sigma_x^*\sigma_x$ is bounded below by
  $2\lambda^2/n$.  
\end{thm}
\begin{proof}
  Consider the moment map restricted to a $G$-orbit, 
  \begin{eqnarray*}
    \phi :  G&\to&\mathfrak{k}\\
    g&\to& \mu(g\cdot x).
  \end{eqnarray*}
  We can compute
  \[ \langle d\phi_e(i\xi),\eta\rangle = \omega_x(\sigma_x(\eta),
  J\sigma_x(\xi))=\langle\eta,\sigma_x^*\sigma_x(\xi)\rangle,
  \]
  so the operator $\sigma_x^*\sigma_x$ is given by the derivative of $\phi$
  in the $i\mathfrak{k}$ directions at the identity. To prove the
  result, we therefore need to show that for all $\xi\in\mathfrak{k}$
  \[
     \langle d\phi_e(i\xi),\xi\rangle\geq \frac{2\lambda^2}{n}\Vert\xi
     \Vert^2, \]
  and it is enough to restrict to the case when $\xi$ generates a
  $\mathbf{C}^*$ action.

  Suppose the line bundle $L$ induces an embedding
  $X\subset\mathbf{P}(V)$ with $\dim V=n$, and
  $\xi$ generates a $\mathbf{C}^*$-action on $V$. Let
  $V=\bigoplus V_i$ be the weight decomposition of $V$, so that the
  action on $V_i$ has weight $w_i$, and $w_1\leq w_2\leq\ldots\leq w_n$.
  Choose an orthonormal basis $\{ e_i\}$ with $e_i\in V_i$. We can
  assume without loss of generality that $x$ is in the orbit of
  $x_0=(1,1,\ldots,1)$ in these 
  coordinates (if there were fewer non-zero coordinates, then we would
  get a sharper inequality in the end). The moment map is given by (see
  Example~\ref{ex:moment})
  \[ \mu(t\cdot x_0) = \frac{\sum w_i e^{2w_i t}}{\sum e^{2w_i t}}. \]
  We want to estimate the derivative of $\mu(t\cdot x_0)$ with respect
  to $t$ at the point $t_0$ for which $x=t_0\cdot x_0$ (recall that
  $\mu(x)=0$). We have
  \[ \frac{d}{dt}\mu(t_0\cdot x_0) = \frac{\sum 2w_i^2 e^{2w_it_0}}{\sum
  e^{2w_it_0}}. \]
  Let us suppose without loss of generality 
  that $t_0\geq 0$. Then $\sum w_i^2e^{2w_it_0}\geq
  w_n^2e^{2w_nt_0}$ and $\sum e^{2w_it_0}\leq ne^{2w_nt_0}$, so we obtain
  \[ \frac{d}{dt}\mu(t_0\cdot x_0)\geq\frac{2w_n^2}{n}.\]
  Since by the definition of the modulus of stability $\lambda$ we have
  $w_n\geq\lambda\Vert\xi\Vert$, the proof is complete.
\end{proof}

\section{Torus actions}\label{sec:torus}

In the case of a torus action, stability can be understood by analysing
the weights of the action.
Let
$T^c=(\mathbf{C}^*)^k$ act on $\mathbf{P}(V)$ via a representation of $T$
on $V$. Choose a basis $\{e_1,\ldots,e_n\}$ for $V$ such that the
action is diagonal, given by weights $\alpha_j\in\mathfrak{t}^*$. The
action is given by
\[ \exp(\xi)e_j = \exp(i\langle \xi,\alpha_j\rangle)e_j, \quad\text{for all
}\xi\in\mathfrak{t}. \]
Let $\{X_1,\ldots,X_n\}$ be the dual basis for $V^*$, on which the
corresponding action is given by the same weights. 
Invariant monomials are given by $\prod_i X_i^{a_i}$ such that $\sum a_i\langle
\xi,\alpha_i\rangle = 0$ for all $\xi$, ie. 
\[ \sum_i a_i\alpha_i = 0. \]
Invariant sections of $O(m)$ over $\mathbf{P}(V)$ are sums of these
monomials with
$\sum a_i = m$.

Let $x\in\mathbf{P}(V)$ and $\hat{x}\in V$ a non-zero lifting. 
We define the \emph{weight polytope} $\Delta_x$ of $x$ to be the closed
convex hull 
of the weights acting nontrivially on $x$:
\[ \Delta_x = \overline{\mbox{co}}\{\alpha_j | X_j(\hat{x})\not= 0\}
\subset\mathfrak{t}^*. \]
Note that
$\Delta_x$ is contained in a proper affine subspace if and only if $x$ has
non-discrete stabiliser. In the following theorem when referring to the
interior of $\Delta_x$, we are considering $\Delta_x$ to be a subset of
the minimal affine subspace containing it.

\begin{thm}
  Let $x\in \mathbf{P}(V)$. We have
  \begin{enumerate}
    \item $x$ is semistable if and only if $\Delta_x$ contains the
      origin.
    \item $x$ is polystable if and only if $\Delta_x$ contains the
      origin in its interior.
    \item $x$ is stable if and only if $\Delta_x$ contains the origin in
      its interior and $\Delta_x$ is not contained in any proper
      subspace.
  \end{enumerate}
\end{thm}
\begin{proof}
  \begin{enumerate}
    \item By definition $x$ is semistable if and only if there is an
      invariant section of $O(m)$ for some $m$ which does not vanish at
      $x$. Invariant monomials which do not vanish at $x$ are products
      $\prod_i X_i^{a_i}$ with $\sum a_i\alpha_i=0$ and $a_i=0$ whenever
      $X_i(\hat{x})=0$. Such a section exists if and only if zero is contained
      in $\Delta_x$. 
    \item Let us first show that if $\Delta_x$ contains the origin in
      its interior then $x$ is polystable. For simplicity let us assume
      that $\hat{x}=(1,1,\ldots,1)\in V$. By the hypothesis we can choose
      non-zero $a_i$'s such that $s=\prod_i X_i^{a_i}$ is an invariant
      monomial. The set where $s$ does not vanish is the affine set 
      $(\mathbf{C}^*)^n$,
      and we need to show that the action of $T^c$
      on this set is closed. It is
      enough to show that the action of the $\mathbf{R}^k$ component of the torus
      is closed since $(S^1)^k$ is compact. Define the map
      \begin{eqnarray*}
	  \psi : (\mathbf{C}^*)^n&\to & \mathbf{R}^n\\
	  (z_1,\ldots,z_n)&\mapsto & (\log|z_1|,\ldots,\log|z_n|).
      \end{eqnarray*}
      The images of the orbits under this map are subspaces of
      $\mathbf{R}^n$, so since $\psi$ is continuous, the orbits are
      closed. 

      Conversely, suppose the origin is on the boundary of $\Delta_x$.
      Choose $\xi$ to be orthogonal to the face containing the origin,
      pointing inwards to $\Delta_x$. Then the lift of
      $\lim_{t\to\infty}\exp(it\xi)x$ in $V$ is given by 
       \[\sum_{j:\langle\xi,\alpha_j\rangle=0}e_j.\]
      Since $x=(1,1,\ldots,1)$, this is not in the orbit of $x$, so this orbit
      is not closed, and $x$ is not polystable. 
    \item This follows from the remark before the theorem, since $x$ is
      stable if and only if it is polystable with discrete stabiliser. 
  \end{enumerate}
\end{proof}

Using the moment map we can give a different description of $\Delta_x$.
The compact torus $T$ acts on the orbit $T^c(x)$, and the interior of 
$\Delta_x$ is the image of the moment map for this action. This follows
from Atiyah's convexity theorem (see~\cite{Ati82}). Part 2 of 
the above theorem thus confirms the Kempf-Ness theorem in the case of a
torus action.

Let us introduce a rational inner product on $\mathfrak{t}$ so that we can
identify $\mathfrak{t}$ with $\mathfrak{t}^*$. A one-parameter subgroup
of $T$ corresponds to an integral element $\xi\in\mathfrak{t}$. The
limit $\lim\limits_{t\to -\infty}\exp(it\xi)\cdot x$ in $\mathbf{P}(V)$
is the sum of 
those $e_i$ for which $\langle \xi,\alpha_i\rangle$ is maximal and
$X_i(\hat{x})\not=0$. The
weight on the central fibre is therefore
\[ F(\xi) = \max_{i:X_i(\hat{x})\not=0}\{\langle\xi,\alpha_i\rangle\}.\]
The Hilbert-Mumford criterion says that $x$ is stable if and only if
$F(\xi)>0$ for each integral $\xi$. This
means that for any rational hyperplane in $\mathfrak{t}$ there are some
$\alpha_i$ on both sides of it. This is equivalent to the origin being
contained in $\Delta_x$, since the $\alpha_i$ are rational. The
argument with semistable and polystable points is similar, and this
gives a proof of the Hilbert-Mumford criterion for torus actions. 

It is easy to see that if $x$ is polystable, 
then the modulus of stability  of $x$ is the distance
of the boundary of $\Delta_x$ from the origin. 
If on the other hand $x$ is unstable, then the worst destabilising configuration 
(in the sense that $-w(x,\alpha)/\Vert\alpha\Vert$ is maximal) is given by
$-\xi$, where $\xi$ is the closest point of $\Delta_x$ to the origin.
The weight of this is $-\Vert\xi\Vert^2$, so we
see that in this case 
\[ \inf_{t\in T^c}\Vert\mu(t\cdot x)\Vert =
\sup_{\alpha}\frac{-w(x,\alpha)}{
\Vert\alpha\Vert}.\]
Note that this is a strengthening of Theorem~\ref{thm:destablwrbound} in
the case of torus actions. In fact this stronger version is true in
general, but we do not need it (see Kirwan~\cite{Kir84}).

Finally we describe relative polystability for a torus action. 
In
Section~\ref{sec:normsquared} we saw that a point $x$ is a critical
point of the norm squared of the moment map if $\mu(x)$ (as an element
in $\mathfrak{t}$) fixes $x$. 
Recall that the image of the $T^c$-orbit of $x$ under the
moment map is the interior of $\Delta_x$. Identifying $\mathfrak{t}$
with $\mathfrak{t}^*$ using the inner product we find that a vector
$\xi\in\mathfrak{t}$ fixes $x$ if and only if $\xi$ is orthogonal to an
affine subspace containing 
$\Delta_x$. We thus have the following

\begin{thm} The point $x$ is relatively polystable if and only if the
  orthogonal projection of the origin onto the minimal affine subspace
  containing $\Delta_x$ is in the interior of $\Delta_x$. 
\end{thm}

\chapter{Extremal metrics}\label{chap:extremal}
Extremal metrics were defined by Calabi~\cite{Cal82} as an attempt to
find canonical metrics in a given K\"ahler class on a K\"ahler manifold.
For the definition let $(M,\omega_0)$ be a K\"ahler manfiold.
For any K\"ahler metric $\omega$ in the same cohomology class as
$\omega_0$, define the \emph{Calabi functional}
\[ f(\omega) = \int_M (S(\omega)-\hat{S})^2 \frac{\omega^n}{n!}, \]
where $S(\omega)$ is the scalar curvature, $\hat{S}$ is its average 
and $n$ is the dimension of
$M$. We will see that $\hat{S}$ is independent of the choice of
$\omega\in[\omega_0]$. 
The metric $\omega$ is called \emph{extremal} if it is a critical
point of this functional. Calabi showed that the Euler-Lagrange
equation for this variational problem is that the gradient of the scalar
curvature is a holomorphic vector field. The problem is the existence
and uniqueness of extremal metrics. The uniqueness problem has been
solved by Mabuchi~\cite{Mab04_2} in the algebraic case and
Chen-Tian~\cite{CT05} in general, in the sense that the extremal
metric is unique up to the action of holomorphic automorphisms. 

In Section~\ref{sec:FMprelim} we will recall some of the more elementary
theory of extremal metrics. Then in
Section~\ref{sec:scalmoment} we explain how the scalar curvature arises
as a moment map for an infinite dimensional symplectic action. The
Calabi functional then appears as the norm squared of the moment map, so
the theory developed in the previous chapter becomes relevant to the
study of extremal metrics. This point of view is used to motivate the
definition of K-stability in Chapter~\ref{chap:stabvar}. Also,
all of the concepts introduced in Section~\ref{sec:FMprelim} can be seen
as special cases of the constructions in
Chapter~\ref{chap:finiteGIT}. 

\section{Futaki invariant and Mabuchi functional}\label{sec:FMprelim}
As above, let $(M,\omega_0)$ be a K\"ahler manifold. For K\"ahler
metrics $\omega$ in the same cohomology class as $\omega_0$ we 
define the following three functionals:
\[ \begin{split}
  &f(\omega)=\int_M (S(w)-\hat{S})^2\frac{\omega^n}{n!},\\
   &g(\omega)=\int_M
  |Ric(\omega)|^2\frac{\omega^n}{n!},\quad h(\omega)=\int_M 
  |Riem(\omega)|^2 \frac{\omega^n}{n!},
\end{split}\]
where $Ric$ is the Ricci curvature, and $Riem$ is the full Riemannian
curvature. Calabi showed that these three functionals differ by
constants depending only on the K\"ahler class, so their critical points are
the same. Calabi showed (see~\cite{Cal82})

\begin{prop} A metric $\omega$ is a critical point of $f$ if and only if the
  gradient of $S(\omega)$ is a holomorphic vector field. Such a metric
  is called an \emph{extremal metric}.
\end{prop}

In particular a metric with constant scalar curvature (cscK) is an extremal
metric, but there are also examples with non-constant scalar curvature
(see eg. Section~\ref{sec:extremalruled}).
In fact if we fix a maximal torus of holomorphic automorphisms of the
manifold, then we can determine a priori what the gradient vector field of
the scalar curvature of an
extremal metric is if one exists. First of all the total scalar curvature
is an invariant of the K\"ahler class, since
\[ \int_M S(\omega)\frac{\omega^n}{n!} = \int_M
\rho\wedge\frac{\omega^{n-1}}{(n-1)!} =
\frac{2\pi c_1(M)\cup[\omega]^{n-1}}{(n-1)!},\]
where $\rho$ is the Ricci form of $\omega$. 
Since the volume is also an invariant of the
K\"ahler class, we see that the average scalar curvature $\hat{S}$ is
fixed. 
In order to refine this, we need to define
the Futaki invariant. This was introduced by Futaki in~\cite{Fut83} as
an obstruction to the existence of a K\"ahler-Einstein
metric. We first need some preliminaries about holomorphic vector
fields. 

Given a complex valued function $f:M\to\mathbf{C}$ and a metric $\omega$, 
we can define a
vector field $X_f$ of type $(1,0)$ by
\[ X_f = \sum_\alpha g^{\overline{\beta}\alpha}\frac{\partial
f}{\partial\overline{z^\beta}}\frac{\partial}{\partial z^\alpha},\]
where $g$ is the metric corresponding to $\omega$. 
This is the $(1,0)$-part of the gradient of $f$. 
\begin{defn} We say $f:M\to\mathbf{C}$ is a \emph{holomorphy potential}
  if $X_f$ is a holomorphic vector field. Denote by $\mathfrak{h}$ the
  Lie algebra of holomorphic vector fields and by $\mathfrak{h}_1$ the
  subspace of holomorphic vector fields of the form $X_f$.
\end{defn}

It is shown in Kobayashi~\cite{Kob95} that when $M$ is a projective
variety, then the space $\mathfrak{h}_1$
coincides with the space of holomorphic vector fields that can be lifted
to an ample line bundle over $M$. We therefore write $\mbox{Aut}(M,L)$
for the group of automorphisms generated by $\mathfrak{h}_1$. 

The Futaki invariant is defined as a functional
\[ \begin{aligned}
  \mathcal{F}_\omega :\mathfrak{h}_1&\to\mathbf{C}\\
  X_f &\mapsto\int_M f(S(\omega)-\hat{S})\frac{\omega^n}{n!},
\end{aligned}
\]
The point is that this functional is independent of the choice of
K\"ahler metric $\omega$ in the class $[\omega_0]$.
\begin{prop}[cf. Calabi~\cite{Cal85}] 
  The functional $\mathcal{F}$ is independent of the choice
  of representative of the K\"ahler class. 
\end{prop}
Thus, if there is a metric $\omega\in[\omega_0]$ which has constant
scalar curvature, then $F(X_f)=0$ for all $X_f\in\mathfrak{h}_1$. This
gives an obstruction to the existence of a cscK metric and was the
original context in which the Futaki invariant was used. 

Choose a maximal compact subgroup $K$ of $\mbox{Aut}(M,L)$, and a maximal torus
inside $K$ with
Lie algebra $\mathfrak{t}\subset\mathfrak{h}_1$. Let
$\mathfrak{t}^\mathbf{C}\subset\mathfrak{h}_1$ be the complexification
of $\mathfrak{t}$. 
We define an inner product on $\mathfrak{t}^\mathbf{C}$, following Futaki and
Mabuchi~\cite{FM95} (they defined the inner product on a larger algebra,
but we do not need that here). We choose a metric $\omega\in[\omega_0]$
which is invariant under $K$, and define
\[ \langle X_f,X_g\rangle = \int_M fg\frac{\omega^n}{n!},\]
where we normalise $f,g$ to have integral zero on $M$. 
It is shown in~\cite{FM95} that this 
is invariant of the representative of the K\"ahler class
chosen.
This inner product is positive definite on
$\mathfrak{t}$, and by duality the Futaki invariant defines a vector
field $\chi\in\mathfrak{t}$. This is called the \emph{extremal vector
field}, and it only depends on the K\"ahler class and the choice of
$K$ (it is in the centre of $\mathfrak{k}$). 
If we change $K$ to a conjugate, the new extremal vector field
is a
conjugate of $\chi$. In particular the norm $\langle\chi,\chi\rangle$ is
an invariant of the K\"ahler class.
From the definition of the Futaki invariant we see that it is given
by $X_{\pi(S(\omega))}$ where $\pi(S(\omega))$ is the $L^2$-orthogonal
projection of the scalar curvature $S(\omega)$ onto the space of
holomorphy potentials. The fact that
$X_{\pi(S(\omega))}\in\mathfrak{h}_1$ lies in $\mathfrak{t}$ is shown
in~\cite{FM95}. This means that the gradient of the scalar curvature of
an extremal metric if it
exists is given by $X_{\pi(S(\omega))}$ 
for any $K$-invariant $\omega$ in the K\"ahler class.

Note that if we normalise $\pi(S(\omega))$ to have zero mean, then we
have
\begin{equation}\label{eq:celwr}
  \begin{split}
  \int_M (S(\omega)-\hat{S})^2\frac{\omega^n}{n!} = &\int_M
  [S(\omega)-\hat{S}-\pi(S(\omega))]^2\frac{\omega^n}{n!} +
\Vert\pi(S(\omega))\Vert^2_{L^2} \\ &\geq 
\langle\chi,\chi\rangle.
\end{split}
\end{equation}
This gives a lower bound on the Calabi functional for $K$-invariant
metrics which 
is achieved by a metric $\omega$
if and only if $\omega$ is an extremal metric. 

We now define the Mabuchi functional (see \cite{Mab86}) which is a
functional on the set of K\"ahler metrics in a fixed K\"ahler class,
whose critical points are constant scalar curvature metrics. Write
$\mathcal{K}$ for the space of metrics in the K\"ahler class
$[\omega_0]$. The tangent space to $\mathcal{K}$ at a metric $\omega$
can be identified as
\[ T_\omega\mathcal{K} = \left\{\phi\in C^\infty(M)\,\left|\, \int_M
\phi\frac{\omega^n}{n!}=0\right.\right\}.\]
We define the Mabuchi functional by its variation as follows:
\[ d\mathcal{M}_\omega(\phi) = -\int_M\phi(S(\omega)-\hat{S})\frac{
\omega^n}{n!}.\]
This defines a closed 1-form on $\mathcal{K}$, and so it defines
$\mathcal{M}$ up to a constant since $\mathcal{K}$ is contractible.
From the definition it is clear that the critical points of the Mabuchi
functional are metrics of constant scalar curvature. The space
$\mathcal{K}$ can be thought of as an infinite dimensional symmetric
space (analogous to $SL(n,\mathbf{C})/SU(n)$), and the Mabuchi
functional is convex along geodesics. Therefore the existence of a
constant scalar curvature metric in $\mathcal{K}$ is expected to be
equivalent to the properness of $\mathcal{M}$. This has been shown in
the case of K\"ahler-Einstein metrics by Tian~\cite{Tian97} (see
also~\cite{PSSW}).

\section{Scalar curvature as a moment map}\label{sec:scalmoment}

In this section we show how the scalar curvature arises as the moment
map in an infinite dimensional symplectic quotient problem. This was
shown by Donaldson in~\cite{Don97}. We follow here the computation in
local coordinates given by Tian~\cite{Tian00}.

Let $(M,\omega)$ be a symplectic manifold. An \emph{almost complex
structure} 
$J$ on $M$ is an endomorphism $J:TM\to TM$ such that $J^2=-Id$. We say
that the almost complex structure is \emph{compatible with $\omega$} if
the tensor $g_J$ defined by 
\[ g_J(u,v)=\omega(u,Jv) \]
is symmetric and positive definite, ie. it defines a Riemannian metric.
The almost complex structure $J$ is \emph{integrable} if we can define
local holomorphic coordinates on $(M,J)$. If $J$ is integrable and
compatible with $\omega$, then together they define a K\"ahler structure
on $M$. 

Let us denote by $\mathcal{J}$ the space of (integrable) complex
structures on $M$, compatible with the symplectic form. This is an
infinite dimensional manifold with tangent space at $J$ given by
\[ T_J\mathcal{J} = \{A:TM\to TM\,:\, AJ+JA=0,\,\omega(u,Av)=\omega
(v,Au)=0\}.\]
For $A\in T_J\mathcal{J}$ define $\mu_A(u,v)=\omega(u,Av)$. We can check
that $\mu_A(Ju,Jv)=-\mu_A(u,v)$ and $\mu_A(u,v)=\mu_A(v,u)$, and
conversely these symmetric, anti $J$-invariant sections of $T^*M\otimes
T^*M$ can be identified with $T_J\mathcal{J}$. This tangent space has a
natural complex structure induced by $J$, namely 
\[ (J\mu)(u,v) = -\mu(Ju,v).\]
This complex structure is integrable
and vector fields on $M$ acting on $\mathcal{J}$ preserve this complex
structure. 
There is also an $L^2$ inner
product on $T_J\mathcal{J}$ 
induced by $g_J$. Together these induce a K\"ahler structure on
$\mathcal{J}$. 

Assume for simplicity that $H^1(M)=0$ and let $K=Symp^0(M,\omega)$ be
the connected component of the identity in the group of
symplectomorphisms of $M$. This acts on
$\mathcal{J}$ preserving its symplectic structure, and we wish to
identify a moment map for this action. The Lie algebra $\mathfrak{k}$
of $K$ can be
identified with smooth functions on $M$ with zero mean, via the
Hamiltonian construction. For an element $J\in\mathcal{J}$ we denote by
$S(J)$ the scalar curvature of $g_J$, and by $\hat{S}$ the average scalar
curvature, which is independent of $J$. We use the complex scalar
curvature which is half of the usual Riemannian one. 
We can now state the result proved
in~\cite{Don97}. 

\begin{prop} \label{prop:scalmoment}
  The map $J\mapsto 4(S(J)-\hat{S})$ is an equivariant moment map for the action
  of $K$ on $\mathcal{J}$, where we have identified $\mathfrak{k}$ with
  its dual via the $L^2$ pairing. 
\end{prop}

To prove this, we need to compute two maps:
\begin{eqnarray*}
  P&:&C^\infty_0(M)\to T_J\mathcal{J},\\
  Q&:&T_J\mathcal{J}\to C^\infty_0(M),
\end{eqnarray*}
where $P$ is the infinitesimal action of $Symp^0(M,\omega)$ on
$\mathcal{J}$ and $Q$ is the infinitesimal change in the scalar
curvature of $g_J$ induced by an element in $T_J\mathcal{J}$. 
To do the computation we will choose local normal 
coordinates $x_1,\ldots,x_{2n}$. Since $g_J$ is K\"ahler, $dJ(0)=0$. 

\begin{prop} Identifying $T_J\mathcal{J}$ with symmetric 2-tensors as
  above, we have
  \[ P(H)_{ij} = J^k_iH_{jk}+J^k_jH_{ik}.\]
\end{prop}
\begin{proof}
  Let us denote by $X_H$ the Hamiltonian vector field corresponding to
  $H$. The components of $X_H$ are given by $X^i_H=-J^i_kg^{jk}H_j$. We
  need to compute $\mathcal{L}_{X_H}J$, and identify it with a symmetric
  2-tensor. 
  
  Let us write
  \[ \mathcal{L}_{X_H}J\left(\frac{\partial}{\partial x^j}\right) =A^i_j
  \frac{\partial}{\partial x^i}.
  \]
  Using $(\mathcal{L}_uJ)(v)=\mathcal{L}_u(J(v))-J\mathcal{L}_uv$ we can
  compute 
  \[ A^i_j=J^k_jJ^i_pg^{pq}H_{qk}+g^{ik}H_{kj}. \]
  Since $A$ is related to $P(H)$ by $P(H)_{ij}=\omega_{ik}A^k_j$, we get
  \[ P(H) = J^k_jH_{ik} + J^k_iH_{jk}.\]
\end{proof}

\begin{prop} We have $Q(\mu) = \frac{1}{2}\mu_{jk,kj}$.
\end{prop}
\begin{proof}
 Let us choose a path of complex structures $J_t$ such that
 $\frac{d}{dt}\left|_{t=0}J_t\right.=\mu$. Then the variation of the
 induced metrics $g_t$ is also $\mu$. Since the Christoffel symbols of
 $g_t$ are of order $t$, we have 
 \[ R^t_{ijkl} = \frac{\partial^2 g_{t,il}}{\partial x^j\partial x^k} -
 \frac{\partial^2 g_{t,ik}}{\partial x^j\partial x^l}+O(t^2).\]
 The scalar curvature is therefore 
 \[ S(g_t) = \frac{1}{2}g^{ik}_tg^{jl}_t \left(
 \frac{\partial^2 g_{t,il}}{\partial x^j\partial x^k} -
 \frac{\partial^2 g_{t,ik}}{\partial x^j\partial x^l}\right)
 +O(t^2).\]
 Differentiating at $t=0$, we obtain 
 \[ Q(\mu) =
 -\mu_{ik}\mathrm{Ric}_{ik}+\frac{1}{2}(\mu_{jk,kj}-\mu_{jj,kk}).\]
 Since $\mathrm{Ric}$ is $J$-invariant and $\mu$ is anti $J$-invariant,
 the first term vanishes. Since $\mu$ is anti J-invariant, the trace
 $\mu_{jj}$ vanishes and so we get the result we wanted.
\end{proof}

Putting the previous two results together, we can verify
proposition~\ref{prop:scalmoment}. We need to check that
\[ 4(Q(\mu),H)_{L^2}=\Omega(P(H),\mu),\]
where $\Omega$ is the symplectic form on $\mathcal{J}$ induced by the
complex structure defined above and the $L^2$ product. We have
\begin{equation*}
  \begin{split}
    (Q(\mu),H)_{L^2} &= \frac{1}{2}\int_M \mu_{jk,kj}H\frac{\omega^n}{n!}\\
    \Omega(P(H),\mu)&= -(P(H),J\mu)_{L^2}.
  \end{split}
\end{equation*}
Since $(J\mu)_{ij} = -J^k_i\mu_{kj}$, we have
\begin{equation*} 
  \begin{split}
    -(P(H),J\mu)_{L^2} &= \int_M (J^k_iH_{jk}+J^k_jH_{ik})J^l_p\mu_{lq}
  g^{ip}g^{jq}\frac{\omega^n}{n!}\\ &= 2\int_M H_{jk}\mu_{lq}g^{kl} g^{jq}
  \frac{\omega^n}{n!}.
\end{split}
\end{equation*} 
Integrating by parts we get the required result. 

We now see that the norm squared of this moment map is the
Calabi functional up to a scalar multiple. We can therefore
hope to apply the results of
Section~\ref{sec:normsquared} to characterise the complexified orbits of
critical points of the functional in terms of stability. Unfortunately
the complexification of $Symp^0(M,\omega)$ does not exist, but we can
think of the orbits of this complexification as follows. The
infinitesimal action of $Symp^0(M,\omega)$ defines a distribution on
$\mathcal{J}$. It can be shown that the complexification of this
distribution is integrable, so it defines a foliation of $\mathcal{J}$.
We think of the leaves of this foliation as the orbits of the
complexified group. 

At this point it is convenient to change our point of view. So far we
have been looking at varying the complex structure on a symplectic
manifold, but in the end we are interested in K\"ahler metrics on a
complex manifold. If $F:M\to M$ is a diffeomorphism then the metric defined
by the pair $(\omega,F^*(J))$ is isometric to the one defined by
$( (F^{-1})^*\omega,J)$. If $\phi\in C^{\infty}_0(M)$, then the
infinitesimal action of the vector field $-JX_\phi$ on $\omega$ is
\[ \mathcal{L}_{-JX_\phi}\omega = -dJd\phi=2i\partial
\overline{\partial}\phi.\]
This shows that at least formally
the orbits of the complexified group can be identified
with the space K\"ahler metrics in a fixed K\"ahler class if we keep the
complex structure fixed instead of the symplectic form by applying
diffeomorphisms.

While we cannot directly use the results of the finite dimensional theory
developed in the previous chapter to characterise K\"ahler classes which
admit extremal metrics, we can use that theory to guide us to some
extent. For example we can now reinterpret the Futaki invariant and the
Mabuchi functional in this framework.
Let us fix a complex structure $J\in \mathcal{J}$. The
stabiliser $\mathfrak{g}_J$ of $J$ in $\mathfrak{g}$ is the space of holomorphic
vector fields on $(M,J)$ which have holomorphy potentials (ie. the space
$\mathfrak{h}_1$ introduced in the previous section), and the inner
product of Futaki and Mabuchi is just the restriction of the $L^2$
product on $\mathfrak{g}$. We can now rewrite the definition of the
Futaki invariant as a functional 
\[\begin{aligned}
  \mathcal{F}:\mathfrak{g}_J&\to\mathbf{C}\\
  \alpha &\mapsto \langle\mu(J),\alpha\rangle,
\end{aligned}
\]
which we recognise to be the weight functional defined in
Section~\ref{sec:normsquared}. 
Similarly, the variation of the Mabuchi functional $\mathcal{M}$
at the metric defined
by $J$ can be written as 
\[\begin{aligned}
  d\mathcal{M}_J:\mathfrak{g}&\to\mathbf{C}\\
  \alpha&\mapsto -\langle\mu(J),\alpha\rangle,
\end{aligned}
\]
which is the same as the variation of the norm functional defined in
Section~\ref{sec:kempfness}. In analogy with the finite dimensional
situation, to test whether a K\"ahler class contains a cscK metric,
we need to look at the asymptotic rate of change of
$\mathcal{M}$ as we tend towards the boundary of the K\"ahler class. One
problem is to identify what this boundary is, and another is to compute
the asymptotics of the Mabuchi functional. In the next chapter we
introduce the tools used to study this problem algebro-geometrically. 

\chapter{Stability of varieties}\label{chap:stabvar}

In the previous chapter we have seen that the Calabi functional can be
interpreted as the norm squared of a moment map, so if we apply the
results of Chapter~\ref{chap:finiteGIT} 
at least on a formal level then we expect that a K\"ahler class admits
an extremal metric if and only if it satisfies some kind of 
stability condition. In this chapter we make this more precise. In
Section~\ref{sec:Kstab} we introduce the notion of K-polystability. A
preliminary version was defined by Tian in \cite{Tian97} aiming to 
characterise
Fano varieties which admit K\"ahler-Einstein metrics. The version given
here is due to Donaldson~\cite{Don02} and was conjectured to be
equivalent to the existence of a constant scalar curvature K\"ahler
metric. We also define relative K-polystability due to the author
(see~\cite{GSz04}) which
was conjectured to characterise K\"ahler classes containing extremal
metrics. These conjectures are now likely to be false, since by an
example in~\cite{ACGT3} one might need to consider test-configurations
which are not algebraic.  In Section~\ref{sec:uniformKstab} we suggest a
way of strengthening the definition of K-stability to what we call
uniform K-stability to address this problem.

It is desirable to generalise the notion of K-stability to pairs
$(X,D)$ where $X$ is a polarised variety, and $D$ a divisor. If $D$ is a
smooth divisor then this would give a criterion to decide when
$X\setminus D$ admits a complete extremal metric, and in general
one expects that
an unstable variety breaks up into such stable pairs. An
example of this is given in Section~\ref{sec:extremallwr}. We propose a
notion of 
K-stability for pairs in Section~\ref{sec:kstabpair} and then we will
consider an example computation on a ruled surface in
Section~\ref{sec:ruledexample}. We will compare the results we obtain in
this case with the explicit construction of extremal metrics in
Section~\ref{sec:extremalruled}.

A fairly simple way to prove that a variety which admits a cscK metric is
K-semistable was given in \cite{Don05}. Here Donaldson showed that a
destabilising test-configuration gives a lower bound on the Calabi
functional. We show that this can be extended to the case of extremal
metrics in Section~\ref{sec:lwrcalabi}, which in particular proves that
a polarised variety that admits an extremal metric is relatively
K-semistable.

\section{K-stability}\label{sec:Kstab}

We would like to motivate the definition of K-stability using the
moment map picture we described in Section~\ref{sec:scalmoment}. 
In Section~\ref{sec:stability} we saw that the stability of a point $x$ in
geometric invariant theory can be verified by looking at the orbits of
the point under one-parameter subgroups and evaluating a numerical
weight on the limiting point. When this is positive for all one-parameter
subgroups which do not fix $x$ then $x$ is polystable. The main problem
in applying this directly to our infinite dimensional setting is
understanding what the one-parameter subgroups are. Alternatively,
we need to make sense of the boundary of the space of K\"ahler
metrics in a fixed K\"ahler class. What we can do instead is to consider
algebro-geometric degenerations of our complex manifold into possibly very
singular schemes. These are the ``test-configurations'' that we will
define below, and they are analogous to the orbits of one-parameter
subgroups in the finite dimensional theory. In fact as was shown
in~\cite{PS06} such a test-configuration can be used to define a weak
geodesic ray of metrics in the K\"ahler class although we will not use
this. 

Given a test-configuration, we need 
to define the weight on the central fibre. 
We know that for trivial degenerations, which are induced
by holomorphic vector fields on the manifold, this weight has to be
the Futaki
invariant we defined in Section~\ref{sec:FMprelim}.

We first recall the definition of the generalised Futaki invariant from
Donaldson~\cite{Don02}.
Let $V$ be a polarised scheme of dimension $n$ with a very ample line
bundle $L$.  Let
$\alpha$ be a $\mathbf{C}^*$-action on $V$ with a lifting to
$L$.  This induces
a $\mathbf{C}^*$-action on the vector space of sections
$H^0(V,L^k)$ for
all integers $k\geqslant 1$. Let $d_k$ be the dimension of
$H^0(V,L^k)$,
and denote the infinitesimal generator of the
action by $A_k$. Denote by $w_k(\alpha)$ the weight of the action on the top
exterior power of $H^0(V,L^k)$, which is the same as the trace
$\mathrm{Tr}(A_k)$. Then $d_k$ and $w_k(\alpha)$ are polynomials
in $k$ of degree $n$ and $n+1$ respectively for $k$ sufficiently large,
so we can write  
\[
\begin{aligned}
d_k&=c_0k^n + c_1k^{n-1} + O(k^{n-2}),\\
w_k(\alpha)=\mathrm{Tr}(A_k)&= a_0k^{n+1} + a_1k^n + O(k^{n-1}).
\end{aligned}
\]
\begin{defn}
  The \emph{Futaki invariant} of the $\mathbf{C}^*$-action $\alpha$ on
  $(V,L)$ is defined to be 
  \[F(\alpha)=\frac{c_1a_0}{c_0}-a_1.\]
\end{defn}

The choice of lifting of $\alpha$ to the line
bundle is not unique, however $A_k$ is defined up to addition of a
scalar matrix. In fact if we embed $V$ into $\mathbf{P}^{d_1-1}$ using
sections of 
$L$, then lifting $\alpha$ is equivalent to
giving a $\mathbf{C}^*$-action on $\mathbf{C}^{d_1}$ which induces
$\alpha$ on $V$ in
$\mathbf{P}^{d_1-1}$. Since
the embedding by sections of a line bundle is not
contained in any hyperplane, two such $\mathbf{C}^*$-actions differ by
an action that acts trivially on $\mathbf{P}^{d_1-1}$ ie. one with a
constant weight, say $\lambda$. We obtain that for
another lifting, the sequence of matrices $A^\prime_k$ are related to
the $A_k$ by 
\[ A^\prime_k = A_k + k\lambda I,\]
where $I$ is the identity matrix. A simple computation now shows
that $F(\alpha)$ is independent of the lifting of $\alpha$ to
$L$. It is shown in~\cite{Don02} that when $V$ is smooth and
the $\mathbf{C}^*$-action is induced by a holomorphic vector field then
this generalised Futaki invariant coincides with the classical Futaki
invariant we defined in Section~\ref{sec:FMprelim} up to a scalar
multiple. More precisely we have
\begin{prop} \label{prop:difffutaki}
  Suppose $\omega$ is a K\"ahler metric in the class
  $2\pi c_1(L)$ and the $\mathbf{C}^*$-action $\alpha$ is generated
  by a vector field $X$ with holomorphy potential $f$. Then
  \[ 2\cdot (2\pi)^n F(\alpha) = -\int_V
  f(S(\omega)-\hat{S})\frac{\omega^n}{n!}.\]
\end{prop}

We next recall the notion of a test-configuration from~\cite{Don02}.
\begin{defn} 
  A \emph{test-configuration for $(V,L)$ of exponent $r$}
  consists of a $\mathbf{C}^*$-equivariant flat family of schemes
  $\pi:\mathcal{V}\to\mathbf{C}$ (where $\mathbf{C}^*$ acts on
  $\mathbf{C}$ by multiplication) and a $\mathbf{C}^*$-equivariant ample
  line bundle $\mathcal{L}$ over $\mathcal{V}$.  We require that the
  fibres $(\mathcal{V}_t,\mathcal{L}|_{\mathcal{V}_t})$ are isomorphic
  to $(V,L^r)$ for $t\not=0$, where $\mathcal{V}_t=\pi^{-1}(t)$. The
  test-configuration is called a \emph{product configuration} if
  $\mathcal{V} = V\times\mathbf{C}$. The Futaki invariant of the induced
  $\mathbf{C}^*$-action on
  $(\mathcal{V}_0,\mathcal{L}|_{\mathcal{V}_0})$ is called the Futaki
  invariant of the test-configuration. 
\end{defn}

With these definitions we can now define when a polarised variety is
K-polystable. 
\begin{defn}
  A polarised variety $(V,L)$ is \emph{K-polystable} if for all
  test-configurations the Futaki invariant is non-negative and
  is zero if and only if the test-configuration is
  a product configuration. 
\end{defn}

We would now like to define relative K-polystability, following the
definitions in Section~\ref{sec:normsquared}. This uses an inner product
on the Lie algebra of the compact group in the moment map picture, which
in our case is the $L^2$ product on $C^{\infty}_0(V)$ when $V$ is
smooth. However we want to compute this algebro-geometrically on the
central fibre of a test-configuration, so we define
an inner product on $\mathbf{C}^*$-actions on a polarised variety
$(V,L)$, which coincides with the $L^2$ product when $V$ is smooth. 
Note that $\mathbf{C}^*$-actions do not naturally form a
vector space and we are really defining an inner product on a
subspace of the Lie algebra of the automorphism group of $(V,L)$. 

Let $\alpha$ and $\beta$ be two
$\mathbf{C}^*$-actions on $V$ 
with liftings to $L$. If we denote the infinitesimal generators of
the actions on $H^0(V,L^k)$ by $A_k, B_k$, then
$\mathrm{Tr}(A_kB_k)$ is a polynomial of degree $n+2$ in $k$. 
\begin{defn} The
inner product $\langle\alpha,\beta\rangle$ is defined 
to be the leading coefficient in
\begin{eqnarray*}
  &&\mathrm{Tr}\left[
  \left(A_k-\frac{\mathrm{Tr}(A_k)}{d_k}I\right)\left(B_k-
  \frac{\mathrm{Tr}(B_k)}{d_k}I\right)
\right] = \\ &&\qquad = \mathrm{Tr}(A_kB_k) -
\frac{w_k(\alpha)w_k(\beta)}{d_k} =  
\langle\alpha,\beta\rangle k^{n+2} + O(k^{n+1})\qquad \mbox{for }k\gg1.
\end{eqnarray*}
\end{defn}
Like the Futaki invariant, 
this does not depend on the particular liftings of $\alpha$ and
$\beta$ to the line bundle since we are normalizing each $A_k$ and $B_k$
to have trace zero.

Let us see what this is when the variety is
smooth. In this case we can consider the algebra of holomorphic vector
fields on $V$ which lift to $L$. This is the Lie algebra
of a group of holomorphic automorphisms of $V$. Inside this group, let
$G$ be the complexification of a maximal compact subgroup $K$. Let
$\mathfrak{g},\mathfrak{k}$ be the Lie algebras of $G,K$. Denoting by
$\mathfrak{k}_\mathbf{Q}$ the elements in $\mathfrak{k}$
which generate circle subgroups, our inner product on
$\mathbf{C}^*$-actions gives an inner
product on $\mathfrak{k}_\mathbf{Q}$. Since this is a dense subalgebra of
$\mathfrak{k}$, the inner product extends to 
$\mathfrak{k}$ by continuity. We further extend
this inner product to $\mathfrak{g}$ by complexification
and compute it differential geometrically. This is analogous
to the proof of Proposition~\ref{prop:difffutaki}
in Donaldson~\cite{Don02}. Let us choose
a $K$-invariant K\"ahler metric $\omega$ in the class $2\pi c_1(L)$. 
Note that
$\mathfrak{g}$ is a space of holomorphic vector fields on $V$ which lift
to $L$. Let $v,w$ be
two holomorphic vector fields on $V$, with liftings $\hat{v}, \hat{w}$
to $L$. 
We can write  
\[ \hat{v} = \overline{v} + if\underline{t}, \qquad
\hat{w}=\overline{w}+ig\underline{t},\]

\noindent where $\overline{v}$ (respectively $\overline{w}$) is the horizontal
lift of $v$ (respectively $w$), $\underline{t}$ is the canonical vector
field on the 
total space of $L$ defined by the action of scalar
multiplication, and $f,g$ are smooth functions on $V$. As
in \cite{Don02} we have that 
\[ \overline{\partial}f = -(i_v(\omega))^{0,1}, \qquad
\overline{\partial}g = -(i_w(\omega))^{0,1}, \]

\noindent so in particular $f$ and $g$ are defined up to an additive
constant, and we can normalise them to have zero integral over $V$. We
would like to show that
\[ \langle v,w\rangle = (2\pi)^{-n}\int_V fg\frac{\omega^n}{n!}, \]

\noindent where we have assumed that $f,g$ have zero integral over $V$.
Making use of the identity $\langle v,w\rangle = \frac{1}{2}(\langle
v+w,v+w\rangle - \langle v,v\rangle - \langle w,w\rangle)$ it is enough
to show this when $v=w$. Furthermore, we can assume that $v$ generates a
circle action since $\mathfrak{k}_\mathbf{Q}$ is dense in
$\mathfrak{k}$.  

We can find the leading coefficients of $d_k,\mathrm{Tr}(A_k),
\mathrm{Tr}(A_kA_k)$ for this circle action using the
equivariant Riemann-Roch formula, in the 
same way as was done in~\cite{Don02}. We find that these leading
coefficients are given by
\[ (2\pi)^{-n}\int_V\frac{\omega^n}{n!},\quad (2\pi)^{-n}
\int_V f\frac{\omega^n}{n!},\quad
(2\pi)^{-n}\int_V f^2\frac{\omega^n}{n!},\]

\noindent respectively. If we normalise $f$ to have zero integral over
$V$, then we obtain the formula for the inner product that we were
after. 

To define relative K-polystability we also need to modify the
definition of a test-configuration slightly.  
\begin{defn}
  We say that the test-configuration $(\mathcal{V},\mathcal{L})$ for
  $(V,L)$ is \emph{compatible with a torus $T$
  of automorphisms of $(V,L)$}, if there is a torus action on
  $(\mathcal{V},\mathcal{L})$ which preserves the fibres of
  $\pi:\mathcal{V}\to\mathbf{C}$, commutes with the
  $\mathbf{C}^*$-action, and restricts to $T$ on
  $(\mathcal{V}_t,\mathcal{L}|_{\mathcal{V}_t})$ for $t\not=0$. 
\end{defn}

Fix a maximal torus of automorphisms of
$(V,L)$, and write $\chi$ for the $\mathbf{C}^*$-action induced by the
extremal vector field. This is defined as in Section~\ref{sec:FMprelim},
by requiring that $F(\alpha)=\langle\chi,\alpha\rangle$ for all
$\mathbf{C}^*$-actions $\alpha$ in the torus. Because of the difference
between the algebraic and differential-geometric definitions of the
Futaki invariant and inner product, this is half of the differential
geometric extremal vector field. 
With these preliminaries we can state the definition of relative
K-polystability.
\begin{defn}\label{def:kstable}
  A polarised variety $(V,L)$ is \emph{K-semistable relative to a
  maximal torus $T$ of automorphisms} if 
  \begin{equation}\label{eq:relativefutaki} 
    F_{\tilde{\chi}}(\tilde{\alpha}):=F(\tilde{\alpha}) -
    \langle\tilde{\chi},\tilde{\alpha}\rangle\geq 0
  \end{equation}
for all
test-configurations compatible with $T$. Here we denote by
$\tilde{\alpha}$
and $\tilde{\chi}$ the $\mathbf{C}^*$-actions induced on the central fibre of
the test-configuration. 
The variety is relatively K-polystable if in addition equality holds only if
the test-configuration is a product configuration.
\end{defn}

\subsection{Uniform K-stability}\label{sec:uniformKstab}

In~\cite{GSz04} we conjectured that a polarised variety is 
K-polystable
relative to a maximal torus of automorphisms if and only if it admits an
extremal metric, analogously to the conjecture in~\cite{Don02} in the
cscK case. In Section~\ref{sec:lwrcalabi} we will prove 
that a variety that admits an extremal metric
is relatively K-semistable. 
As we mentioned in the introduction,
an example of~\cite{ACGT3} shows that the converse statement is likely
to be false and the
conjectures need to be refined. Their example is a ruled manifold which
is destabilised by a test-configuration with a non-algebraic
polarisation. This is possible in the framework of slope-stability
(see~\cite{RT06}).

A natural approach to remedy this situation
is to strengthen the definition of K-polystability to
\emph{uniform K-polystability} as follows, extending the notion of the
modulus of stability (see Section~\ref{sec:modulus}) to this setting. We
choose a maximal torus of 
automorphisms $T$ of $(V,L)$, and define $(V,L)$ to be uniformly
K-polystable if there exists a positive
constant $\lambda>0$ such that for all test-configurations compatible
with $T$,
\[ F(\alpha) \geq \lambda\Vert\alpha - \pi(\alpha)\Vert. \]
Here $\alpha$ is the $\mathbf{C}^*$-action induced on the central
fibre and $\pi(\alpha)$ is the orthogonal 
projection of $\alpha$ onto $T$. 
Then $\alpha-\pi(\alpha)$ might not generate a
$\mathbf{C}^*$-action, but we can still define its norm. 

In finite dimensional GIT such a $\lambda$ exists for all polystable
points, but in the case
of varieties this is no longer necessarily the case. In
section~\ref{sec:toricsurface} we will show
that in the case of toric surfaces K-polystability implies uniform
K-polystability, however for higher dimensional varieties we need to
change the definition of the norm to an $L^{\frac{n}{n-1}}$ 
norm instead of the $L^2$
norm we used above if we want the definition to make sense ($n$ is the
complex dimension of the variety). 
We can define the $L^p$-norm of a $\mathbf{C}^*$-action as we have
defined the $L^2$-norm, looking at the asymptotics of
$\mathrm{Tr}(|A_k|^p)$ where the $A_k$ are endomorphisms induced on
$H^0(V,L^k)$ by the $\mathbf{C}^*$-action as before. However unless $p$ is an
even integer, the equivariant
Riemann-Roch formula can no longer be used to show that this coincides with
the $L^p$-norm of a Hamiltonian. In any case it is tempting to
conjecture that uniform K-polystability is the correct condition
characterising the existence of cscK metrics. 
The stronger assumption of uniform K-polystability should make it easier
to make analytic deductions. In particular we saw in
Section~\ref{sec:modulus} that control of the modulus of stability can
be used to control the first eigenvalue of an operator $\sigma^*\sigma$.
In the infinite dimensional setting this is the Lichnerowicz
operator and controlling its first eigenvalue is crucial to the analysis
in trying to prove the existence of an extremal metric. The bound we
gave in Section~\ref{sec:modulus} is not good enough for this purpose,
it is just an indication of what might be possible. 

\subsection{K-stability of a pair $(V,D)$} \label{sec:kstabpair}

In this subsection we propose a definition of K-stability for a pair
$(V,D)$ consisting of a smooth polarised variety $(V,L)$ and a smooth divisor
$D\subset V$. The aim is to find a stability condition for the existence
of a complete extremal metric on $V\setminus D$ which is asymptotically
hyperbolic near $D$. A well-known case is the polarisation $K_V+D$, where
$K_V$ is the canonical bundle of $V$ and we assume that $K_V+D$ is
ample. In this case it was shown by Cheng-Yau~\cite{CY80} (see
also Tian-Yau~\cite{TY87}) that a complete K\"ahler-Einstein metric exists on
$V\setminus D$, which is asymptotically hyperbolic near
$D$. 

It is not yet clear what the precise class of
metrics is that one should consider, but let us use the following as a
preliminary definition:

\begin{defn}\label{defn:cuspsing}
  A complete K\"ahler metric $g$ on $V\setminus D$ is
  \emph{asymptotically hyperbolic near $D$} if near $D$ it is asymptotic
  to a metric of the form
  \[ g_0 = K\cdot\frac{|dz|^2}{(|z|\log |z|)^2}+h,\]
  where $K$ is a smooth positive function on $V$, 
  the symmetric 2-tensor $h$ is
  a smooth extension of a metric on $D$, and $z$ is a local defining
  holomorphic function for $D$.
  By $g$ being asymptotic to $g_0$ near $D$, we mean that for all
  $i\in\mathbf{N}$,
  \[ \lim_{z\to 0}\big\Vert \nabla^i_{(g_0)} (g-g_0)\big\Vert_{(g_0)}=0.\]
\end{defn}

Note that $|dz|^2/(|z|\log |z|)^2$ is the standard hyperbolic cusp
metric on the punctured disk up to a scalar factor. 
The function $K$ is necessary because we do not want to prescribe the
curvature near $D$ in the normal directions to $D$. 
For such a metric on $V\setminus D$ we define a
K\"ahler class in $H^2(V)$ as follows. The K\"ahler form $\omega$
corresponding to $g$ defines an
$L^2$-cohomology class in $H^2_{L^2}(V\setminus D,g)$. Since $g$ is
quasi-isometric to the fibred cusp metrics of~\cite{HHM04}, according to
Corollary 2 in that paper, this $L^2$-cohomology group 
is naturally isomorphic to the de Rham cohomology $H^2(V)$. The class in
$H^2(V)$ defined in this way by $\omega$ is the K\"ahler class of our
metric. For example for a metric on $\mathbf{P}^1$ which is
asymptotically hyperbolic near a point, the K\"ahler class is simply
given by the total area as an element in $H^2(\mathbf{P}^1)\cong
\mathbf{R}$. 

We incorporate the divisor $D$ into the definition
of a test-configuration as follows.

\begin{defn} A test-configuration for $(V,D,L)$ is a test-configuration
  $(\mathcal{V},\mathcal{L})$ for $(V,L)$ with a
  $\mathbf{C}^*$ invariant Cartier divisor
  $\mathcal{D}\subset\mathcal{V}$ which is flat over $\mathbf{C}$ and
  restricts to $D$ on the non-zero fibres.  
\end{defn}

The central fibre of such a test-configuration is a polarised scheme
$(V_0,L_0)$ with a $\mathbf{C}^*$-action and a divisor $D_0\subset V_0$
fixed by the $\mathbf{C}^*$-action. 
We define a modification of the Futaki invariant for this situation.

Let $(V,L)$ be a polarised scheme with a $\mathbf{C}^*$-action $\alpha$ and
$D\subset V$ a divisor fixed by the $\mathbf{C}^*$-action. Let us write 
\[ H^0(V,L^k\otimes\mathcal{O}(-D))\subset H^0(V,L^k) \]
for the sections which vanish along $D$. The inclusion is induced by
a section of $\mathcal{O}(D)$ which vanishes along $D$.
The assumption that $D$ is invariant
under the $\mathbf{C}^*$-action means that this subspace is preserved by
the action. As before, let us write $d_k, w_k$ for the dimension of
$H^0(V,L^k)$ and for the total weight of the action on this space. Let
us also write $\tilde{d}_k, \tilde{w}_k$ for the dimension and weight of
the action on the space $H^0(V,L^k\otimes\mathcal{O}(-D))$. Define
$c_0,c_1,a_0,a_1$ by
\begin{equation}\label{eq:cuspfutaki} \begin{split}
  \frac{d_k+\tilde{d}_k}{2} &= c_0k^n+c_1k^{n-1}+O(k^{n-2}),\\
  \frac{w_k+\tilde{w}_k}{2} &= a_0k^{n+1}+a_1k^n+O(k^{n-1}).
\end{split} \end{equation}
The Futaki invariant of the $\mathbf{C}^*$-action $\alpha$
on the pair $(V,D)$
is then 
\[ F(\alpha) = \displaystyle{\frac{c_1}{c_0}a_0-a_1}.\] 

\noindent Let us also define $\alpha_1, \alpha_2$ by
\[ \dim H^0(D,L^k\big\vert_D) = \alpha_1 k^{n-1}+\alpha_2 k^{n-2}+
O(k^{n-3}).\]

\begin{defn}
We say that the triple $(V,D,L)$, where $L$ is an ample line bundle over
$V$ and $D\subset V$ is a divisor, is K-polystable if the Futaki
invariant of every test-configuration for $(V,D,L)$ is 
non-negative and is zero only for product configurations. In addition we
require that $c_1/c_0<\alpha_2/\alpha_1$. 
\end{defn}

Let us briefly explain the last condition. In
Section~\ref{sec:extremalruled} we will construct complete extremal
metrics on the complement of a divisor on a ruled surface, and we will
find that for a range of polarisations $m<k_2$ (where $m$ parametrises
the polarisation and $k_2$ is a constant - see
Section~\ref{sec:extremalruled} for details) we obtain metrics which are
asymptotically hyperbolic as in Definition~\ref{defn:cuspsing}. For
$m=k_2$ we also obtain an extremal metric but it no longer has the
asymptotic behaviour prescribed in Definition~\ref{defn:cuspsing} but
instead the fibre metrics behave like
\[ \frac{|dz|^2}{|z|^2(\log |z|)^{3/2}} \]
near the divisor. The non-degeneracy condition $c_1/c_0 <
\alpha_2/\alpha_1$ is aimed to rule out this possibility. The condition
arises when looking at deformation to the normal cone of the divisor
$D$ (see Section~\ref{sec:ruledexample} for the definition of
deformation to the normal cone or~\cite{RT06} for more details). This gives a
family of test-configurations parametrised by $c\in(0,\epsilon)$ where
$\epsilon$ is a small positive number. It
turns out that the Futaki invariant $F(c)$ of these test-configurations
satisfies $F(0)=F'(0)=0$ and our non-degeneracy condition is $F''(0)>0$. 

It is straightforward to extend the notion of relative K-polystability
to pairs as well. We only consider automorphisms of $(V,L)$ which fix
$D$ (but can induce a nontrivial automorphism of $D$), 
and define the extremal $\mathbf{C}^*$-action in this group. The
modified Futaki invariant is defined as before. The non-degeneracy
condition is defined as above using deformation to the normal cone of
$D$, but with the modified Futaki invariant. 

We conjecture that if $D$ is a smooth divisor then $(V,D,L)$ is 
relatively K-polystable (with a positive modulus of stability)
if and only if there exists a complete extremal
metric on $V\setminus D$ in the cohomology class $c_1(L)$,
which is asymptotically hyperbolic near
$D$. 

Our definition of a
test-configuration was chosen because it seems natural, and it is
satisfactory for the example that we compute in the next section. The
definition of the Futaki invariant is motivated by the calculations
involved in the explicit
construction of extremal metrics in Section~\ref{sec:extremalruled}.
Those calculations also suggest
that different combinations of $d_k$ and $\tilde{d}_k$ in
Equation~\ref{eq:cuspfutaki}
could be used to characterise incomplete metrics with edge singularities
along $D$ with various angles. 

\section{Relative K-polystability of a ruled surface}
\label{sec:ruledexample}

The aim of this section is to work out the stability
criterion in a special case. Let $\Sigma$ be a genus two curve, and
$\mathcal{M}$ a line
bundle on it with degree one (the calculations also work for genus
greater than two and a line bundle with degree greater than one). 
Define $X$ to be the ruled surface
$\mathbf{P}(\mathcal{O}\oplus\mathcal{M})$ over $\Sigma$.
T\o{}nnesen-Friedman~\cite{TF97} constructed a family of extremal metrics on
$X$, which does not exhaust the entire K\"ahler cone (see also
Section~\ref{sec:extremalruled}). We will show that
$X$ is K-unstable (relative to a maximal torus of automorphisms) for the remaining
polarisations (it was shown in~\cite{ACGT3} that $X$ does not admit an
extremal metric for these unstable polarisations). We will
also look at K-stability of the pairs $(X,S_0)$ and
$(X,S_\infty)$ where $S_0$ and $S_\infty$ are the zero and infinity
sections of $X$.  

Since there are no non-zero holomorphic vector fields on $\Sigma$, a
holomorphic vector field on $X$ must preserve the fibres. Thus, the
holomorphic vector fields on $X$ are given by sections of
$\mathrm{End}_0(\mathcal{O}\oplus\mathcal{M})$. Here $\mathrm{End}_0$
means endomorphisms with trace zero. The vector field given by the
matrix
\[ \left( \begin{array}{cc}
              -1 & 0 \\
	      0 & 1 
	  \end{array} \right)
\]

\noindent generates a $\mathbf{C}^*$-action $\beta$, which is a maximal
torus of automorphisms (see Maruyama~\cite{Mar71} for
proofs). Therefore this must be a multiple of the extremal vector field,
which is then given by $\chi =
\frac{F(\beta)}{\langle\beta,\beta\rangle}\beta$. 

The destabilising test-configuration is an example of deformation to
the normal cone of a subvariety, studied by Ross and Thomas~\cite{RT04}
(see also Section~\ref{sec:metricdegen}),
except we need to take into account the extremal $\mathbf{C}^*$-action
as well. We consider the polarisation $L=C+mS_0$ where $C$ is the
divisor given by a fibre, $S_0$ is the zero section (ie. the image of
$\mathcal{O}\oplus\{0\}$ in $X$, so that $S_0^2=1$) 
and $m$ is a positive constant. We
denote by $S_\infty$ the infinity section, so that $S_\infty=S_0-C$. 
Note that $\beta$ fixes $S_\infty$ and acts on the normal
bundle of $S_\infty$ with 
weight 1. We make no distinction between divisors and their associated
line bundles, and use the multiplicative and additive notations
interchangeably, so for example $L^k=kC+mkS_0$ for an integer $k$.

The deformation to the normal cone of $S_\infty$ is given by the blowup
\[\mathcal{X}:=\widetilde{X\times\mathbf{C}}\xrightarrow{\pi}X\times\mathbf{C}
\]
in the subvariety $S_\infty\times\{0\}$. Denoting the
exceptional divisor by $E$, the line bundle $\mathcal{L}_c=\pi^*L-cE$ is
ample for $c\in(0,m)$. For these values of $x$ we therefore 
obtain a test-configuration
$(\mathcal{X},\mathcal{L}_c)$ with the $\mathbf{C}^*$ action induced by
$\pi$ from the product of the trivial action on $X$ and the usual
multiplication on $\mathbf{C}$. Denote the restriction of this
$\mathbf{C}^*$-action to the central fibre $(X_0,L_0)$ by $\alpha$.

Since the $\mathbf{C}^*$-action $\beta$ fixes $S_\infty$ we obtain
another action on the test-configuration, induced by $\pi$ from the
product of the $\mathbf{C}^*$-action $\beta$ on $X$ and the trivial
action on $\mathbf{C}$. Let us call the induced action on the central
fibre $\beta$ as well. We wish to calculate $F_\chi(\alpha)$ where
$\chi$ is a scalar multiple of $\beta$ as above. 
For this we need the weight decomposition of the space $H^0(X_0,L^k_0)$.
According to~\cite{RT04} we have 
\[ H^0(X_0,L^k_0) = H^0_X(kL-ckS_\infty)\oplus \bigoplus_{j=1}^{ck}t^j
\frac{H^0_X(kL-(ck-j)S_\infty)}{H^0_X(kL-(ck-j+1)S_\infty)}, \]
for $k$ large, 
with $t$ being the standard coordinate on $\mathbf{C}$. This gives the
weight decomposition for the action of $\alpha$. For the action $\beta$,
we need to further decompose $H^0_X(kL-ckS_\infty)$ into weight spaces
as follows:
\[ H^0_X(kL-ckS_\infty)=H^0_X(kL-mkS_\infty)\oplus\bigoplus_{i=1}^{mk-ck}
\frac{H^0_X(kL-(mk-i)S_\infty)}{H^0_X(kL-(mk-i+1)S_\infty)},\]
for $k$ large. This holds because of the following cohomology vanishing
lemma.
\begin{lem} $H^1(X,kC+lS_0)=0$ for $k\gg 0$ and $l\geqslant 0$. 
\end{lem}
\begin{proof}
  Let $f:X\to\Sigma$ be the projection map. Since
  \[\mathcal{O}_X(C)=f^*(\mathcal{O}_\Sigma(P))\] where $P\in\Sigma$ is
  a point, we have
  \[ R^1f_*\mathcal{O}_X(kC+lS_0)=\mathcal{O}_\Sigma(kP)\otimes R^1
  f_*\mathcal{O}_X(lS_0).\]
  The restriction of $\mathcal{O}_X(lS_0)$ to a fibre is
  $\mathcal{O}_{\mathbf{P}^1}(l)$ which for $l\geq 0$ has
  trivial $H^1$. This shows that $R^1 f_*\mathcal{O}_X(kC+lS_0)=0$
  for $l\geq 0$. The Leray-Serre spectral sequence now shows
  that $H^1(X,kC+lS_0)=H^1(\Sigma, \mathcal{O}_\Sigma(kP)\otimes
  f_*\mathcal{O}_X(lS_0))$. Since $\mathcal{M}$ has degree one, 
  each summand in
  \[ f_*\mathcal{O}_X(lS_0)=\bigoplus_{i=0}^l \mathcal{M}^{\otimes i} \]
  has non-negative degree, so for $k$ large (in fact for $k>2$) we have
  \[ H^1(\Sigma, \mathcal{O}_\Sigma(kP)\otimes
  f_*\mathcal{O}_X(lS_0))=0\]
  by Serre duality.
  This completes the proof. 
\end{proof}
In sum we obtain the decomposition
\begin{equation}\label{eq:decomp}
  \begin{split}
H^0(X_0,L_0^k)=&H^0_X(kL-mkS_\infty)\oplus\bigoplus_{i=1}^{mk-ck}
\frac{H^0_X(kL-(mk-i)S_\infty)}{H^0_X(kL-(mk-i+1)S_\infty)}\oplus\\
&\bigoplus_{j=1}^{ck}t^j
\frac{H^0_X(kL-(ck-j)S_\infty)}{H^0_X(kL-(ck-j+1)S_\infty)}.
 \end{split}
\end{equation}
According to~\cite{RT04} $\alpha$ acts with weight $-1$ on $t$
that is, it acts with weight $-j$ on the summand of index $j$ above.
Also, $\beta$ acts on
\[\frac{H^0_X(kL-lS_\infty)}{H^0_X(kL-(l+1)S_\infty)}\]
with weight $l$, plus perhaps a constant independent of $l$
which we can neglect, since the matrices are normalized to have trace zero
in the formula for the modified Futaki invariant.
The dimension of this space is $k+l-1$ by the Riemann-Roch theorem.
Writing $A_k, B_k$ for the infinitesimal generators of the actions
$\alpha$ and $\beta$ on $H^0(X_0,L_0^k)$ and $d_k$ for the dimension of
this space, we can now compute
\begin{equation}\label{eq:futakicomp}
  \begin{split}
d_k &= \frac{m^2+2m}{2}k^2 + \frac{2-m}{2}k + O(1),\\
\mathrm{Tr}(A_k) &= -\frac{c^3+3c^2}{6}k^3 + \frac{c^2-c}{2}k^2 + O(k),\\
\mathrm{Tr}(B_k) &= \frac{2m^3 + 3m^2}{6}k^3 + \frac{m}{2}k^2 + O(k),\\
\mathrm{Tr}(A_kB_k) &= -\frac{c^4+2c^3}{12}k^4 + O(k^3),\\
\mathrm{Tr}(B_kB_k) &= \frac{3m^4+4m^3}{12}k^4 + O(k^3).
\end{split}
\end{equation}

\noindent Using these, we can compute 
\[ F_\chi(\alpha) = F(\alpha)-\langle\alpha,\chi\rangle = 
F(\alpha)-\frac{\langle\alpha,\beta\rangle}{\langle
\beta,\beta\rangle}F(\beta).\]
We obtain
\[
F_\chi(\alpha)=\frac{c(m-c)(m+2)}{4(m^2+6m+6)}\Big[(2m+2)c^2-(m^2-4m-6)c
+m^2+6m+6\Big]. 
\]

\noindent If $F_\chi(\alpha)\leqslant0$
for a rational 
$c$ between 0 and $m$, then the variety is K-unstable (relative to
a maximal torus of automorphisms). We see that the variety is relatively
K-unstable for $m\geq k_1\cong 18.889$, where $k_1$ is the only positive
real root of the quartic $m^4-16m^3-52m^2-48m-12$. 

Let us now look at relative K-stability of the pair $(X,S_\infty)$ and use
the same test-configuration as above, ie. deformation to the normal cone
of $S_\infty$ with parameter $c\in(0,m)$. We
have the same decomposition of $H^0(X_0,L_0^k)$ as in
Equation~\ref{eq:decomp}, and we can 
identify the quotient of $H^0(X_0,L_0^k)$ by the space of sections
which vanish along $S_\infty$ with 
$H^0(kL)/H^0(kL-S_\infty)$. We need to subtract half of the
contribution of this space from the formulae in~\ref{eq:futakicomp} to
calculate the modified Futaki invariant for the pair. The dimension of
this space is
$k-1$, the weight of $\alpha$ on it is $-ck$ and the weight of $\beta$ is
$0$. The new formulae are therefore 
\begin{equation*}
  \begin{split}
    d^\prime_k &= \frac{m^2+2m}{2}k^2+\frac{1-m}{2}k+ O(1),\\
    \mathrm{Tr}(A_k^\prime) &= -\frac{c^3+3c^2}{6}k^3 + \frac{c^2}{2}k^2
    + O(k),
  \end{split}
\end{equation*}
while the other expansions remain unchanged.
The new modified Futaki invariant is then
\[ F_\chi(\alpha) = 
\frac{c^2(m-c)}{2m^2(m^2+6m+6)}\left[ c(2m^2+4m+3)-m^3+3m^2+9m+6\right]. \]
The pair $(X,S_\infty)$ is relatively K-unstable if $F_\chi(\alpha)\leqslant
0$ for some rational $c\in(0,m)$ or if the order of vanishing of
$F_\chi(\alpha)$ at $c=0$ is greater than 2. This happens if $m\geqslant 
k_2\cong 5.0275$, where
$k_2$ is the only positive real root of the cubic $m^3-3m^2-9m-6$. 

The calculation for the pair $(X,S_0)$ is essentially the same,
except we use deformation to the normal cone of $S_0$ in that case. 
These results should be compared to the results in
Section~\ref{sec:extremalruled} where we construct extremal metrics in
the K\"ahler classes which are not destabilised by the
test-configurations we considered here. 

\section{Lower bound on the Calabi functional}\label{sec:lwrcalabi}

In~\cite{Don05} Donaldson showed that a destabilising test-configuration
gives a lower bound for the Calabi functional. The precise statement is
the following.
\begin{thm} 
\label{thm:lowerCalabi}
  Let $\alpha$ be a destabilising test-configuration with
  Futaki invariant $F(\alpha)<0$. Then for any metric $\omega$ in the
  class of our polarisation,
  \[ \Vert
  S(\omega)-\hat{S}\Vert^2_{L^2}\geq
  4\cdot(2\pi)^n\frac{F(\alpha)^2}{\Vert\alpha\Vert^2}.
  \]
\end{thm}
The constant $4\cdot (2\pi)^n$ arises from the difference between the
differential-geometric and algebro-geometric Futaki invariants and inner
products. 
This is analgous to the finite dimensional result,
Theorem~\ref{thm:destablwrbound}. 
In the same way as we did there, we can extend the result to the case of
extremal metrics. We simply need to modify the test-configuration in
such a way as to obtain the optimal inequality. This should be compared
with Inequality~\ref{eq:celwr}.
\begin{thm}\label{thm:extremallwr}
Let $T$ be a maximal torus of automorphisms of $(X,L)$ with
corresponding extremal vector field $\chi$. Let $\mathcal{X}$ be a
test-configuration compatible with $T$ such that $F_\chi(\alpha)<0$ for
the $\mathbf{C}^*$-action $\alpha$ induced on the central fibre.
Then for any metric $\omega\in 2\pi c_1(L)$,
\[ \Vert S(\omega)-\hat{S}\Vert^2_{L^2}\geq2\cdot(2\pi)^n
\frac{F_\chi(\alpha)^2}{
\Vert\alpha\Vert^2} + \Vert\chi\Vert_{L_2}^2. \]
\end{thm}
\noindent Here $\Vert\chi\Vert_{L^2}$ 
is the differential-geometric norm of the differential geometric
extremal vector field, which is $2\cdot (2\pi)^{n/2}$ times the
algebraic norm of the algebraic extremal vector field. We will write
$\Vert\chi\Vert$ without the $L^2$ subscript for the latter, hoping that
it does not cause confusion. 
\begin{proof}
  Since the test-configuration is compatible with $T$,
  there is a $\mathbf{C}^*$-action $\tilde{\chi}$ on $\mathcal{X}$
  fixing the base $\mathbf{C}$,
  which restricts to $\chi$ on the nonzero fibres. Write
  $\tilde{\alpha}$ for the $\mathbf{C}^*$-action on $\mathcal{X}$
  induced by the test-configuration. We can modify
  the test-configuration by multiplying $\tilde{\alpha}$ by a multiple
  of $\tilde{\chi}$. The Futaki invariant of the new test-configuration
  will be $F(\alpha)+l F(\chi)$ for some integer $l$, 
  where $F(\chi)$ is the Futaki
  invariant of the vector field $\chi$ on $X$. Note that
  $F(\chi)=\Vert\chi\Vert^2$ by the definition of the extremal vector
  field.
  
  We can also pull back the test-configuration $\mathcal{X}$ under  
  a map 
  \[\begin{split} \mathbf{C}&\to\mathbf{C}\\
    \lambda&\mapsto\lambda^k,\end{split} \]
  for positive integers $k$,
  which changes the Futaki
  invariant to $kF(\alpha)$. This means that we can
  construct a test-configuration with Futaki invariant equal to
  $kF(\alpha) + lF(\chi)$ for any integers $k,l$ with $k>0$. 
  Since we are only interested in the quotient of the Futaki invariant
  by the norm of the $\mathbf{C}^*$-action, we will assume that we can
  also have rational $k$ and $l$. 
  
  For irrational $k,l$ we can still
  define $\Vert k\alpha+l\chi\Vert$ and $F(k\alpha+l\chi)$ by
  continuity. We follow the proof of Theorem~\ref{thm:destablwrbound}.
  Let $\overline{\alpha}$ be the component of $\alpha$ orthogonal to
  $\chi$, ie. $\overline{\alpha}=\alpha - \lambda\chi$ for some
  $\lambda$ such that $\langle\overline{\alpha},\chi\rangle=0$. By our
  assumption $F(\overline{\alpha})=F_\chi(\alpha)$ is negative, so we can
  choose a positive constant $\mu$ such that 
  $F(\mu\overline{\alpha})=-\Vert\mu
  \overline{\alpha}\Vert^2$. Now define $\gamma = \mu\overline{\alpha}
  - \chi$. We have
  \[ 
    F(\gamma) = -\Vert\mu\overline{\alpha}\Vert^2 - \Vert\chi\Vert^2 =
    - \Vert\gamma\Vert^2, \]
  which is negative, and
  \[ \frac{F(\gamma)^2}{\Vert\gamma\Vert^2} =
  \frac{F_\chi(\alpha)^2}{\Vert\overline{\alpha}\Vert^2} +
  \Vert\chi\Vert^2.
  \]
  We can approximate $\gamma$ with a rational linear combination
  $k\alpha+l\chi$, and apply Theorem~\ref{thm:lowerCalabi} to obtain
 \[ \Vert S(\omega)-\hat{S}\Vert^2_{L^2}\geq 4\cdot (2\pi)^n
    \frac{F_\chi(\alpha)^2}{
    \Vert\overline\alpha\Vert^2} + 4\cdot(2\pi)^n\Vert\chi\Vert^2. \]
 Since $\Vert\overline{\alpha}\Vert\leq\Vert\alpha\Vert$, we get the
 required result.
\end{proof}
Note that this proves that if $(X,L)$ admits an extremal metric then it
is relatively K-semistable, since the extremal metric would satisfy
$\Vert S(\omega)-\hat{S}\Vert_{L^2} = \Vert\chi\Vert_{L^2}$ (see
Section~\ref{sec:FMprelim}).

\chapter{Toric varieties}\label{chap:toric}

In \cite{Don02} Donaldson developed the theory of K-stability for
toric varieties. 
The main result is that on a toric surface the Mabuchi
functional on torus invariant metrics 
is bounded from below if and only if the surface is K-polystable with respect
to toric degenerations. It remains to be seen whether this implies the
existence of a cscK metric on the variety, but much progress on this has
been made in~\cite{Don05_1}, giving interior a priori estimates for the
PDE in the case of toric surfaces. 

In Section~\ref{sec:kstabtoric} we recall the construction of
toric test-configurations from~\cite{Don02}, 
but we generalise it to test-configurations of bundles
of toric varieties. We will use this generalisation in the next chapter. 

In Section~\ref{sec:toricsurface} we concentrate on
toric surfaces. We first 
show that a K-polystable toric surface is uniformly K-polystable. 
Using the results of \cite{Don02} this 
boils down to the statement that for a
positive convex function on a polygon the integral on the boundary
controls the $L^2$-norm on the interior. We will see that for this to
hold in higher dimensions we need to use the $L^\frac{n}{n-1}$-norm instead
of the $L^2$-norm, where $n$ is the dimension. We then give an alternative
proof of a result in~\cite{Don02}, using the notion of measure
majorisation from convex geometry. Finally we use the same technique
to prove that a
semistable polygon can be decomposed into stable subpolygons.  This is
analogous to the Jordan-H\"older filtration of a semistable vector bundle. 

\section{K-stability of toric varieties}\label{sec:kstabtoric}
Let $\Delta\subset\mathbf{R}^n$ be the moment polytope of a smooth
polarised toric variety.
The polytope is defined by a finite number of
linear inequalities $h_k(x)\geq c_k$, where the $h_k$ are linear maps
from $\mathbf{R}^n$ to $\mathbf{R}$ which induce primitive maps from the
integer lattice $\mathbf{Z}^n$ to $\mathbf{Z}$. Let $d\mu$ be the
standard Euclidean volume form on $\Delta$, and define a measure
$d\sigma$ on
the boundary $\partial \Delta$  as follows. On the face defined by the
equation $h_r(x)=c_r$ we let $d\sigma$ be the constant $(n-1)$-form such
that $dh_r\wedge d\sigma$ is, up to sign, $d\mu$. 

Recall the following result from~\cite{GS06}
\begin{thm} \label{thm:trap}
  Let $Q:\Delta\to\mathbf{R}$ be a continuous function. We
  have for $k\gg 1$, 
  \[ \sum_{\alpha\in k\Delta\cap\mathbf{Z}^n} Q(\alpha) = k^n\int_\Delta Q\,d\mu +
  \frac{k^{n-1}}{2}\int_{\partial\Delta} Q\,d\sigma + O(k^{n-2}).
  \]
\end{thm}

\subsection*{Test-configurations for toric bundles}
We now construct test-configurations for toric bundles, extending the
construction for toric varieties in~\cite{Don02}. 
Let us first define toric bundles. The data is a principal
$T=(\mathbf{C}^*)^n$-bundle $P\to M$ over a projective variety $M$ of
dimension $m$, and
an $n$-dimensional 
polarised toric variety $(V,\mathcal{O}_V(1))$ with corresponding
moment polytope $\Delta\subset\mathfrak{t}^*$. Define a bundle of toric
varieties $\pi: X\to M$ by
\[ X = P\times_T V.\]
Let $L_M\to M$ be an ample line bundle and define a line bundle $L$ over
$X$ by
\[ L = \pi^* L_M\otimes \left(P\times_T\mathcal{O}_V(1)\right). \]
Let us assume that it is ample for our choice of data, so that the pair
$(X,L)$ is a polarised variety. 

For each $\alpha\in\mathfrak{t}^*\cap\mathbf{Z}^n$ we define a line
bundle $F_\alpha$ over $M$ with transition functions induced by the map 
\[ \begin{aligned} T &\to  \mathbf{C}^*\\
  \exp(\xi) &\mapsto \exp(i\alpha(\xi)),\quad\text{for
  }\xi\in\mathfrak{t}. \end{aligned}
\]
The pushforward $\pi_* L^k$ is then given by
\[ \pi_* L^k=L_M^k\otimes \bigoplus_{\alpha\in k\Delta\cap\mathbf{Z}^n}
F_\alpha.\]
Let us define the functions $Q_1,Q_2:\Delta\to\mathbf{R}$ by
\[\begin{aligned} Q_1(\alpha)&=c_1(L_M\otimes F_\alpha)^m,\\
  Q_2(\alpha)&=\frac{1}{2}c_1(L_M\otimes F_\alpha)^{m-1}\cup c_1(TM),
\end{aligned}
\]
first for rational $\alpha$ then extending by continuity. We then have
\[ \dim H^0(L_M^k\otimes F_{k\alpha})=k^mQ_1(\alpha)+k^{m-1}Q_2(\alpha)
+ O(k^{m-2}),\]
for rational $\alpha\in\Delta$ and $k\gg 1$.
Since $H^0(X,L^k)=H^0(M,\pi_*(L^k))$ it follows using
Theorem~\ref{thm:trap} that
\[\begin{split} \dim H^0(L^k)=&k^{m+n}\int_\Delta Q_1\,d\mu+k^{m+n-1}\left(
\frac{1}{2}\int_{\partial\Delta}Q_1\,d\sigma+\int_\Delta Q_2\,d\mu
\right) \\ &+ O(k^{m+n-2}). \end{split}\]
We can now state the main result. 
\begin{thm} \label{thm:toricbtc}
  A rational piecewise-linear convex function $f$ on $\Delta$
  defines a test-configuration for $(X,L)$ with Futaki invariant
  \[ \frac{1}{2}\int_{\partial\Delta} fQ_1\,d\sigma +
  \int_\Delta fQ_2\,d\mu- \frac{a_1}{a_0}\int_\Delta fQ_1\,d\mu,\]
  where 
  \[ \begin{split} a_0 &= \int_{\Delta} Q_1\,d\mu,\\
    a_1 &= \frac{1}{2}\int_{\partial\Delta} Q_1\,d\mu+ \int_\Delta Q_2\,
    d\mu.
  \end{split}\]
  The norm of the test-configuration (ie. the norm of the induced
  $\mathbf{C}^*$-action on the central fibre) is given by 
  \[ \left(\int_\Delta (f-\overline{f})^2 Q_1\,d\mu\right)^\frac{1}{2},\]
  where $\overline{f}$ is the average of $f$ over $\Delta$ 
  with respect to the measure $Q_1\,d\mu$. 
\end{thm}
\begin{proof}
We define the test-configuration in the same way as was done
in~\cite{Don02} for toric varieties. Suppose $f<R$ for some integer $R$,
and let $\Delta^\prime$ be the polytope 
\[ \Delta^\prime = \{(x,t):x\in\Delta,\quad 0<t<R-f(x)\}\subset
\mathfrak{t}^*\times\mathbf{R}.\]
Let us assume that $\Delta^\prime$ is an integral polytope, otherwise we
could replace it by $k\Delta^\prime$ for an integer $k$. The polytope
$\Delta^\prime$ defines a polarised toric variety
$(W,\mathcal{O}_W(1))$. The face $\overline{\Delta^\prime}\cap
(\mathfrak{t}^*\times\{0\})$ is a copy of $\Delta$ so we get a natural
embedding $i:V\to W$ such that the restriction of $\mathcal{O}_W(1)$ to
$V$ is isomorphic to $\mathcal{O}_V(1)$. Write the $n+1$ torus action on
$W$ as $T\times\mathbf{C}^*$, where the $T$ action restricts in the
obvious way to $i(V)$. We now form the toric bundle
\[ Y = P\times_T W, \]
with the line bundle
\[ \mathcal{L} = \pi^*L_M\otimes\left(P\times_T\mathcal{O}_W(1)\right),\]
as we have done for $X$. Note that we have only twisted $W$ using the
first $n$ torus components (only $T$, not $T\times\mathbf{C}^*$). This
means that the corresponding $Q_1,Q_2:\Delta^\prime\to\mathbf{R}$ are
just the same as for $\Delta$, composed with the projection
$\Delta^\prime\to\Delta$. 

In~\cite{Don02} 
Donaldson showed that there is a $\mathbf{C}^*$-equivariant map 
\[ p: W\to \mathbf{P}^1,\]
with $p^{-1}(\infty)=i(V)$ such that the restriction of $p$ to
$W\setminus i(V)$ is a test configuration for $(V,\mathcal{O}_V(1))$.
Since the map $p$ is $T$-invariant, we can define
\[ \begin{split} 
     \tilde{p}&:Y\to \mathbf{P}^1\\
     (x,w)&\mapsto p(w).
   \end{split}
\]
and this defines a test configuration when restricted to
$Y\setminus \tilde{p}^{-1}(\infty)$. We can compute the Futaki invariant
of this test configuration in the same way as was done in~\cite{Don02}.
We have divisors $X_0=\tilde{p}^{-1}(0)$ and
$X\cong X_\infty=\tilde{p}^{-1}(\infty)$ defined by the vanishing of sections
$\sigma_0,\sigma_1$ of the line bundle $\tilde{p}^{-1}(\mathcal{O}(1))$
over $X$. When $k$ is large we therefore get the following exact
sequences:
\[\begin{split}
  & 0\to
  H^0(Y,\mathcal{L}^k(-1))\stackrel{\cdot\sigma_0}{\longrightarrow}
  H^0(Y,\mathcal{L}^k)\to H^0(X_0,
  \mathcal{L}^k|_{X_0})\to 0\\
  & 0\to
  H^0(Y,\mathcal{L}^k(-1))\stackrel{\cdot\sigma_1}{\longrightarrow}
  H^0(Y,\mathcal{L}^k)\to
  H^0(X_\infty,
  \mathcal{L}^k|_{X_\infty})\to 0
\end{split}
\]
The inclusion maps are multiplication by $\sigma_0$, $\sigma_1$. We
first see that the dimension $d_k$ of $H^0(X_0,\mathcal{L}^k)$ is the
same as that of 
$H^0(X,L^k)$. The $\mathbf{C}^*$-action acts with weight 0
on $\sigma_0$ and with weight 1 on $\sigma_1$, so the weight $w_k$ of
the action on $\bigwedge^{d_k} H^0(X_0,\mathcal{L}^k)$ is given by the
weight of the action on $\bigwedge^{d_k} H^0(X_\infty,\mathcal{L}^k)$ plus
the dimension of $H^0(Y,\mathcal{L}^k(-1))$. Since the action on
$H^0(X_\infty,\mathcal{L}^k)$ is trivial, we obtain
\[ w_k = \dim H^0(Y,\mathcal{L}^k(-1)) = \dim H^0(Y,\mathcal{L}^k) -
\dim H^0(X,L^k).\]
Using Theorem~\ref{thm:trap} we have
\[ \begin{split}
   w_k =& k^{m+n+1}\int_\Delta (R-f)Q_1\,d\mu  \\ &+k^{m+n}\left( \frac{1}{2}
\int_{\partial\Delta} (R-f)Q_1\, d\sigma+\int_\Delta (R-f)Q_2\, d\mu
\right) + O(k^{m+n-1}).
\end{split}\]
Recall that if
\[\begin{split} 
  d_k = a_0k^{m+n} + a_1 k^{m+n-1} + O(k^{m+n-2}), \\
  w_k = b_0k^{m+n+1}+b_1k^{m+n}+O(k^{m+n-1}), 
\end{split}
\]
then the Futaki invariant is $\frac{a_1}{a_0}b_0-b_1$. This gives the
required formula for the Futaki invariant. The formula for the norm of
the test configuration can be shown in the same way. 
\end{proof}

Let us assume now that $M$ is just a point, so that $(X,L)$ is a
polarised toric variety with corresponding polytope $\Delta$. In this
case $Q_1=1$ and $Q_2=0$, so we get back Donaldson's result
in~\cite{Don02}: A rational
piecewise linear convex function $f$ on $\Delta$ defines a
test-configuration for $X$ whose Futaki invariant is given by
$\frac{1}{2}\mathcal{L}(f)$, where
\[ \mathcal{L}(f)= \int_{\partial\Delta} f\,d\sigma - a\int_\Delta f\,
d\mu,\]
and $a=\frac{Vol(\partial\Delta,d\sigma)}{Vol(\Delta,\mu)}$.
Note that this is slightly different from the formula in~\cite{Don02}
because of our different convention for the definition of the Futaki
invariant. Let us rescale the Futaki invariant to be $\mathcal{L}(f)$ to
avoid unnecessary factors of two below. 

Let us see what the modified Futaki invariant is (see
Section~\ref{sec:Kstab}). The maximal torus of
automorphisms is the standard torus action on the toric variety. The
toric test-configurations we have defined are compatible with this torus
by definition. 
An affine linear function on $\Delta$ corresponds to a holomorphic vector
field on $X$ by the Hamiltonian construction. The inner product of two
vector fields is just the $L^2$ inner product of their Hamiltonians
normalised to have mean zero (see Section~\ref{sec:Kstab}). This means
that the normalised  Hamiltonian $B$ of the extremal vector field
satisfies
\[ \mathcal{L}(H) = -\int_\Delta BH\,d\mu \]
for all affine linear functions $H$ (the minus sign appears because the
test-configuration corresponding to $H$ by the above construction is
the product configuration corresponding to $-H$). 
There is a unique such $B$ and it
can be computed easily for specific toric varieties. 
Given a convex piecewise-linear function $f$ and a vector field
with Hamiltonian $H$, the inner product of the induced
$\mathbf{C}^*$-actions on the central fibre of the test-configuration
induced by $f$ is the $L^2$ product of $-f$ and $H$ normalised to have
zero mean. This is because the weights of the $\mathbf{C}^*$-action on
the central fibre
induced by the test-configuration are just the values of $-f$ plus some 
constant. This means that the modified Futaki invariant is
\[ \mathcal{L}(f)+\int_\Delta fB\,d\mu = \int_{\partial\Delta} f\,d\sigma
-\int_\Delta (a-B)f\, d\mu.\]
By the definition of $B$, this is zero for all affine linear $f$. If we
define $A=a-B$, then we see that the modified Futaki invariant is given
by 
\[ \mathcal{L}_A(f) = \int_{\partial\Delta} f\,d\sigma - \int_\Delta
Af\, d\mu. \]
Donaldson defines $\mathcal{L}_A$ with this formula 
for all bounded $A$ and conjectured
that if $\mathcal{L}_A(f)>0$ for all non-affine convex functions
$f$, then there is a K\"ahler metric on the toric variety with scalar
curvature given by the function $A$. There is a unique affine linear $A$
such that $\mathcal{L}_A(H)=0$ for all affine linear $H$, and our
discussion shows that for this $A$ the condition $\mathcal{L}_A(f)>0$
for all  non-affine convex functions $f$ means that the toric variety is
relatively K-polystable with respect to toric test-configurations. 

Let us now see what uniform K-polystability corresponds to.
Define the projection map $\pi:C(\Delta)\to C(\Delta)$ onto the
$L^2$-orthogonal 
complement of the space of affine linear functions. By definition 
$(X,L)$ is uniformly K-polystable with respect to toric
degenerations, if there exists a $\lambda>0$ for which
\begin{equation*}
  \mathcal{L}(f)\geq\lambda\Vert\pi(f)\Vert_{L^\frac{n}{n-1}},
\end{equation*}
for all convex $f$. The choice of the $L^{n/(n-1)}$-norm will become
clear at the end of Section~\ref{sec:toricmod}.
We can summarise all this as follows. 

\begin{prop} Let $\Delta$ be a polytope corresponding to the polarised
  toric variety $(X,L)$ of dimension $n$. 
  \begin{itemize}
    \item $(X,L)$ is \emph{K-semistable} for toric test-configurations
       if $\mathcal{L}(f)\geq0$ for all rational piecewise-linear
       convex functions $f$. 
     \item If in addition
  $\mathcal{L}(f)=0$ if and only if $f$ is affine linear, then $(X,L)$
  is \emph{K-polystable}.
\item $(X,L)$ is \emph{uniformly K-polystable} if there exists
  $\lambda>0$ such that for all convex functions $f$,
  \begin{equation}\label{eq:unifineq}
    \mathcal{L}(f)\geq\lambda\Vert\pi(f)\Vert_{L^\frac{n}{n-1}}. 
  \end{equation}
\item Let $A$ be the unique affine linear function such that
  $\mathcal{L}_A(H)=0$ for all affine linear $H$. Then $(X,L)$ is
  \emph{relatively K-polystable} if $\mathcal{L}_A(f)\geq0$ for all
  rational piecewise-linear
  convex functions $f$, with equality only if $f$ is affine linear. 
\end{itemize}
\end{prop}

\section{Toric surfaces}\label{sec:toricsurface}

We now restrict attention to toric surfaces. We first prove the
following
\begin{thm} \label{thm:unifKstab}
  A K-polystable toric surface is uniformly K-polystable. 
\end{thm}
Let the toric surface correspond to the polygon $P$ containing the
origin. Call a convex
function $f$ normalised if $f(0)=0$ and $f\geq 0$ on $P$. 
In~\cite{Don02} Donaldson proved that on a $K$-polystable toric surface
there exists $\lambda>0$, such that
\[ \mathcal{L}(f)\geq\lambda\int_{\partial P}f\,d\sigma \]
for all normalised convex functions $f$. To prove our result, it is
therefore enough to show the following, which we will prove in the next
subsection. 
\begin{prop}\label{prop:ineq1}
  There exists a constant $C$ such that for all non-negative
  continuous convex functions $f$ on $P$, 
  \[ \Vert f\Vert_{L^2(P)}\leq C\int_{\partial P} f\, d\sigma.\]
\end{prop}
Together with Donaldson's result, this shows that on a K-polystable
toric surface there exists $\mu>0$ such that
\[ \mathcal{L}(f)\geq\mu\ \Vert f\Vert_{L^2(P)}, \]
for all normalised convex functions $f$. 
This implies the inequality~\ref{eq:unifineq} for all convex functions
$f$ since $\Vert \pi(f)\Vert\leq\Vert f\Vert$, and both sides of the
inequality are invariant under adding affine linear functions to $f$. 

In
Section~\ref{sec:measuremaj} we will reprove part of Donaldon's
result, namely the fact that on a toric surface 
if $\mathcal{L}(f)\geq 0$ for all convex
$f$ and $\mathcal{L}(f)=0$ for some convex $f$ which is not affine linear,
then $\mathcal{L}(h)=0$ for a simple piecewise linear convex function
$h$.
This is a piecewise linear function with one ``crease''. Then in
Section~\ref{sec:semistab} we prove a decomposition theorem for
K-semistable polygons. 

\subsection{Uniform K-stability}\label{sec:toricmod}
The aim of this section is to prove Proposition~\ref{prop:ineq1}. 
Before giving the proof we need two lemmas.
Define a \emph{simple piecewise
linear function} on $\mathbf{R}^2$ to be a function of the form

\[ f(x) = \max(\lambda(x) + c,0),\]

\noindent where $\lambda :\mathbf{R}^2\to\mathbf{R}$ is a linear
function. We call the line $\lambda(x)=-c$ the \emph{crease} of $f$. 

\begin{lem}\label{lem:ineq1}
  Let $h$ be a simple piecewise linear convex function on the triangle 
  $\Delta=\{(x,y)\,|\,x,y\geq 0,\,x/a+y/b\leq 1\}$
    with $h(a,0)=h(0,b)=0$. Then there is a constant $C$ independent of
    $a$ and $b$ such that
    
    \[ \left(\int_\Delta h^2\, d\mu\right)^{1/2}\leq C\left(\int_0^a h(x,0)\,
    dx +  \int_0^b h(0,y),dy\right). 
    \]
\end{lem}

\begin{proof}
  Suppose the crease of $h$ is the segment $(c,0),(0,d)$. The inequality
  is invariant under multiplying $h$ by a constant so we can assume
  $h(0,0)=1$, so that on $\{(x,y)|x,y\geq 0,\, dx+cy\leq cd\}$ we have
  $h(x,y)=1-\frac{x}{c}-\frac{y}{d}$. 

  We have
  \[ \left(\int_\Delta h^2\, d\mu\right)^{1/2} =
  \sqrt{cd}\left(\int_0^1\int_0^{1-x}(1-x-y)^2\,dy\,dx\right)^{1/2} =
  C\sqrt{cd},
  \]
  for some constant $C$, and

  \[ \int_0^a h(x,0)\,dx +  \int_0^b h(0,y)\,dy = \frac{c+d}{2}. \]

  \noindent The result follows since $2\sqrt{cd}\leq c+d$.
\end{proof}

\begin{lem} \label{lem:ineq2}
  Let $f$ be a non-negative convex funcion on the triangle
  $\Delta = \{(x,y)\,|\, x,y\geq 0, x/a+y/b\leq 1\}$, such that $f(x,0),f(0,y)$ are
  non-increasing. There exists a constant $C$, independent of $a,b$ such that 
  \[ \begin{split} \left(\int_\Delta f^2\, d\mu\right)^{1/2}\leq
    C\Big(&\int_0^a
  f(x,0)\,dx + \int_0^b 
  f(0,y)\, dy \\ +&\sqrt{ab}\cdot\max\{f(a,0),f(0,b)\}\Big).\end{split}\]
\end{lem}
\begin{proof}
  It is enough to prove the assertion for piecewise linear convex
  functions. Let $f$ be a piecewise linear convex function on
  $\Delta$, and denote by $x_0=a> x_1>\ldots > x_k = 0$ and
  $y_0=b >y_1>\ldots>y_l=0$ the points where the restriction of 
  $f$ to the $x$ and $y$ axes is non-linear. Define the points
  $X_i=(x_i,0)$ and $Y_j=(0,y_j)$. We will define a new convex
  function $\tilde{f}$ which is equal to $f$ on the edges $X_0X_k$ and
  $Y_0Y_l$
  and dominates it in the interior of $\Delta$. We will then prove the inequality
  for $\tilde{f}$. 

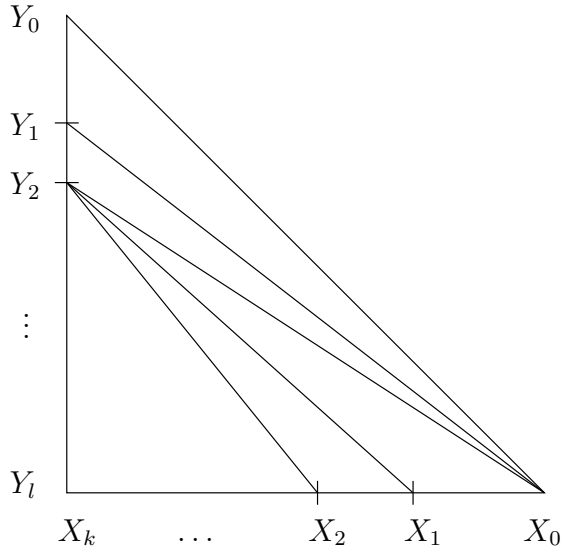
  \begin{figure}[htbp]
    \begin{center}
      \input{lemma424.pstex_t}
    \end{center}
  \end{figure}
  
  We define $\tilde{f}$ by induction on rectangles
  $\overline{X_iY_jY_0X_0}$. We start by defining $\tilde{f}$ on the
  ``rectangle'' $\overline{X_0Y_0Y_0X_0}$ to be the linear interpolation
  between the values of $f$ at $X_0$, and $Y_0$. 
  Supposing we have defined $\tilde{f}$ for
  some $i,j$, consider the triangles $\overline{X_iX_{i+1}Y_j}$ and
  $\overline{X_iY_{j+1}Y_j}$. Define linear 
  functions $u$ and $v$ on these
  triangles, which are equal to $f$ on $\partial\Delta$. Extend
  $\tilde{f}$ by the function which is smaller on the intersection of
  the two triangles (or either function if they are equal). Note that
  $\tilde{f}$ is the greatest convex function equal to $f$ on the
  two orthogonal edges of $\Delta$. 

  Define $g = \tilde{f}-g_0$, where $g_0$ is a linear function on
  $\Delta$ such that $g_0(a,0)=f(a,0)$, $g_0(0,b)=f(0,b)$ and
  $g_0(0,0)=\max\{f(a,0),f(0,b)\}$. 
  We can write $g$ as a sum of non-negative simple piecewise linear
  functions $g_1,\ldots,g_k$ with creases inside $\Delta$. 

  We have
  \[ \Vert \tilde{f} \Vert_{L^2(\Delta)}\leq\sum_{i=1}^k \Vert
  g_i\Vert_{L^2(\Delta)} + \Vert g_0\Vert_{L^2(\Delta)}.\]

  \noindent The sum is handled by Lemma~\ref{lem:ineq1} and the last
  term is bounded above by $\sqrt{ab}\cdot\max\{f(a,0),f(0,b)\}$ by definition.
\end{proof}

\begin{proof}[Proof of Proposition \ref{prop:ineq1}]
  Denote by $P_1,\ldots,P_k$ the vertices of the polytope (also let
  $P_{k+1}=P_1$) and by $e_i$
  the edge joining $P_i,P_{i+1}$ for $i=1,\ldots,k$.
  On each edge
  $e_i$ choose a point $Q_i$ where the restriction of $f$ to $e_i$ is
  minimal. Let us assume for simplicity that none of the $Q_i$ coincide
  with a vertex (we can achieve this by perturbing $f$ slightly).
  We first restrict attention to the triangles $Q_iP_{i+1}Q_{i+1}$. The
  property of these triangles that we need is that $f$ is non-decreasing
  on the edges $Q_iP_{i+1}$ and $Q_{i+1}P_{i+1}$. 

  
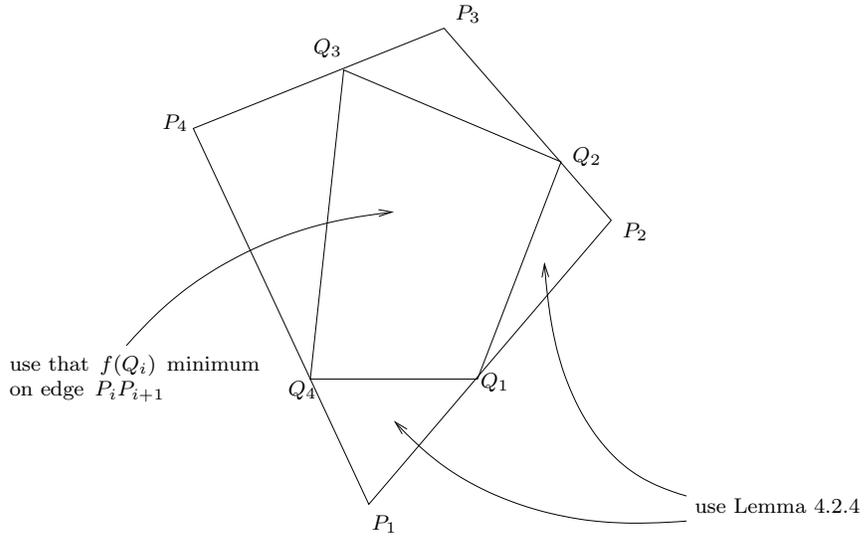
\begin{figure}[htbp]
    \begin{center}
      \input{prop422.pstex_t}
    \end{center}
  \end{figure}

  We apply Lemma~\ref{lem:ineq2} to the triangle $Q_iP_{i+1}Q_{i+1}$ to
  get that the $L^2$ norm of $f$ on this triangle is bounded above by
  $C$ times the integral on the segments $Q_iP_{i+1}$ and
  $P_{i+1}Q_{i+1}$ plus $\sqrt{Vol(P)}\max\{f(Q_i),f(Q_{i+1})\}$. 
  The $L^2$ norm of $f$ on the interior of the polygon $Q_1Q_2\ldots
  Q_k$ is bounded above by $\sqrt{Vol(P)}
  \max\{f(Q_1),\ldots,f(Q_k)\}$. In sum we obtain that for some constant
  $C_1$,
  \[ \Vert f\Vert_{L^2(P)}\leq C_1\left(\int_{\partial P}f\, d\sigma +
  \max\{f(Q_1),\ldots,f(Q_k)\}\right). \]

  \noindent Since $f(Q_i)$ is the minimum of $f$ on the edge $e_i$, we
  have 
  \[ f(Q_i)\leq \frac{1}{Vol(e_i)}\int_{P_i}^{P_{i+1}} f\,d\sigma \leq
  C_2\int_{\partial P} f\, d\sigma. \]
 
  \noindent for some $C_2$. This completes the proof of the result.
\end{proof}

This result does not hold in higher dimensions because of the way the
$L^p$ norms scale. For $a<1$, consider the
convex function
\[ h(x_1,\ldots,x_n) = \begin{cases}
  1-\frac{1}{a}(x_1+\ldots+ x_n),\quad &\text{if } x_1+\ldots+x_n < a,
  \\
  0, \quad \text{otherwise}.\end{cases}\]
on the set $\{(x_1,\ldots,x_n)\,|\,\sum x_i \leq 1\}$. The $L^p$-norm on
the interior of the set is $a^{n/p}C_1$ for some constant $C_1$ and the
$L^1$-norm on the boundary is $a^{n-1}C_2$ for another constant $C_2$.
Therefore the natural inequality to consider for $n>2$ is 
\[ \Vert f\Vert_{L^p(P)}\leq C\int_{\partial P} f\,d\sigma, \]
with $\displaystyle{p=\frac{n}{n-1}}$. This is the reason for our choice
of norm in the definition of uniform K-stability in
Section~\ref{sec:uniformKstab}. It is an interesting question to
see whether the inequality holds in this form for $n>2$.

\subsection{Measure majorisation} \label{sec:measuremaj}

 The aim of this section is to prove the following result. Recall that
 $P$ is a polygon corresponding to a toric surface. 
\begin{thm}  Suppose that 
\[ \mathcal{L}(f) = \int_{\partial P} f\, d\sigma - a\int_P f\, d\mu 
\geq 0\] 
for all continuous convex functions $f$ on $P$, but there
is a continuous convex function $u$ on $P$ which is not affine linear and
$\mathcal{L}(u)=0$. Then there is a simple piecewise linear function $f$
with crease passing through $P$, such that $\mathcal{L}(f)=0$. 
\end{thm}

This is proved in Donaldson~\cite{Don02}, but we give a slightly 
different proof
based on a result in~\cite{CFM}.
The fact that $\mathcal{L}(f)\geq 0$ for all continuous convex functions
means that the measure $d\sigma$ \emph{majorises} $a\,d\mu$. In this case
(see~\cite{CFM}) there exists a family $\{T_x\}_{x\in P}$ of probability
measures on $P$ such that the barycentre of $T_x$ is $x$, and 

\begin{equation}\label{eq:measuredecomp}
  \sigma = a\int_P T_x\, d\mu(x).
\end{equation}

\noindent Note that this implies that $T_x$ is supported on $\partial
P$. Let us denote the convex hull of its support by $l_x$ since we will
normally think of them as line segments. For $f$
convex, the Jensen inequality implies $T_x(f)\geq f(x)$ with equality if
and only if $f$ is linear when restricted to $l_x$. Hence 

\[ \int_{\partial P} f\, d\sigma = a\int_P T_x(f)\, d\mu(x) \geq a\int_P
f(x)\, d\mu(x) \] 

\noindent with equality if and only if $f$ is linear when
restricted to $l_x$ for $\mu$-almost every $x$. From this we can
immediately see a case when there is a simple piecewise linear function
giving equality. If there is a line $L$ through $P$ such that the
set of $x$ with $l_x$ intersecting $L$ transversely\footnote{By 
two convex sets intersecting transversally, we mean that
their interiors (in the case of a line segment the complement of its
endpoints) intersect.} has measure zero (with
respect to $\mu$), then any simple piecewise linear function with crease
$L$ will do. We wish to show that if there is no such $L$, then the only
convex functions giving equality are the linear ones. We need the
following

\begin{lem} If a convex function $f$ is linear when restricted to the
  convex sets $l_x$ and $l_y$ which intersect transversally, 
  then $f$ is linear when restricted to the convex hull of $l_x\cup l_y$.  
\end{lem}

\begin{proof} Suppose $l_x$ and $l_y$ are line segments.
  Let us denote the convex hull of $l_x\cup l_y$ by $S$.  By
  subtracting a linear function from $f$, we can assume that $f$
  restricted to $l_x\cup l_y$ is zero. Since $f$ is convex, it follows
  by definition, that $f$ is non-positive on $S$. Also, note that for
  any point $p$ in $S$, we can find a point $q$ in $l_x\cup l_y$ such
  that the segment $pq$ intersects $l_x\cup l_y$ in a point $r$ with $q$
  lying between $p$ and $r$ (see Figure~\ref{fig:pqr}). Then,
  since $f(q)=f(r)=0$ and $f(p)\leq 0$, we must have $f(p)=0$ since $f$
  is convex when restricted to $pr$. Thus, $f$ is identically zero on
  $S$.  

  \begin{figure}[htbp]
    \begin{center}
      \input{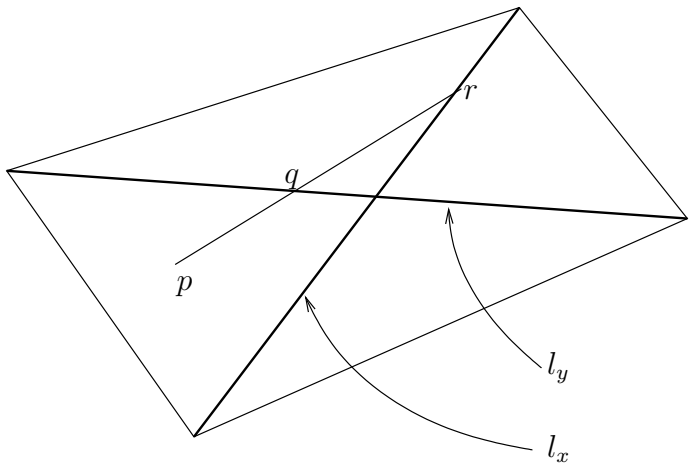}
      \caption{Since $f(q)=f(r)=0$, by convexity $f(p)=0$.\label{fig:pqr} }
    \end{center}
  \end{figure}

  In general if $l_x$ and $l_y$ are convex sets, then we can apply the
  previous argument to all pairs of line segments contained in $l_x$ and
  $l_y$ which intersect transversally.  
\end{proof}

Let $\mathcal{L}(f)=0$ and let $E\subset P^o$ (the interior of $P$)
be the set of $x$ such that $f$
is linear when restricted to $l_x$. The complement of $E$ in $P$ has
measure zero. For such an $l_x$ if there is another $l_y$ which
intersects it transversally, then $f$ is linear on the convex hull of
$l_x\cup l_y$ and thus linear on a neighbourhood of $x$. We obtain that
either there is a line $L$ as above, or $f$ is linear on $P$. The proof
is thus complete.

\subsection{Semistable surfaces} \label{sec:semistab}

Suppose we decompose $P$ into subpolygons $Q_i$. On each $Q_i$ we have
the Lebesgue measure $d\mu$ and also a measure $d\sigma_i$ on $\partial
Q_i$ which is the restriction of $d\sigma$ (ie. it is equal to $d\sigma$
on edges of $Q_i$ which are
subsets of edges of $P$ and is zero on edges of $Q_i$
which lie on the interior of $P$). For any bounded function $A$ on $Q_i$ 
we can define the functional 
\[ \mathcal{L}_A(f) = \int_{\partial Q_i} f\, d\sigma_i - \int_{Q_i}
Af\, d\mu.
\]
There is a unique affine linear function $A$ for which
$\mathcal{L}_A(H)=0$ for all affine linear $H$. Let us say that
$(Q_i,d\sigma_i)$ is relatively K-polystable 
if $\mathcal{L}_A(f)>0$ for all non-affine convex functions $f$. 
This is the same as relative K-polystability of the pair $(V_{Q_i},D_i)$
where $V_{Q_i}$ is the variety corresponding to $Q_i$ and $D_i$ is the
divisor corresponding to the edges of $Q_i$ where $d\sigma_i$ vanishes. 
In~\cite{Don02}, Donaldson conjectured 
that if $(Q_i,d\sigma_i)$ is
relatively K-polystable then there exists a complete extremal metric on
$V_{Q_i}\setminus D_i$. 
Donaldson then suggested that if the toric variety corresponding to a
polygon $P$ is not K-polystable,
then there should be a canonical decomposition of $P$ into subpolygons
$Q_i$, with measures $d\sigma_i$ induced by the measure on $\partial P$
as above and such that 
each $Q_i$ should either be relatively K-polystable, or be a
parallelogram in which two opposite edges lie on edges of $P$. We prove
this in the case of a K-semistable polygon. 
\begin{thm}\label{thm:semidecomp}
  A K-semistable polygon has a canonical decomposition into
  rational subpolygons, each of which is either K-polystable or is a
  parallelogram in which two opposite edges lie on edges of $P$. 
\end{thm}

Suppose that $P$ is K-semistable, so that $\mathcal{L}(f)\geq 0$ for all
convex functions on $P$. Recall that we have associated to each point
$x\in P$ a convex set $l_x$ containing $x$, 
which is the convex hull of a subset of
$\partial P$.  
Let $F$ be the set of line segments $l$ joining points on the boundary of
$P$ (and passing through the interior of $P$) such that
\[ \mu(\{x\in P,\, l\text{ intersects }l_x\text{ transversally}\})=0. \]
From the discussion in the previous subsection we see that $F$ is the
set of possible creases of a simple piecewise linear function $h$ such
that $\mathcal{L}(h)=0$. 
\begin{lem} \label{lem:semistab1}
  If $l_1,l_2$ are line segments in $F$ 
  then $l_1$ and $l_2$ cannot intersect.  
\end{lem}
\begin{proof}
  Let $C$ be the
  convex hull of the union $l_1\cup l_2$. 
  Suppose $l_1$ and $l_2$ intersect transversally in the interior. 
  For any $x$ in the interior of
  $C$ which doesn't lie on $l_1$ or $l_2$ we have that $l_x$ intersects
  $l_1$ or $l_2$ transversally. Therefore we cannot have both $l_1,l_2\in
  F$. 

  If $l_1$ and $l_2$ intersect on the boundary of $P$, in a point $y$,
  say, then for any $x$ in the interior of $C$ we have that $l_x$ either
  passes through $y$ or intersects $l_1$ or $l_2$ transversally. The set
  of $x$ with $l_x$ passing through $y$  but not intersecting $l_1$ or
  $l_2$ transversally must have measure zero,
  otherwise Equation~\ref{eq:measuredecomp} could not hold since the
  $\sigma$-measure of $y$ is zero. Therefore again, we
  cannot have both $l_1,l_2\in F$. 
\end{proof}

The line segments $l\in F$ are therefore a set of disjoint line segments in $P$.
Suppose now that $h$ is a simple
piecewise linear function with crease $l$, so that $\mathcal{L}(h)=0$.
In~\cite{Don02} (Section 6) Donaldson shows the
following:
\begin{enumerate}
  \item If one of the endpoints $x$ of $l$ is a vertex of $P$ then either
    the other endpoint is rational, or the other endpoint lies on an
    edge $J$ of $P$ such that every other line segment joining $x$ and
    $J$ is in $F$. 
  \item If $l$ joins two edges of $P$ which are not parallel, then
    $l$ is a rational line (its endpoints are rational). 
  \item If $l$ joins two parallel edges then either it is a rational
    line, or all other line segments parallel to $l$ joining the same
    two edges are in $F$. 
\end{enumerate}
The second possibility in case (1) is not possible because of
Lemma~\ref{lem:semistab1}. 
On the space of line segments joining two fixed edges of the polygon, the
functional $\mathcal{L}$ is a polynomial. Therefore it can only have
finitely many isolated zeroes. It follows that the set $F$ consists of
finitely many rational line segments, and a finite number of families
joining parallel edges of $P$ as in case (3). We therefore obtain the
decomposition we were after into rational subpolygons $Q_i$. Note that the
measure decomposition in Equation~\ref{eq:measuredecomp} can be
restricted to 
the $Q_i$ because for almost every point $x\in Q_i$ the line segment
$l_x$ lies inside $Q_i$ by our construction. This shows that each pair
$(Q_i,d\sigma_i)$ is 
K-polystable. The argument also shows that the decomposition is
canonical. 

\chapter{Ruled manifolds}\label{chap:ruled}

In this section we study extremal metrics on ruled manifolds. In
Section~\ref{sec:momentumsummary} we 
summarise the
momentum construction of Hwang-Singer~\cite{HS02} for writing down
circle invariant metrics on a ruled manifold starting with a function of
one variable (the momentum profile). Using this description in
Section~\ref{sec:metricdegen} we will
write down a sequence of degenerating 
metrics which models the deformation to the
normal cone of the zero section differential geometrically. To justify
this we compute the asymptotics of the Mabuchi functional along this
sequence of metrics and compare it to the Futaki invariant of the
deformation to the normal cone. 

We then concentrate on a ruled surface and construct explicit extremal
metrics on it in Section~\ref{sec:extremalruled}, 
including complete extremal metrics on the complement of
a divisor. The results fit in with the stability calculation in
Section~\ref{sec:ruledexample}. We then use these extremal metrics to
compute the infimum
of the Calabi functional for the unstable polarisations by writing down
a degenerate metric which achieves this infimum. The fact that it
is the infimum follows because it gives equality in Donaldson's
Theorem~\ref{thm:lowerCalabi}. 

\section{Summary of the momentum construction}
\label{sec:momentumsummary}
We briefly recall the momentum construction of circle invariant metrics
on line bundles. The reference for this section is
Hwang-Singer~\cite{HS02}. 
Let $(M,\omega_M)$ be a K\"ahler manifold of dimension $m$
and $(L,h)$ a Hermitian
holomorphic line bundle over $M$. Let $\gamma=-\sqrt{-1}\partial\bar{
\partial}\log h$ be the curvature form of $h$. Let $t=\log h$ be the
logarithm of the fibrewise norm function. We want to consider K\"ahler
metrics
on the total space of $L$ of the form
\begin{equation}\label{eq:calabiansatz}
  \omega = p^*\omega_M+2\sqrt{-1}\partial\bar{\partial} f(t),
\end{equation}
where $p:L\to M$ is the projection map and $f$ is a suitably convex
smooth function. Let us define $\tau=f^\prime(t)$. For each $f$, the
metric $\omega$ is invariant under the $S^1$ action rotating the fibres
of $L$, and $\tau$ is the moment map for this action. Let
$I\subset\mathbf{R}$ be the image of this moment map. Let $X$ be the
generator of the $S^1$ action normalised so that $\exp(2\pi X)=1$. The
function $\Vert X\Vert_\omega$ is constant on level sets of $\tau$ so we
can define a function $\phi: I\to [0,\infty)$ such that
$\phi(\tau)=\Vert X\Vert^2_\omega$. This function $\phi$ is called the
\emph{momentum profile} of the metric. We can reconstruct $f$ from
$\phi$; in fact $t$ and $\tau$ are related by the Legendre transform
with respect to $f$ and
the Legendre transform $F$ of $f$ satisfies $F^{\prime\prime}=1/\phi$.
This means that $t=F^\prime(\tau)$, and
\[ F(\tau) + f(t) = t\tau.\]
We can also express the metric on the fibres using $\phi$, namely it is
\[ \phi(\tau)\frac{\vert dz\vert^2}{\vert z\vert^2}, \]
where $z$ is a coordinate on the fibre. In other words, $\phi$ gives the
conformal factor relating the restriction of $\omega$ 
to the fibres, to the cylindrical metric. 
The advantage of this transformation is that in 
terms of $\phi$ the scalar curvature of $\omega$ is a second order
linear 
differential expression, so we can compute with it conveniently. 

We would now like to start with a momentum profile on an interval, and
define a metric. For this we first need the following data.
\begin{defn} \emph{Horizontal data} $(p:(L,h)\to (M,\omega_M), I)$
  consists of a Hermitian holomorphic line bundle over a K\"ahler
  manifold as above, together with a \emph{compatible} momentum interval
  $I\subset\mathbf{R}$.
  The interval $I$ is compatible if for all $\tau\in I$ the form
  $\omega_M(\tau)=\omega_M-\tau\gamma$ is positive.
\end{defn}

Given horizontal data, define a \emph{momentum profile} to be a smooth
function $\phi: I\to [0,\infty)$, which is positive on the interior of
$I$. Define the ruled manifold $X=\mathbf{P}(L\oplus\mathcal{O})$, and
let $S_0$, $S_\infty$ be the zero and infinity sections. By the Legendre
transform as above we obtain a function $f$ from $\phi$, and using
\ref{eq:calabiansatz} we define a
metric $\omega_\phi$ on a subset of $X$ with properties as follows. We
assume for simplicity that the momentum profile is a rational function
since that is all that we need for applications.  

\begin{thm}[see~\cite{HS02}]
  Let $I=[a,b]$, and suppose $\phi(a),\phi(b)=0$ and $\phi$ is
  a rational function positive on $(a,b)$. 
  The metric corresponding to the momentum
  profile $\phi$ has the following properties depending on the boundary
  conditions of $\phi$:
  \begin{eqnarray*}
    \phi^\prime(a)=2,\,\phi^\prime(b)=-2 & & \mbox{smooth metric on }
    X,\\
    \phi^\prime(a)=0,\,\phi^\prime(b)=-2 & & \mbox{complete metric on }
    X\setminus S_0,\\
    \phi^\prime(a)=2,\,\phi^\prime(b)=0 & & \mbox{complete metric on }
    X\setminus S_\infty,\\
    \phi^\prime(a)=0,\,\phi^\prime(b)=0 & & \mbox{complete metric on }
    X\setminus\{S_0\cup S_\infty\}.
  \end{eqnarray*}
\end{thm}

We now show that in the case of a complete metric if the order of vanishing
of the momentum profile is precisely 2, then the metric is
asymptotically hyperbolic as in Definition~\ref{defn:cuspsing}.

\begin{thm} Suppose $I=[0,1]$ and we have a momentum profile $\phi$ such
  that $\phi^\prime(0)=0$ and $\phi^{\prime\prime}(0)>0$. Then near the
  divisor $\tau^{-1}(0)$ the metric is asymptotically hyperbolic.
\end{thm}
\begin{proof}
  Let us assume for simplicity that $\phi^{\prime\prime}(0)=2$. We can
  then write
  \[ \phi(\tau) = \tau^2+ a\tau^3 + O(\tau^4) \]
  for small $\tau$ with some constant $a$. 
  In this proof by $O(\tau^k)$ we do not just mean
  bounded by $C\tau^k$ for some constant $C$, but that the function has
  a convergent Taylor expansion with terms of order at least $k$ for
  small $\tau$. We
  have
  \[ \frac{1}{\phi(\tau)}=\frac{1}{\tau^2}(1-a\tau)+O(1).\]
  Integrating this, by the definition of the Legendre transform we get 
  \begin{equation}\label{eq:def_t} 
    t=-\frac{1}{\tau}-a\log\tau + O(\tau).
  \end{equation}
  Integrating once more gives 
  \[ F(\tau) = -\log\tau-a\tau\log\tau+a\tau + O(\tau^2).\]
  The K\"ahler potential of the fibre metric
  is the Legendre transform of $F$ so we obtain
  \[ f(t) = t\tau - F(\tau) = \log\tau -1 + O(\tau).\]
  From Equation~\ref{eq:def_t} we obtain
  \[ \log t = -\log\tau + \log(1+a\tau\log\tau+O(\tau^2)).\]
  With the previous formula for $f(t)$ this implies
  \begin{equation}\label{eq:ftlogt}
    f(t)+\log t=-1+\log(1+a\tau\log\tau+O(\tau^2))+O(\tau)
  \end{equation}
  with slight abuse of notation. Since $-\log t$ is the K\"ahler
  potential of the hyperbolic cusp, to prove that our metric is
  asymptotically hyperbolic we need to show that the terms on the right
  hand side of Equation~\ref{eq:ftlogt} have covariant derivatives whose
  norms
  tend to zero as $z\to0$ (recall that $t=\log|z|$). The norm squared
  $\Vert dt\Vert^2$ with respect to the hyperbolic metric is $c t^2$ for
  some constant $c$. The worst term on the right of
  Equation~\ref{eq:ftlogt} is $\tau\log\tau$, since the other terms are
  products of powers of $\tau\log\tau$ with powers of $\tau$. To 
  prove the result it is therefore enough to show that for all $k>0$
  \[ t^k\frac{\partial^k}{\partial t^k}(\tau\log\tau)\to 0,\quad 
  \text{ as } \tau\to0.\]
  By the chain rule $\frac{\partial
  t}{\partial\tau}\frac{\partial}{\partial t}
  =\frac{\partial}{\partial\tau}$ and since $\frac{\partial
  t}{\partial\tau} = \frac{1}{\phi(\tau)}$, we get
  $\frac{\partial}{\partial t} = \phi(\tau)\frac{\partial
  }{\partial\tau}$. 
  By induction one can show that the term with smallest order of
  vanishing in $\frac{\partial^k}{\partial t^k}(\tau\log\tau)$ is
  $\tau^{k+1}\log\tau$. From Equation~\ref{eq:def_t} we see that $t\tau$
  is bounded as $\tau\to0$ and also $\tau\log\tau$ tends to zero as
  $\tau\to0$, hence
  \[ t^k\tau^{k+1}\log\tau = (t\tau)^k\cdot \tau\log\tau\to
  0,\quad\text{ as } t\to0.\]
  This completes the proof.
\end{proof}

Note that if the order of vanishing of $\phi$ is greater than two, then
the resulting metric is no longer asymptotically hyperbolic, but instead
the fibre metrics are asymptotic to
\[ \frac{ \vert dz\vert^2}{\vert z\vert^2 (\log\vert
z\vert)^{k/(k-1)}}\]
for some $k>2$. 

Define $Q:I\times M\to\mathbf{R}$ by
$Q(\tau)=\omega_M(\tau)^m/\omega_M^m$. The area of the fibres of $X$ is
$2\pi(b-a)$ and the volume of the zero section is $Q(a)Vol(M,\omega_M)$. 
As for the scalar curvature of this metric $\omega_\phi$, we have

\begin{thm}\label{thm:momentscal}
  Let $S_M(\tau)$ denote the scalar curvature of the
  metric $\omega_M(\tau)$. Then the scalar curvature of
  $\omega_\phi$ is given by 
  \[ S(\omega_\phi) = S_M(\tau) - \frac{1}{2Q}\frac{\partial^2}{
  \partial\tau^2}(Q\phi)(\tau).\]
\end{thm}

Since $\tau: X\to I$ is a moment map, we can consider the symplectic
reductions $M_c = \tau^{-1}(c)/S^1$ for $c$ in the interior of $I$. From
this point of view $\omega_M(c)$ and $S_M(c)$ give the induced K\"ahler form
and its scalar curvature on $M_c$. 

\section{A metric degeneration}\label{sec:metricdegen}
We will first construct a family of metrics on $\mathbf{P}^1$ using
K\"ahler potentials, and then
use their momentum profiles
to define a sequence of metrics on a ruled manifold and we compute the
asymptotic rate of change of the Mabuchi functional. 
We define the sequence using K\"ahler
potentials on $\mathbf{C}\setminus\{0\}$. Let $t=\log{|z|}$, and let
$g(t)$ be the K\"ahler potential of a cusp metric on $\mathbf{P}^1$
minus a point. 
We can take $g(t)$ to be a strictly convex smooth real valued function
such that $g(t) = -\log t$ for $t\gg 1$ and $g(t)$ is asymptotically
$-ct$ as $t\to -\infty$ (meaning that $g(t)+ct$ converges to a constant as 
$t\to -\infty$), 
where $2\pi c$ is the area. 

Let
\begin{equation*}
  g_s(t) = \begin{cases} g(t+s)-g(s+\frac{1}{s}),\quad &t<0 \\
		0, &t>1/s \end{cases}
\end{equation*}
and let $g_s$ be smooth and strictly convex on $(0,1/s)$. 
Define $h(t)$ to be smooth, $h(t)=0$ for $t<0$, and strictly convex for
$t>0$, such that $h(t)$ is asymptotically $dt$ for $t\gg 1$. Let
\begin{equation*} f_s(t) = h(t) + g_s(t) + ct. 
\end{equation*}
This is a strictly convex function on $\mathbf{R}$ and defines a metric
of area $2\pi (c+d)$ on $\mathbf{P}^1$. 
As $s\to\infty$, the potential $f_s$
approaches the potential of a cusp metric of area $2\pi c$ on the
interval $(-\infty,0)$, so the $\mathbf{P}^1$ breaks up into two pieces. 

Let us see what the corresponding momentum profiles look like. By
definition
$\tau_s=\frac{df_s}{dt}$, which changes with $s$, but we will normally
drop the subscript $s$. The momentum profile $\phi$ is defined by
$\phi_s(\tau)=1/(F_s)^{\prime\prime}$, 
where $F_s$ is the Legendre transform of $f_s$.
For each $s$, the momentum profile is a non-negative smooth function on
the interval $[0,c+d]$, positive on the interior and $\phi_s^\prime(0)=2,
\phi_s^\prime(c+d)=-2$. As $s\to\infty$, we have $\phi_s(c), \phi_s^\prime
(c)\to 0$.

We can use these momentum profiles to define metrics on our ruled
manifold $X$. We assume for simplicity that $c+d$ is
small enough so that the interval $[0,c+d]$
is a compatible momentum interval. 
\begin{rem}\label{rem:sesh} In
general we expect that we could allow $c+d$ to be as 
large as the Seshadri constant of
the zero section by possibly changing our choice
of $\omega_M$ and $h$, but we do not wish to discuss this in this thesis.
\end{rem}
Let $\omega_s$ be the metric corresponding to
$\phi_s$, and $S(\omega_s)$ the scalar curvature given by
Theorem~\ref{thm:momentscal} as
\begin{equation}\label{eq:scal}
  S(\omega_s) = S_M(\tau) - \frac{1}{2Q}\frac{\partial^2}{
\partial\tau^2}(Q\phi_s)(\tau).
\end{equation}
We observe that the scalar curvature is uniformly bounded as $t$ varies
since the $\phi_s$ are bounded in $C^2$. 

We would now like to compute the rate of change of the Mabuchi
functional as $s\to\infty$. 
The change in the Mabuchi functional is defined by
\[ \frac{d}{ds}\mathcal{M}(\omega_s) = -\int_X
\frac{df_s}{ds}(S(\omega_s) - \hat{S})
\frac{\omega_s^n}{n!}.\]

\noindent where $\hat{S}$ is the average scalar curvature of $X$, and
$n$ is the dimension of $X$, ie. $n=m+1$.  

We have 
\begin{equation*}
  \frac{df_s(t)}{ds} = \begin{cases}
    g^\prime(t+s)-\left(1-s^{-2}\right)g^\prime(s+s^{-1}), \quad & t <
    0, \\ 
	    0, & t > 1/s,
	  \end{cases}
\end{equation*}
and also
\begin{equation*}
  \tau_s = \frac{df_s(t)}{dt} = \begin{cases}
    g^\prime(t+s)+c,\quad & t<0, \\
    h^\prime(t)+c, & t>1/s.
  \end{cases}
\end{equation*}
We see that as $s\to\infty$, the limit of the integrand is (writing
$\tau$ for $\tau_\infty$)
\begin{align*}
  (\tau-c)(S(\omega_\infty)-\hat{S}), &\quad \tau < c \\
  0, &\quad \tau > c.
\end{align*}
Here $S(\omega_\infty)$ is defined by the formula (\ref{eq:scal}) for
$\phi_\infty=\lim_{t\to\infty}\phi_t$, considered as functions on
$[0,c+d]$ (note that $\tau_t$ changes with $t$). Although $\omega_\infty$
is a singular metric, $S(\omega_\infty)$ is a continuous function on
$X$ and the volume form is given by 
\[ \frac{\omega_\infty^n}{n!} = Q(\tau)d\tau\wedge d\theta
\wedge\frac{\omega_M^m}{m!},\]
which is bounded ($d\theta$ is the angular measure on the fibres). 
Since the convergence of the integrands is uniform (it is important here
that the scalar curvature remains uniformly bounded as the metric
degenerates), we can simply integrate the limit. We therefore find that
\begin{equation}\label{eq:limitmabuchi}
  \lim_{s\to\infty} \frac{d}{ds}\mathcal{M}(\omega_s) = \int_X
  (c-\tau)^+
  (S(\omega_\infty)-\hat{S})\frac{\omega_\infty^n}{n!},
\end{equation}
where $(c-\tau)^+=\max\{c-\tau,0\}$. 
We would like to show that up to a scalar multiple
this is the Futaki invariant of an algebraic
test-configuration. 

\subsection*{Deformation to the normal cone}
Let $(X,\mathcal{L})$ be a polarised variety (the line bundle here is
curly $\mathcal{L}$ to differentiate it from the line bundle $L$ over
$M$ in the previous section). Deformation to the normal cone of a
subscheme of $X$ was studied by Ross and Thomas~\cite{RT06}. Let us
recall the construction. Suppose $Z\subset X$ is a subscheme, and define
a test configuration $\mathscr{X}\to\mathbf{C}$ obtained by blowing up
$X\times\mathbf{C}$ along $Z\times\{0\}$. Denote the exceptional divisor
by $P$. The $\mathbf{C}^*$-action on $\mathscr{X}$ is induced by the
product action on $X\times\mathbf{C}$ acting trivially on $X$ and by
multiplication on $\mathbf{C}$. The central fibre of this test
configuration can be written as 
$\hat{X}\cup_E P$, where $\hat{X}$ is the blowup of $X$
along $Z$ with exceptional divisor $E$. If both $X$ and $Z$ are smooth,
then $P=\mathbf{P}(\nu\oplus\mathbf{C})$, the projective completion of
the normal bundle $\nu$ of $Z$ in $X$, and $E=\mathbf{P}(\nu)$. 

There is a choice of line bundles on $\mathscr{X}$. Let $\pi$ denote the
composition
\[ \pi : \mathscr{X}\to X\times\mathbf{C}\to X.\]
For a positive rational number $c$ let $\mathscr{L}_c$ be the
$\mathbf{Q}$-line bundle $\pi^*(\mathcal{L})-cP$. The restriction of this to the
general fibre of the test configuration is $\mathcal{L}$ and it is ample for
sufficiently small $c$. In fact it is ample for $c <\epsilon(Z)$, where
$\epsilon(Z)$ is the Seshadri constant of $Z$ (see~\cite{RT06}). 

In~\cite{RT06} the Futaki invariant of this test-configuration is
computed. Before stating the result we need some more definitions. Let
$\mathscr{I}_Z$ be the ideal sheaf of $Z$, and for a fixed
$x\in\mathbf{Q}_{>0}$ define $\alpha_i(x)$ by
\[ \chi(\mathcal{L}^k\otimes\mathscr{I}^{xk}_Z/\mathscr{I}^{xk+1}_Z) =
\alpha_1(x)k^{n-1}+\alpha_2(x)k^{n-2}+O(k^{n-3}),\quad k\gg 0,
xk\in\mathbf{N}.
\]
Define the \emph{slope} of $X$ by 
\[ \mu(X) = -\frac{nK_X\cdot \mathcal{L}^{n-1}}{2\mathcal{L}^n}.\]
The Futaki invariant is then
\[ F(\mathscr{X})= 
\int_0^c (c-x)\alpha_2(x)\,dx+\frac{c}{2}\alpha_1(0)- \left( \int_0^c
(c-x)\alpha_1(x)\,dx\right)\mu(X).\]
When $Z$ is a divisor we can use the Riemann-Roch formula to compute
\[ \alpha_1(x) = \frac{Z.(\mathcal{L}-xZ)^{n-1}}{(n-1)!},\quad \alpha_2(x) =
-\frac{Z.(K_X+Z).(\mathcal{L}-xZ)^{n-2}}{2(n-2)!}.\]

Recall now the ruled manifold $X=\mathbf{P}(L\oplus\mathcal{O})$ we
defined before, where $p:L\to M$ is a line bundle over a K\"ahler manifold
$(M,\omega_M)$. Let the $\mathbf{Q}$-polarisation $\mathcal{L}$ 
over $X$ be given by
a metric with momentum interval $[0,c+d]$ for rational $c,d>0$.   
Consider the deformation to the normal cone
$\mathscr{X}$ of the zero
section $S_0$ with parameter $c$. 
We now show that its Futaki invariant is up to a positive
multiple the asymptotic rate of change of the Mabuchi functional
in Equation~\ref{eq:limitmabuchi}.

\begin{thm}
  Using the notation from Equation~\ref{eq:limitmabuchi}, we have
  \[ \lim_{s\to\infty} \frac{d}{ds}\mathcal{M}(\omega_s) = 
  \int_X (c-\tau)^+ (S(\omega_\infty)-\hat{S})\frac{\omega_\infty^n}{n!} =
  2(2\pi)^n F(\mathscr{X}).\]
\end{thm}
\begin{proof}
  Let us write
  \begin{gather*} 
    A = \int_X (c-\tau)^+\frac{\omega_\infty^n}{n!},\quad B = \int_X
    (c-\tau)^+S_M(\tau)\frac{\omega_\infty^n}{n!},\\
    C = \int_X (c-\tau)^+\frac{1}{Q}\frac{\partial^2}{\partial \tau^2}
    (Q\phi_\infty)(\tau)\frac{\omega_\infty^n}{n!},
  \end{gather*} 
  so we need to show $F(\mathscr{X})=B-C/2-\hat{S}A$.
  
  First of all the average scalar curvature $\hat{S}$ is $2\mu(X)$. 
  Let us compute A. 
  The volume form is $\frac{\omega_M(\tau)^{n-1}}{(n-1)!}\wedge
  d\tau\wedge d\theta$, where 
  $d\theta$ is the angular measure on the fibres. The integrand is
  constant on the $S^1$ fibres and also on the level sets of $\tau$, so
  \[ A = 2\pi\int_0^c (c-\tau)\mathrm{Vol}(M,\omega_M(\tau))\,d\tau.\]
  The volume $\mathrm{Vol}(M,\omega_M(\tau))$ is the volume of the zero
  section with respect to the metric $\omega_M-\tau\gamma$, ie.
  $(2\pi)^{n-1}S_0.(\mathcal{L}-\tau S_0)^{n-1}/(n-1)!$, which is just
  $(2\pi)^{n-1}\alpha_1(\tau)$, since 
  $Z=S_0$. We therefore have
  \begin{equation}\label{eq:A} 
    A = (2\pi)^n\int_0^c (c-x)\alpha_1(x)\,dx.
  \end{equation}

  Now let us move on to B. We have 
  \[ B = 2\pi
  \int_0^c(c-\tau)\int_M S_M(\tau)\frac{\omega_M(\tau)^{n-1}}{(n-1)!}\,
  d\tau.\]
  The integral over $M$ is the total scalar curvature of $M$ with the
  polarisation $\mathcal{L}-\tau S_0$, which using the adjunction formula is 
  \[ \begin{split}
  -(2\pi)^{n-1}\frac{K_{S_0}.(\mathcal{L}-\tau S_0)^{n-2}}{(n-2)!}&=
  -(2\pi)^{n-1}\frac{(K_X+S_0).S_0.(\mathcal{L}-\tau S_0)^{n-2}}{
  (n-2)!}\\ &=2(2\pi)^{n-1}\alpha_2(\tau),
  \end{split}\]
  so that we have
  \begin{equation}\label{eq:B}
    B = 2(2\pi)^n\int_0^c (c-x)\alpha_2(x)\,dx.
  \end{equation}

  Finally let us compute C.
  \[ C = 2\pi\int_M\int_0^c (c-\tau)\frac{\partial^2}{\partial\tau^2}
  (Q\phi_\infty) (\tau)
 \, d\tau \frac{\omega_M^{n-1}}{(n-1)!},\]
 using that $\omega_M(\tau)^{n-1}=Q(\tau)\omega_M^{n-1}$. We can
 integrate the inner integral by parts remembering that
 $\phi_\infty(0)=\phi_\infty(c)=0, \phi_\infty^\prime(0)=2$ and
 $\phi_\infty^\prime(c)=0$. We get 
 \begin{equation}\label{eq:C}
   C = -4\pi c\int_M Q(0)\frac{\omega_M^{n-1}}{(n-1)!} = -2(2\pi)^n c\alpha_1(0).
 \end{equation}
 Putting together equations \ref{eq:A}, \ref{eq:B} and \ref{eq:C} 
 we obtain the required result.
\end{proof}

This result shows that for a ruled manifold if the deformation to the
normal cone of the zero section destabilises for $c$ sufficiently small
(cf. Remark~\ref{rem:sesh}), then the Mabuchi
functional is not bounded from below. In particular the manifold cannot
admit a cscK metric, although we knew this from the result in
Section~\ref{sec:lwrcalabi} already. With this approach however we have
a direct relationship between a metric degeneration and the asymptotics
of the Mabuchi functional, and a corresponding algebro-geometric
test-configuration and its Futaki invariant. As we mentioned before,
understanding this relationship in general is an important problem.
It should be possible to extend the calculation here to deformation to
the normal cone of a smooth divisor in a general manifold, transferring
the metrics we have constructed here to a suitable tubular neighbourhood
of the
divisor. At the time of writing this thesis I have not yet worked out
how to do this. 

We now perform a similar calculation but with the more general
test-con\-fi\-gu\-ra\-tions for toric bundles we constructed in
Section~\ref{sec:kstabtoric}. For this, note that $X$ is a toric bundle
with base $M$ and fibre $(\mathbf{P}^1,\mathcal{O}(l))$ with moment
``polytope'' $[0,l]$. Let the principal
$\mathbf{C}^*$-bundle $P$ on $M$ be the complement of the zero section
in $L^{-1}$ so that the polarisation we defined in
Section~\ref{sec:kstabtoric}
coincides with the one obtained from the momentum construction for the
same interval. The ample line bundle $L_M$ over $M$ is a holomorphic line
bundle with first Chern class $[\frac{1}{2\pi}\omega_M]$. 
The functions $Q_1,Q_2:[0,l]\to\mathbf{R}$ are given by 
\[ \begin{split}
  Q_1(\tau) &= \frac{1}{(2\pi)^m m!}\int_M (\omega_M-\tau\gamma)^m =
  (2\pi)^{-m}\int_M
      Q(\tau)\frac{\omega_M^m}{m!},\\
      Q_2(\tau) &= \frac{1}{2(2\pi)^{m-1}(m-1)!}
      \int_M (\omega_M-\tau\gamma)^{m-1}\wedge
      c_1(M)\\ &
      =\frac{1}{2(2\pi)^m}\int_M S_M(\tau)Q(\tau) \frac{\omega_M^m}{m!}.
    \end{split}\]
Now according to Theorem~\ref{thm:toricbtc} any rational
piecewise linear convex function $h$ on $[0,l]$ 
defines a test-configuration $\mathscr{X}$ for $X$, with
Futaki invariant equal to 
\[ F(\mathscr{X}) = \frac{h(0)Q_1(0)+h(l)Q_1(l)}{2}+\int_0^l
h(\tau)Q_2(\tau)\, d\tau -
\frac{a_1}{a_0}\int_0^l h(\tau)Q_1(\tau)\, d\tau,\]
where $\hat{S}=2\frac{a_1}{a_0}$ is the average scalar curvature of $X$.
Let us define $F(h)$ with the same formula for any piecewise
\emph{smooth} function $h$. 

\begin{thm}\label{thm:ruledfutaki}
  Let $h$ be any piecewise smooth convex function on
  $[0,l]$ which is smooth on the intervals $[l_i,l_{i+1}]$ for some
  $0=l_0<l_1<\ldots<l_N=l$. Let $\phi\in C^2([0,l])$ be non-negative,
  satisfying 
  \[ \begin{split}
        \phi(0)=\phi(l_i)=\phi(l)=0\, \text{for all }i,\\
	\phi^\prime(0)=2, \phi^\prime(l)=-2.
      \end{split} \]
      Suppose in addition that $h$ is linear on any interval
      $[l_i,l_{i+1}]$ on which $\phi$ does not vanish identically.
  We then
  have
  \begin{equation}\label{eq:singfutaki} 
    F(h) = \frac{1}{2(2\pi)^n}\int_X
    h(\tau)(S(\omega_\phi)-\hat{S})\frac{\omega_\phi^n}{n!},
  \end{equation}
  where $\omega_\phi$ is the singular metric corresponding to the
  ``momentum profile'' $\phi$. 
\end{thm}
Note that $h$ does not define a test-configuration since it is not
piecewise linear, and also $\phi$ does not define a metric since it
vanishes on a subset of $(0,l)$. On the other hand $h$ can be uniformly
approximated with piecewise linear functions, and $\phi$ can be
approximated in $C^2$ with momentum profiles which do define metrics.
Equation~\ref{eq:singfutaki} will then hold in the limit for such
approximating sequences of test-configurations and metrics. This is the
setting in which the
result will be used in Section~\ref{sec:extremallwr}. 
\begin{proof}
  As in the previous proof, we write
  \[ \begin{split}
    A = \int_X h(\tau)\frac{\omega_\infty^n}{n!},\quad B=\int_X
    h(\tau)S_M(\tau)\frac{\omega_\infty^n}{n!}, \\
    C = \int_X h(\tau)\frac{1}{Q}\frac{\partial^2}{\partial \tau^2}
    (Q\phi_\infty)\tau\frac{\omega_\infty^n}{n!},
  \end{split}
  \]
  and we need to show $F(\mathscr{X}) = B-C/2-\hat{S}A$. 
  In the same way as above, we get
  \[ A = (2\pi)^n\int_0^l h(\tau)Q_1(\tau)\,d\tau,\]
  and also 
  \[ B = 2(2\pi)^n\int_0^l h(\tau)Q_2(\tau)\,d\tau,\]
  since $Q_2(\tau)$ is the total scalar curvature of $M$ with polarisation
  $c_1(L_M)-\tau c_1(L)$. 

  Also as before,
  \[ C = 2\pi\int_M\int_0^l h(\tau)\frac{\partial^2}{\partial\tau^2}
  (Q\phi_\infty)(\tau)\, d\tau\frac{\omega^{n-1}}{(n-1)!}.\]
  Integrating by parts, using the assumption that $h$ is linear on
  $[l_i,l_{i+1}]$ if $\phi$ does not vanish there, we get that
  \[ \int_{l_i}^{l_{i+1}} h(\tau)\frac{\partial^2}{\partial\tau^2}
  (Q\phi_\infty)\, d\tau = \begin{cases} -2Q_1(0)h(0)\,&\text{ if }i=0,\\
    0\,&\text{ if }0<i<N-1,\\
    -2Q_1(l)h(l)\,&\text{ if } i=N-1.
  \end{cases}\]
  Summing up, we obtain
  \[ C = -2(2\pi)^n(Q_1(0)h(0)+Q_1(l)h(l)).\]
  Putting everything together, we obtain the result.
\end{proof}

\section{Extremal metrics on ruled surfaces}\label{sec:extremalruled}
We now specialise to a ruled surface. The base manifold $M$ is now a
genus $2$ curve equipped with a metric $\omega_M$ of constant scalar 
curvature and area $2\pi$, 
and $L$ is a degree $-1$ line bundle (this gives the same variety as a
degree $1$ line bundle that we used 
in Section~\ref{sec:ruledexample} but enables us to
use momentum profiles on $[0,m]$ instead of $[-m,0]$). The variety $X$ is
the ruled surface $\mathbf{P}(L\oplus\mathcal{O})$. We pick a Hermitian
metric on $L$ with curvature form $i\omega_M$. 
Using the description of circle-invariant metrics on ruled manifolds in
Section~\ref{sec:momentumsummary}
we now construct extremal metrics on $X$. This was done by
T\o{}nessen-Friedman in~\cite{TF97}, but we also construct complete
metrics on $X\setminus{S_0}$ and $X\setminus{S_\infty}$.  We choose the
polarisation $\mathcal{L}=C+mS_\infty$ where $C$ is a fibre as in
Section~\ref{sec:ruledexample}. This is equivalent to working on the
momentum interval $[0,m]$. The function $Q$ needed in the formula for
the scalar curvature (Theorem~\ref{thm:momentscal}) is given by
$Q(\tau)=1+\tau$, so the expression for the scalar curvature is 
\[ S(\omega_\phi) = \frac{1}{2(1+\tau)}(-4-[(1+\tau)\phi]^{
\prime\prime}),\]
where $\omega_\phi$ is the metric corresponding to a momentum profile
$\phi:[0,m]\to\mathbf{R}$. 

In order to find extremal metrics, we therefore need to find momentum
profiles $\phi:[0,m]\to\mathbf{R}$ satisfying various boundary
conditions, and solving the ODE 
\[ \frac{1}{2(1+\tau)}(-4-[(1+\tau)\phi]^{\prime\prime
})=A\tau+B\]
for some constants $A,B$ since the gradient of a function $h(\tau)$ is
holomorphic if and only if $h$ is linear. More explicitly,
\begin{eqnarray*}
  \left[(1+\tau)\phi\right]^{\prime\prime}&=&-2A\tau^2-2(A+B)\tau-2B-4,\\
  \left[(1+\tau)\phi\right]^\prime&=&-\frac{2A\tau^3}{3}-(A+B)\tau^2-2B\tau
  -4\tau+C,\\
  (1+\tau)\phi&=&-\frac{A\tau^4}{6}-\frac{(A+B)\tau^3}{3}-
  B\tau^2-2\tau^2+C\tau+D,\\
\end{eqnarray*}  
for some constants $C,D$. 
The resulting function defines a metric if it is positive on $(0,m)$. 

\begin{figure}[htbp]
    \label{fig:smooth}
    \begin{center}
      \input{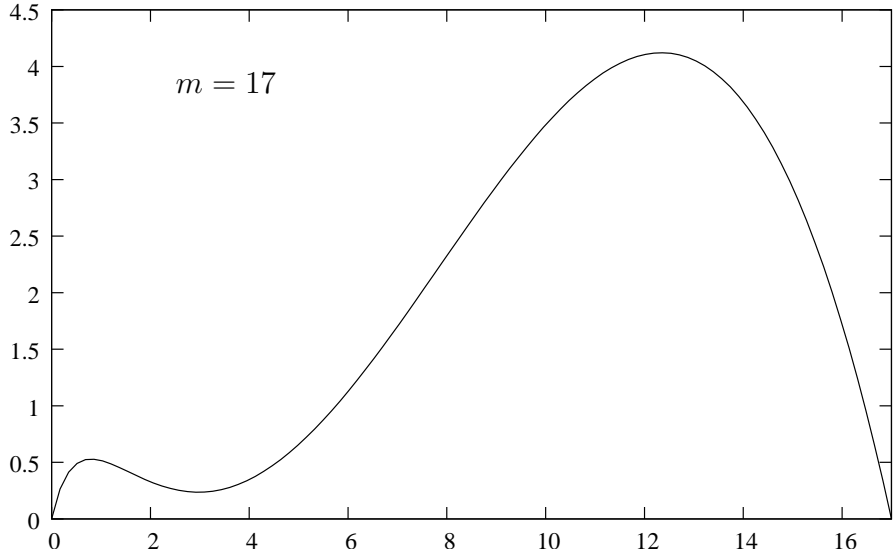}
      \caption{Momentum profile of an extremal metric on $X$, where
      $m=17$}
    \end{center}
  \end{figure}

Let us start with the case of a smooth metric on $X$. The boundary
conditions are
$\phi(0)=\phi(m)=0,\,\phi^\prime(0)=2,\,\phi^\prime(m)=-2$.
Solving the resulting system of linear equations for $A,B,C,D$, we obtain
\[
\begin{split}
  \phi(\tau) =
  \frac{2\tau(m-\tau)}{m(m^2+6m+6)(1+\tau)}\big(&\tau^2(2m+2)+\tau(-m^2+4m+6)+\\ 
  &+m^2+6m+6\big).
\end{split}
\]
This will be positive on $(0,m)$, if and only if the quadratic
expression (in $\tau$) in
brackets is positive on this interval. This is the case for $m<k_1$,
where $k_1$ is the only positive real root of the quartic
$m^4-16m^3-52m^2-48m-12$. Approximately $k_1\cong 18.889$. This is also
the result obtained by T\o{}nessen-Friedman~\cite{TF97}.
Figure~\ref{fig:smooth} shows a plot of $\phi$ for $m=17$.
The scalar curvature is
\[ S(\phi)(\tau)=
 \frac{24(m+1)}{m(m^2+6m+6)}\tau -\frac{6(3m^2+2m-2)}{m(m^2+6m+6)}.
\]

\begin{figure}[htbp]
    \label{fig:singSinf}
    \begin{center}
      \input{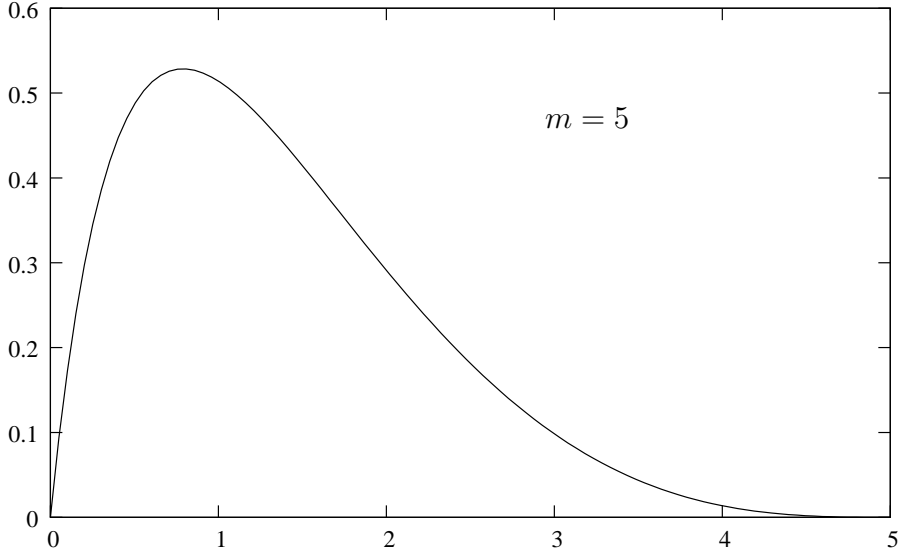}
      \caption{Momentum profile of an asymptotically hyperbolic extremal
      metric on $X\setminus S_\infty$, where
      $m=5$}
    \end{center}
  \end{figure}

We now move on to the case of complete metrics. 
To get a complete metric on $X\setminus S_\infty$ the boundary
conditions are $\phi(0)=\phi(m)=0,\,\phi^\prime(0)=2,\,\phi^\prime(m)=0$.
We get
\[ 
  \phi(\tau) = \frac{2\tau(m-\tau)^2}{m^2(m^2+6m+6)(1+\tau)}\big(
  \tau(-m^2+2m+3)+m^2+6m+6\big).
\]
This is positive on $(0,m)$ if the linear expression in brackets is
positive on this interval. 
This happens for $m\leq k_2$, where $k_2$ is the only
positive real root of the cubic 
$m^3-3m^2-9m-6$. Approximately $k_2\cong 5.0275$.
Figure~\ref{fig:singSinf} shows a plot of $\phi$ for $m=5$.
For $m<k_2$ the order
of vanishing of 
$\phi(\tau)$ at $\tau=m$ is precisely 2 so the metric is asymptotically
hyperbolic near $S_\infty$. 
The scalar curvature is
\begin{equation} \label{eq:scalS0}
  S(\phi)(\tau) = 
  \frac{12(m^2-2m-3)}{m^2(m^2+6m+6)}\tau-\frac{6(2m^2-m-4)}{m(m^2+6m+6)}. 
\end{equation}

For a complete metric on $X\setminus S_0$ the boundary conditions
are $\phi(0)=\phi(m)=0,\,\phi^\prime(0)=0,\,\phi^\prime(m)=-2$. We
obtain
\[ 
  \phi(\tau) =
  \frac{2\tau^2(m-\tau)}{m^2(m^2+6m+6)(1+\tau)}\big(\tau(2m^2+4m+3)
  -m^3+3m^2+9m+6\big).
\]
This is positive on $(0,m)$ if the linear term in brackets is positive
on this interval. 
This is the case for $m\leq k_2$, where $k_2$ is the same as above.
Again, for $m<k_2$ the order of vanishing of $\phi(\tau)$ at $\tau=0$ is
exactly 2, so the resulting metric is asymptotically hyperbolic near
$S_0$. 
The scalar curvature is
\begin{equation}
  S(\phi)(\tau) = 
  \frac{12(2m^2+4m+3)}{m^2(m^2+6m+6)}\tau-\frac{6(3m^2+5m+2)}{m(m^2+6m+6)}.
\end{equation}
Note that if $\phi$ is a solution with these boundary conditions, then
$\psi$ defined by
\[ \psi(\tau) = (a+1)\,\phi\left(\frac{\tau-a}{a+1}\right) \]
is a solution with the same boundary conditions on the interval
$[a,(a+1)m+a]$. This also gives a complete extremal metric on
$X\setminus S_0$, just in a different K\"ahler class. The corresponding
scalar curvature is given by
\begin{equation} \label{eq:scalSinfty}
 S(\psi)(\tau) = \frac{1}{a+1}\,S(\phi)\left(\frac{\tau-a}{a+1}
\right).
\end{equation}
We now summarise these results in a proposition which we
will use in the next section.
\begin{prop} \label{prop:extremalcomplete}
  There exists a complete extremal metric on $X\setminus S_\infty$ with
  momentum profile on $[0,m]$ for $m\leq k_2$ and scalar curvature given
  by~\ref{eq:scalS0}. When $m<k_2$ the resulting metric is
  asymptotically hyperbolic. There exists a complete extremal metric on
  $X\setminus S_0$ with momentum profile on $[c,m]$ for any
  positive $c$ and $m>c$ such that
  \[ \frac{m-c}{c+1} \leq k_2, \]
  with scalar curvature given by~\ref{eq:scalSinfty}. If the
  inequality is sharp, the resulting metric is asymptotically hyperbolic.
\end{prop}

Note that the polarisations for which we have not obtained
asymptotically hyperbolic extremal metrics are K-unstable 
according to the calculations in
Section~\ref{sec:ruledexample}. 

\section{The infimum of the Calabi functional}\label{sec:extremallwr}
In this section we compute the infimum of the Calabi functional on
the ruled surface considered above, for the unstable polarisations.
According to Theorem~\ref{thm:lowerCalabi}, a test-configuration $\chi$
which has negative Futaki invariant $F(\chi)<0$ gives a lower bound
\[ \Vert S(\omega)-\hat{S}\Vert_{L^2}\geq 4\pi
\frac{-F(\chi)}{\Vert\chi\Vert}
\]
on the Calabi functional. Donaldson conjectured in~\cite{Don05} that the
supremum of this lower bound over all test-configurations gives the
infimum of the Calabi functional. We will show that for our ruled
surface this is indeed the case. 
\begin{thm} For the ruled surface $X$ we have 
  \[ \inf_\omega \Vert S(\omega) - \hat{S}\Vert_{L^2} = 4\pi
  \sup_\chi
  \frac{-F(\chi)}{\Vert\chi\Vert}, \]
  where $\chi$ runs over all test-configuration for $X$. 
\end{thm}
\begin{proof}
Let the polarisation of $X$ be $L=C+mS_\infty$ working on the momentum
interval $[0,m]$ as before, and $m\geq k_1$ so
that $(X,L)$ is relatively K-unstable. 
We will define a sequence of momentum profiles $\phi_i$ and a
sequence of test-configurations $\chi_i$ corresponding to rational 
piecewise linear convex functions $h_i$ such that
\begin{equation}\label{eq:limiteq}
  \lim_{i\to\infty}\Vert S(\omega_{\phi_i})-\hat{S}\Vert_{L^2} =
  \lim_{i\to\infty} 4\pi \frac{-F(\chi_i)}{\Vert\chi_i\Vert}.
\end{equation}
It will follow that this limit is the infimum of the Calabi functional.
We will define these sequences by writing down their limits. In the case
of $\phi_i$ this means a $C^2$ momentum profile $\phi$ which may be zero on a
subset of $(0,m)$ and in the case of $\chi_i$ it means a
continuous convex function $h$
which is not necessarily rational and piecewise linear. 

There are two cases to consider. The first is when $k_1\leq m\leq
k_2(k_2+2)\cong 35.33$ 
(recall the constants $k_1,k_2$ defined in the previous section). In
this case we define a constant $c=\sqrt{m+1}-1$, and define
\begin{equation*}
  \phi(\tau) = \begin{cases} \psi_1(\tau)& \mbox{ for }\tau\in[0,c],\\
    \psi_2(\tau)&\mbox{ for }\tau\in[c,m].
  \end{cases}
\end{equation*}
where $\psi_1,\psi_2$ are the momentum profiles of the complete extremal
metrics given by
Proposition~\ref{prop:extremalcomplete}, so that $\phi$ satisfies 
\[ \phi(0)=0, \phi^\prime(0)=2,\quad \phi(c)=\phi^\prime(c)=0,\quad
\phi(m)=0, \phi^\prime(m)=-2.\]
Note that the assumption on $m$ and $c$ ensures
that the intervals $[0,c]$ and $[c,m]$ satisfy the conditions of
Proposition~\ref{prop:extremalcomplete} so that $\psi_1,\psi_2$ exist.
We can check explicitly using the formulae in the previous section, that
our choice of $c$ ensures that
$\psi_1^{\prime\prime}(c)=\psi_2^{\prime\prime}(c)$, so $\phi$ is in
$C^2$.  
The key point is that the scalar curvature
$S(\omega_\phi)$ is concave (by construction it is linear on the
intervals $[0,c]$ and $[c,m]$), which we can also
verify explicitly (here the assumption $m\geq k_1$ is used). 
This allows us to define $h(\tau)=\hat{S}-S(\omega_\phi)$ which is a
convex piecewise linear function. To verify
Equation~\ref{eq:limiteq} we apply Theorem~\ref{thm:ruledfutaki} which
implies that the normalised Futaki invariant is
\[ \frac{\frac{1}{2(2\pi)^2}\int_X
(S(\omega_\phi)-\hat{S})^2\frac{\omega_\phi^n}{n!}}{\frac{1}{2\pi}
\Vert S(\omega_\phi)
- \hat{S}\Vert_{L^2}} = \frac{1}{4\pi}\Vert S(\omega_\phi)-\hat{S}\Vert_{L^2}.\]

The second case is when $m > k_2(k_2+2)$. We now define a constant
\[ c=\frac{m+1}{k_2+1}-1, \]
and we define
\[ \phi(\tau) = \begin{cases} 
  \psi_1(\tau),\quad &\tau\in [0,k_2], \\
  0, &\tau\in [k_2,c], \\
  \psi_2(\tau),\quad &\tau\in [c,m],
\end{cases} \]
where again, $\psi_1$ and $\psi_2$ are given by
proposition~\ref{prop:extremalcomplete}. Once again our choice of the
intervals guarantees that $\phi$ is non-negative and in $C^2$, and also
the scalar curvature $S(\omega_\phi)$ is concave. On the interval
$[k_2,c]$ it is given by $-\frac{2}{1+\tau}$, so it is not piecewise
linear. We choose a sequence $h_i$ of piecewise linear convex
functions converging to $h=\hat{S}-S(\omega_\phi)$.
Theorem~\ref{thm:ruledfutaki} implies Equation~\ref{eq:limiteq} once
again. 
\end{proof}

In the proof we have exhibited the ``worst destabilising
test-configuration'' for each unstable polarisation which breaks up the
manifold into pieces. This is analogous
to the Harder-Narasimhan filtration of an unstable vector bundle. 
Note that most of
the time (except for when $\sqrt{m+1}$ is rational in the first case)
these are not algebraic test configurations. The first case is
deformation to the normal cone with a real parameter, and the unstable
manifold is split into two pieces which admit (complete) extremal
metrics. The second case is more complicated, and the middle piece does
not admit an extremal metric. Instead its cylindrical fibres become
infinitely long and thin.

\newpage
\bibliographystyle{amsplain} 
\addcontentsline{toc}{chapter}{Bibliography}
\bibliography{../mybib}

\providecommand{\bysame}{\leavevmode\hbox to3em{\hrulefill}\thinspace}
\providecommand{\MR}{\relax\ifhmode\unskip\space\fi MR }
\providecommand{\MRhref}[2]{%
  \href{http://www.ams.org/mathscinet-getitem?mr=#1}{#2}
}
\providecommand{\href}[2]{#2}
\begin{thebibliography}{10}

\bibitem{ACGT3}
V.~Apostolov, D.~M.~J. Calderbank, P.~Gauduchon, and C.~W.
  T\o{}nnesen-Friedman, \emph{Hamiltonian 2-forms in {K}\"ahler geometry {III},
  extremal metrics and stability}, preprint (2005).

\bibitem{Ati82}
M.~F. Atiyah, \emph{Convexity and commuting {H}amiltonians}, Bull. London Math.
  Soc. \textbf{14} (1982), no.~1, 1--15.

\bibitem{Aub78}
T.~Aubin, \emph{{\'E}quations du type {M}onge-{A}mp\`ere sur les variet\'es
  k\"ahl\'eriennes compactes}, Bull. Sci. Math. (2) \textbf{102} (1978), no.~1,
  63--95.

\bibitem{Cal82}
E.~Calabi, \emph{Extremal {K}\"ahler metrics}, Seminar on Differential Geometry
  (S.~T. Yau, ed.), Princeton, 1982.

\bibitem{Cal85}
\bysame, \emph{Extremal {K}\"ahler metrics {II}}, Differential geometry and
  complex analysis, Springer, 1985, pp.~95--114.

\bibitem{CFM}
P.~Cartier, J.~Fell, and P.~Meyer, \emph{Comparaison des mesures port\'ees par
  un ensemble convexe compact}, Bull. Soc. Math. France \textbf{92} (1964),
  435--445.

\bibitem{CT05}
X.~Chen and G.~Tian, \emph{Uniqueness of extremal {K}\"ahler metrics}, C. R.
  Math. Acad. Sci. Paris \textbf{340} (2005), no.~4, 287--290.

\bibitem{CY80}
S.~Y. Cheng and S.~T. Yau, \emph{On the existence of a complete {K}\"ahler
  metric on noncompact complex manifolds and the regularity of {F}efferman's
  equation}, Comm. Pure Appl. Math. \textbf{33} (1980), no.~4, 507--544.

\bibitem{Don97}
S.~K. Donaldson, \emph{Remarks on gauge theory, complex geometry and
  four-manifold topology}, Fields Medallists' Lectures (Atiyah and Iagolnitzer,
  eds.), World Scientific, 1997, pp.~384--403.

\bibitem{Don02}
\bysame, \emph{Scalar curvature and stability of toric varieties}, J.
  Differential Geom. \textbf{62} (2002), 289--349.

\bibitem{Don05_1}
\bysame, \emph{Interior estimates for solutions of {A}breu's equation},
  preprint (2005).

\bibitem{Don05}
\bysame, \emph{Lower bounds on the {C}alabi functional}, preprint (2005).

\bibitem{Fuj92}
A.~Fujiki, \emph{Moduli space of polarized algebraic manifolds and {K}\"ahler
  metrics}, Sugaku Expositions \textbf{5} (1992), no.~2, 173--191.

\bibitem{Fut83}
A.~Futaki, \emph{An obstruction to the existence of {E}instein-{K}\"ahler
  metrics}, Invent. Math. \textbf{73} (1983), 437--443.

\bibitem{FM95}
A.~Futaki and T.~Mabuchi, \emph{Bilinear forms and extremal {K}\"ahler vector
  fields associated with {K}\"ahler classes}, Math. Ann. \textbf{301} (1995),
  199--210.

\bibitem{GS06}
V.~Guillemin and S.~Sternberg, \emph{Riemann sums over polytopes}, preprint
  (2006).

\bibitem{HHM04}
T.~Hausel, E.~Hunsicker, and R.~Mazzeo, \emph{Hodge cohomology of gravitational
  instantons}, Duke Math. J. \textbf{122} (2004), no.~3, 485--548.

\bibitem{HS02}
A.~Hwang and A.~M. Singer, \emph{A momentum construction for circle-invariant
  {K}\"ahler metrics}, Trans. Amer. Math. Soc. \textbf{354} (2002), no.~6,
  2285--2325.

\bibitem{Kir84}
F.~C. Kirwan, \emph{Cohomology of quotients in symplectic and algebraic
  geometry}, Princeton University Press, 1984.

\bibitem{Kob95}
S.~Kobayashi, \emph{Transformation groups in differential geometry},
  Springer-Verlag, Berlin, 1995.

\bibitem{Mab86}
T.~Mabuchi, \emph{K-energy maps integrating {F}utaki invariants}, Tohoku Math.
  J. \textbf{38} (1986), no.~4, 575--593.

\bibitem{Mab04_1}
\bysame, \emph{Stability of extremal {K}\"ahler manifolds}, Osaka J. Math.
  \textbf{41} (2004).

\bibitem{Mab04_2}
\bysame, \emph{Uniqueness of extremal {K}\"ahler metrics for an integral
  {K}\"ahler class}, Internat. J. Math. \textbf{15} (2004), no.~6, 531--546.

\bibitem{Mal45}
A.~Malcev, \emph{On the theory of the {L}ie groups in the large}, Rec. Math.
  [Mat. Sbornik] N. S. \textbf{16(58)} (1945), 163--190.

\bibitem{Mar71}
M.~Maruyama, \emph{On automorphism groups of ruled surfaces}, J. Math. Kyoto
  University \textbf{11-1} (1971), 89--112.

\bibitem{MFK94}
D.~Mumford, J.~Fogarty, and F.~Kirwan, \emph{Geometric invariant theory},
  Springer-Verlag, 1994.

\bibitem{PSSW}
D.~H. Phong, J.~Song, J.~Sturm, and B.~Weinkove, \emph{The {M}oser-{T}rudinger
  inequality on {K}\"ahler-{E}instein manifolds}, preprint (2006).

\bibitem{PS06}
D.~H. Phong and J.~Sturm, \emph{Test configurations for {K}-stability and
  geodesic rays}, preprint (2006).

\bibitem{RT04}
J.~Ross and R.~P. Thomas, \emph{A study of the {H}ilbert-{M}umford criterion
  for the stability of projective varieties}, preprint (2004).

\bibitem{RT06}
\bysame, \emph{An obstruction to the existence of constant scalar curvature
  {K}\"ahler metrics}, J. Diff. Geom. \textbf{72} (2006), 429--466.

\bibitem{GSz04}
G.~Sz\'ekelyhidi, \emph{Extremal metrics and ${K}$-stability}, preprint (2004).

\bibitem{Thomas06}
R.~P. Thomas, \emph{Notes on {GIT} and symplectic reduction for bundles and
  varieties}, preprint (2006).

\bibitem{Tian97}
G.~Tian, \emph{K\"ahler-{E}instein metrics with positive scalar curvature},
  Invent. math \textbf{137} (1997), 1--37.

\bibitem{Tian00}
\bysame, \emph{Canonical metrics in {K}\"ahler geometry}, Lectures in
  {M}athematics {ETH} {Z}\"urich, Birkh\"auser Verlag, Basel, 2000.

\bibitem{TY87}
G.~Tian and S.-T. Yau, \emph{Existence of {K}\"ahler-{E}instein metrics on
  complete {K}\"ahler manifolds and their applications to algebraic geometry},
  Mathematical aspects of string theory (San Diego, Calif., 1986), Adv. Ser.
  Math. Phys., vol.~1, World Sci. Publishing, Singapore, 1987, pp.~574--628.

\bibitem{TF97}
C.~T\o{}nnesen-Friedman, \emph{Extremal {K}\"ahler metrics on ruled surfaces},
  Ph.D. thesis, Odense University, 1997.

\bibitem{Yau78}
S.-T. Yau, \emph{On the {R}icci curvature of a compact {K}\"ahler manifold and
  the complex {M}onge-{A}mp\`ere equation, {I}.}, Comment. Pure Appl. Math.
  \textbf{31} (1978), 339--411.

\bibitem{ZZ06}
B.~Zhou and X.~Zhu, \emph{Relative {K}-stability and modified {K}-energy on
  toric manifolds}, preprint (2006).

\end{thebibliography}

\end{document}